\documentclass[a4paper]{article}

\usepackage[cp1250]{inputenc}
\usepackage[IL2]{fontenc}

\usepackage{a4wide}
\usepackage[english]{babel}
\usepackage{euscript}
\usepackage{amstext,amsbsy,amscd,amssymb}
\usepackage{amsmath,enumitem}
\usepackage{amsfonts}
\usepackage{graphics}
\usepackage{graphicx}
\usepackage{microtype}
\usepackage{xcolor}
\usepackage{dynkin-diagrams}
\usepackage{tikz,caption,subcaption}
\usetikzlibrary{arrows,decorations.pathmorphing,backgrounds,positioning,fit,petri}
\include{diagxy}
\allowdisplaybreaks

\let\rarr=\rightarrow

\let\veps=\varepsilon
\let\mcal=\mathcal
\let\mfrak=\mathfrak
\let\eus=\EuScript
\let\bra=\langle
\let\ket=\rangle

\def\N{\mathbb{N}}
\def\Z{\mathbb{Z}}
\def\R{\mathbb{R}}
\def\C{\mathbb{C}}

\def\Q{\mathbb{Q}}
\def\vac{|0\rangle}

\def\normOrd #1{{\mathop{:}\nolimits\!#1\!\mathop{:}\nolimits}}

\DeclareMathSymbol{\squares}{\mathord}{AMSa}{"03}

\DeclareMathOperator{\Ann}{Ann}
\DeclareMathOperator{\Fun}{Fun}
\DeclareMathOperator{\Vect}{Vect}
\DeclareMathOperator{\Pol}{Pol}

\DeclareMathOperator{\Specm}{Specm}

\def\End{\mathop {\rm End} \nolimits}

\def\Hom{\mathop {\rm Hom} \nolimits}

\def\ad{\mathop {\rm ad} \nolimits}
\def\gr{\mathop {\rm gr} \nolimits}
\def\im{\mathop {\rm im} \nolimits}

\def\GL{\mathop {\rm GL} \nolimits}

\def\id{{\rm id}}

\def\rank{\mathop {\rm rank} \nolimits}
\def\Ad{\mathop {\rm Ad} \nolimits}

\def\tr{\mathop {\rm tr} \nolimits}
\def\Res{\mathop {\rm Res} \nolimits}
\def\Ind{\mathop {\rm Ind} \nolimits}

\def\htt{\mathop {\rm ht} \nolimits}

\long\def\proof #1{\noindent \emph{Proof.}\ #1 \hfill $\squares$

\medskip}

\newcounter{num}[section]
\numberwithin{equation}{section}
\numberwithin{num}{section}

\long\def\definition #1 {\refstepcounter{num} \noindent {\bf Definition \thenum.} #1

\medskip}

\long\def\theorem #1{\refstepcounter{num} \noindent \ifnum\value{section}>0 {\bf Theorem \thenum.} #1 \else {\bf Theorem \Alph{num}.} #1 \fi

\medskip}

\long\def\lemma #1{\refstepcounter{num}  \noindent {\bf Lemma \thenum.} #1

\medskip}

\long\def\proposition #1{\refstepcounter{num}  \noindent {\bf Proposition \thenum.} #1

\medskip}

\long\def\corollary #1{\refstepcounter{num}  \noindent {\bf Corollary \thenum.} #1

\medskip}

\long\def\remark #1{\noindent {\bf Remark.}\ #1}

\makeatletter
\newcommand*\if@single[3]{%
  \setbox0\hbox{${\mathaccent"0362{#1}}^H$}%
  \setbox2\hbox{${\mathaccent"0362{\kern0pt#1}}^H$}%
  \ifdim\ht0=\ht2 #3\else #2\fi
  }
\newcommand*\rel@kern[1]{\kern#1\dimexpr\macc@kerna}
\newcommand*\widebar[1]{\@ifnextchar^{{\wide@bar{#1}{0}}}{\wide@bar{#1}{1}}}
\newcommand*\wide@bar[2]{\if@single{#1}{\wide@bar@{#1}{#2}{1}}{\wide@bar@{#1}{#2}{2}}}
\newcommand*\wide@bar@[3]{%
  \begingroup
  \def\mathaccent##1##2{%
    \if#32 \let\macc@nucleus\first@char \fi
    \setbox\z@\hbox{$\macc@style{\macc@nucleus}_{}$}%
    \setbox\tw@\hbox{$\macc@style{\macc@nucleus}{}_{}$}%
    \dimen@\wd\tw@
    \advance\dimen@-\wd\z@
    \divide\dimen@ 3
    \@tempdima\wd\tw@
    \advance\@tempdima-\scriptspace
    \divide\@tempdima 10
    \advance\dimen@-\@tempdima
    \ifdim\dimen@>\z@ \dimen@0pt\fi
    \rel@kern{0.6}\kern-\dimen@
    \if#31
      \overline{\rel@kern{-0.6}\kern\dimen@\macc@nucleus\rel@kern{0.4}\kern\dimen@}%
      \advance\dimen@0.4\dimexpr\macc@kerna
      \let\final@kern#2%
      \ifdim\dimen@<\z@ \let\final@kern1\fi
      \if\final@kern1 \kern-\dimen@\fi
    \else
      \overline{\rel@kern{-0.6}\kern\dimen@#1}%
    \fi
  }%
  \macc@depth\@ne
  \let\math@bgroup\@empty \let\math@egroup\macc@set@skewchar
  \mathsurround\z@ \frozen@everymath{\mathgroup\macc@group\relax}%
  \macc@set@skewchar\relax
  \let\mathaccentV\macc@nested@a
  \if#31
    \macc@nested@a\relax111{#1}%
  \else
    \def\gobble@till@marker##1\endmarker{}%
    \futurelet\first@char\gobble@till@marker#1\endmarker
    \ifcat\noexpand\first@char A\else
      \def\first@char{}%
    \fi
    \macc@nested@a\relax111{\first@char}%
  \fi
  \endgroup
}
\makeatother

\newcommand\rsmraise[1]{%
  \ifx#1\displaystyle .8\else
    \ifx#1\textstyle .8\else
      \ifx#1\scriptstyle .6\else
        .45%
      \fi
    \fi
  \fi}

\title{Positive energy representations of affine vertex algebras}

\author{Vyacheslav Futorny, Libor K\v{r}i\v{z}ka}

\AtEndDocument{\bigskip{\footnotesize%
  (V.\,Futorny) \textsc{Instituto de Matem\'atica e Estat\'istica, Universidade de S\~ao Paulo,  S\~ao Paulo, Brazil, and  International Center for Mathematics, SUSTech, Shenzhen, China} \par
  \textit{E-mail address}: \texttt{futorny@ime.usp.br} \par
  \addvspace{\medskipamount}
  (L.\,K\v{r}i\v{z}ka) \textsc{Instituto de Matem\'atica e Estat\'stica, Universidade de Sa\~{o} Paulo, Sa\~{o} Paulo,  Brasil} \par
  \textit{E-mail address}: \texttt{krizka.libor@gmail.com} \par
}}

\date{}

\begin{document}

\maketitle

\begin{abstract}
We construct new families of positive energy representations of affine vertex algebras together with their free field realizations by using localization technique. We introduce the \emph{twisting functor} $T_\alpha$ on the category of modules over \emph{affine Kac--Moody algebra} $\widehat{\mfrak{g}}_\kappa$ of level $\kappa$ for any positive root $\alpha$ of $\mfrak{g}$, and the \emph{Wakimoto functor} from a certain category of $\mfrak{g}$-modules to the category of \emph{smooth} $\widehat{\mfrak{g}}_\kappa$-modules. These two functors commute (taking a proper restriction of $T_\alpha$ on $\mfrak g$-modules) and the image of the Wakimoto functor consists of \emph{relaxed Wakimoto} $\widehat{\mfrak{g}}_\kappa$-modules. In particular, applying the twisting functor $T_\alpha$ to the relaxed Wakimoto $\widehat{\mfrak{g}}_\kappa$-module whose top degree component is isomorphic to the Verma $\mfrak{g}$-module $M^\mfrak{g}_\mfrak{b}(\lambda)$, we obtain the relaxed Wakimoto $\widehat{\mfrak{g}}_\kappa$-module whose top degree component is isomorphic to the $\alpha$-Gelfand--Tsetlin $\mfrak{g}$-module $W^\mfrak{g}_\mfrak{b}(\lambda, \alpha)$. We show that the relaxed Verma module and relaxed Wakimoto module whose top degree components are such $\alpha$-Gelfand--Tsetlin modules, are isomorphic generically. This is an analogue of the result of E.\,Frenkel for Wakimoto modules both for critical and non-critical level. For a parabolic subalgebra $\mfrak{p}$ of $\mfrak{g}$ we construct a large family of admissible $\mfrak{g}$-modules as images under the twisting functor of generalized Verma modules induced from $\mfrak{p}$. In this way, we obtain new simple positive energy representations of simple affine vertex algebras.
\medskip

\medskip


\end{abstract}

\thispagestyle{empty}

\tableofcontents

%%%%%%%%%%%%%%%%%%%%%%%%%%%%%%%%%%%%%%%%%%%%%%%%%%%%%%%%%%%%%%%%%%%%%%%%%%%%%%%%%%%%%%%%%%%
%%%%%%%%%%%%%%%%%%%%%%%%%%%%%%%%%%%%%%%%%%%%%%%%%%%%%%%%%%%%%%%%%%%%%%%%%%%%%%%%%%%%%%%%%%%

\section*{Introduction}
\addcontentsline{toc}{section}{Introduction}

Classical Wakimoto modules provide a free field realization for affine Kac--Moody algebras. They were constructed in \cite{Wakimoto1986}, \cite{Feigin-Frenkel1988}, \cite{deBoer-Feher1997}, \cite{Szczesny2002} among the others. Simple Wakimoto modules are isomorphic to Verma modules (for the \emph{standard Borel subalgebra}) with the same highest weight. Free field realizations for affine Kac--Moody algebras associated with nonstandard Borel subalgebras were studied in \cite{Cox2005} (\emph{Imaginary Verma modules}), \cite{Cox-Futorny2004} (\emph{intermediate Wakimoto modules}),
\cite{Futorny-Krizka-Somberg2019} (\emph{generalized Imaginary Verma modules}). Parabolic versions of Wakimoto modules,
\emph{generalized Wakimoto modules}, were introduced in  \cite{Frenkel2005}.

Wakimoto modules are the key objects in the theory of rational affine vertex algebras. However, there is a growing interest in the study  of non-highest weight representations of $\widehat{\mfrak{g}}$ that might play a significant role for non-rational affine vertex algebras. For instance, \emph{relaxed Verma modules} for simple affine $\widehat{\mfrak{sl}}_2$ vertex algebra were used to study the $N=2$ conformal field theory in \cite{Feigin-Semikhatov-Tipunin1998} and the $N=4$ conformal field theory in \cite{Adamovic2016}. For recent developments we refer to \cite{Arakawa-Futorny-Ramirez2017}, \cite{Auger-Creutzig-Ridout2018}, \cite{Kawasetsu-Ridout2019}, \cite{Kawasetsu-Ridout2019b}.

The goal of the current paper is to give an explicit free field construction of new families of positive energy representations of the universal affine vertex algebras and simple affine vertex algebras using two main technical tools: the \emph{twisting functor} on the category of modules over the (untwisted) affine Kac--Moody algebra $\widehat{\mfrak{g}}_\kappa$ of level $\kappa$ assigned to an arbitrary positive root of a semisimple Lie algebra $\mfrak{g}$, and the \emph{Wakimoto functor} from a certain subcategory of $\mfrak{g}$-modules to the category of positive energy $\widehat{\mfrak{g}}_\kappa$-modules. In particular, the Wakimoto functor applied to Verma $\mfrak{g}$-modules gives Verma $\widehat{\mfrak{g}}_\kappa$-modules generically.

In the finite-dimensional case the twisting functor was defined in \cite{Futorny-Krizka2019b} following the work of Deodhar \cite{Deodhar1980}. If $\alpha$ is a simple root, then the twisting functor $T_\alpha$ is related to the Arkhipov's twisting functor \cite{Arkhipov2004} on the category $\mcal{O}(\mfrak{g})$. By applying the twisting functor $T_\alpha$ to the generalized Verma $\mfrak{g}$-module $M^\mfrak{g}_\mfrak{p}(\lambda)$ induced from the simple finite-dimensional $\mfrak{p}$-module with highest weight $\lambda$ one obtains the $\alpha$-Gelfand--Tsetlin $\mfrak{g}$-module $W^\mfrak{g}_\mfrak{p}(\lambda,\alpha)$ with finite $\Gamma_\alpha$-multiplicities, where $\Gamma_\alpha$ is the commutative subalgebra of $U(\mfrak{g})$ generated by the Cartan subalgebra $\mfrak{h}$ and by the center of $U(\mfrak{s}_\alpha)$, where $\mfrak{s}_\alpha$ is the Lie subalgebra of $\mfrak{g}$ given by the $\mfrak{sl}_2$-triple for the root $\alpha$. The $\mfrak{g}$-modules $W^\mfrak{g}_\mfrak{b}(\lambda, \alpha)$ are cyclic weight modules with respect to the Cartan subalgebra $\mfrak{h}$ with infinite-dimensional weight subspaces if $\alpha$ is not a simple root.  On the other hand, if $\alpha$ is a simple root, then these modules are twisted Verma $\mfrak{g}$-modules up to conjugation of the action of $\mfrak{g}$, see \cite{Andersen-Lauritzen2003}, \cite{Krizka-Somberg2015b}, \cite{Musson2019}.

Modifying the construction given in \cite{Futorny-Krizka2019} and \cite{Futorny-Krizka2019b} in the finite-dimensional setting, we define the twisting functor $T_\alpha$ assigned to a positive root $\alpha \in \Delta_+$ of $\mfrak{g}$ and obtain an endofunctor on the category of $\widehat{\mfrak{g}}_\kappa$-modules. The most important properties of the twisting functor $T_\alpha$ are collected in the theorem below (cf.\ Theorems \ref{thm:twisting functor cat E}, \ref{thm:highest weight modules condition}, \ref{thm-And}, \ref{thm:twisting functor intertwining}).
\medskip

\theorem{\label{thm-main1}
\begin{enumerate}[topsep=0pt,itemsep=0pt]
\item[i)] For $\alpha \in \smash{\widehat{\Delta}}^{\rm re}$, the twisting functor $T_\alpha$ preserves the category $\mcal{E}(\widehat{\mfrak{g}}_\kappa)$. Moreover, the twisting functor $T_\alpha$ preserves also the category $\mcal{E}_+\!(\widehat{\mfrak{g}}_\kappa)$ provided $\alpha \in \Delta$.
\item[ii)] Let $\alpha \in \smash{\widehat{\Delta}}^{\rm re}_+$ and let $M$ be a smooth weight $\widetilde{\mfrak{g}}_\kappa$-module on which the central element $c$ acts as the identity. Then the $\widetilde{\mfrak{g}}_\kappa$-module $T_\alpha(M)$ is a Gelfand--Tsetlin module with finite $\Gamma_\alpha$-multiplicities if and only if the first cohomology group $H^1(\mfrak{s}_\alpha^-;M)$ is a weight $\smash{\widetilde{\mfrak{h}}}$-module with finite-dimensional weight spaces.
\item[iii)] For $\alpha \in \Delta \subset \smash{\widehat{\Delta}}^{\rm re}$ there exists a natural isomorphism between $ T_\alpha \circ \mathbb{M}_{\kappa,\mfrak{g}}$ and $\mathbb{M}_{\kappa,\mfrak{g}} \circ\, T_\alpha^\mfrak{g}$, where $T_\alpha^\mfrak{g} \colon \mcal{M}(\mfrak{g}) \rarr \mcal{M}(\mfrak{g})$ is the twisting functor for $\mfrak{g}$ assigned to $\alpha$. In particular, we have
    \begin{align*}
        T_\alpha(\mathbb{M}_{\kappa,\mfrak{g}}(M^\mfrak{g}_\mfrak{p}(\lambda))) \simeq \mathbb{M}_{\kappa,\mfrak{g}}(W^\mfrak{g}_\mfrak{p}(\lambda,\alpha))
    \end{align*}
    for $\lambda \in \Lambda^+(\mfrak{p})$ and
    $$\alpha \in \Delta_+^\mfrak{u}=  \{\alpha \in \Delta_+;\, \mfrak{g}_\alpha \subset \mfrak{u}\},$$ where $\mfrak{p}$ is a standard parabolic subalgebra of $\mfrak{g}$ with the nilradical $\mfrak{u}$,  $\Lambda^+(\mfrak{p})$ is the set of $\mfrak{p}$-dominant integral weights. Moreover, the twisting functor $T_\alpha$ commutes with tensoring by Weyl modules.
\end{enumerate}}

Here, the commutative subalgebra $\Gamma_\alpha$ of $U_c(\widetilde{\mfrak{g}}_\kappa)$ is generated by the Cartan subalgebra $\smash{\widetilde{\mfrak{h}}}$ and by the center $Z(\mfrak{s}_\alpha)$ of $U(\mfrak{s}_\alpha)$, where $\mfrak{s}_\alpha$ is the Lie subalgebra of $\widetilde{\mfrak{g}}_\kappa$ given by the $\mfrak{sl}_2$-triple for the root $\alpha$. Besides, we denote $\mfrak{s}_\alpha^- = \mfrak{s}_\alpha \cap \smash{\widehat{\widebar{\mfrak{n}}}}_{\rm st}$, where $\widetilde{\mfrak{g}}_\kappa = \widehat{\mfrak{n}}_{\rm st} \oplus \smash{\widetilde{\mfrak{h}}} \oplus \smash{\widehat{\widebar{\mfrak{n}}}}_{\rm st}$ is the triangular decomposition of the extended affine Kac--Moody algebra $\widetilde{\mfrak{g}}_\kappa$. Further, $\mcal{E}(\widehat{\mfrak{g}}_\kappa)$ and $\mcal{E}_+\!(\widehat{\mfrak{g}}_\kappa)$ stand for the categories of smooth $\smash{\widehat{\mfrak{g}}}_\kappa$-modules and positive energy $\smash{\widehat{\mfrak{g}}}_\kappa$-modules, respectively, on which the central element $c$ acts as the identity.

In particular, by applying the twisting functor $T_\alpha$ to the relaxed Verma module $\mathbb{M}_{\kappa,\mfrak{g}}(M^\mfrak{g}_\mfrak{p}(\lambda))$ we obtain a new class of positive energy $\widehat{\mfrak{g}}_\kappa$-modules, \emph{relaxed Verma modules} induced from the  $\alpha$-Gelfand--Tsetlin $\mfrak{g}$-module $W^\mfrak{g}_\mfrak{p}(\lambda,\alpha)$.

The Feigin--Frenkel homomorphism between the universal affine vertex algebra $\mcal{V}_\kappa(\mfrak{g})$ and the tensor product of the Weyl vertex algebra $\mcal{M}_{\widebar{\mfrak{n}}}$ with the Heisenberg vertex algebra $\mcal{V}_{\kappa-\kappa_c}(\mfrak{h})$ gives an explicit free field construction of Wakimoto modules (\cite{Frenkel2005}). We use the Feigin--Frenkel homomorphism to get a free field realization of relaxed Verma modules in Theorem \ref{thm-FF-hom} and Theorem \ref{thm-realiz}. The obtained $\widehat{\mfrak{g}}_\kappa$-modules are \emph{relaxed Wakimoto modules}, they are images of Verma modules $M^\mfrak{g}_\mfrak{b}(\lambda)$ and $\alpha$-Gelfand--Tsetlin modules $W^\mfrak{g}_\mfrak{b}(\lambda, \alpha)$ under the Wakimoto functor $\mathbb{W}_{\kappa,\mfrak{g}}$. The most important properties of the Wakimoto functor are collected in the theorem below (cf.\ Theorem \ref{thm:twisting functor intertwining Wakimoto}, Corollary \ref{thm:alpha Verma-Wakimoto isomorphism}, Corollary \ref{thm:alpha Verma-Wakimoto isomorphism critical}).
\medskip

\theorem{\label{thm-main2}
Let  $\lambda \in \mfrak{h}^*$, $\alpha \in \Delta_+ \subset \smash{\widehat{\Delta}}^{\rm re}_+$  a positive root.
\begin{enumerate}[topsep=0pt,itemsep=0pt]
\item[i)] The functors $T_\alpha \circ \mathbb{W}_{\kappa,\mfrak{g}}$ and $\mathbb{W}_{\kappa,\mfrak{g}} \circ\, T_\alpha^\mfrak{g}$ are naturally isomorphic. In particular, we have
    \begin{align*}
       T_\alpha(\mathbb{W}_{\kappa,\mfrak{g}}(M^\mfrak{g}_\mfrak{b}(\lambda))) \simeq \mathbb{W}_{\kappa,\mfrak{g}}(W^\mfrak{g}_\mfrak{b}(\lambda,\alpha)).
    \end{align*}

\item[ii)] If $\kappa$ is a non-critical level, then
    \begin{align*}
        \mathbb{M}_{\kappa,\mfrak{g}}(W^\mfrak{g}_\mfrak{b}(\lambda,\alpha))\simeq \mathbb{W}_{\kappa,\mfrak{g}}(W^\mfrak{g}_\mfrak{b}(\lambda,\alpha)),
    \end{align*}
     provided the Verma module $\mathbb{M}_{\kappa,\mfrak{g}}(M^\mfrak{g}_\mfrak{b}(\lambda))$ is a simple $\widehat{\mfrak{g}}_\kappa$-module.
\item[iii)] For the critical level $\kappa_c$, we have
    \begin{align*}
        \mathbb{M}_{\kappa_c, \mfrak{g}}(W^\mfrak{g}_\mfrak{b}(\lambda,\alpha))\simeq \mathbb{W}_{\kappa_c,\mfrak{g}}(W^\mfrak{g}_\mfrak{b}(\lambda,\alpha))
    \end{align*}
    if $\lambda$ satisfies $\langle \lambda +\rho, \gamma^\vee \rangle \notin -\N$ for all $\gamma \in \Delta_+$, i.e.\ $\lambda$ is a dominant weight.
\end{enumerate}}

We see that generically the relaxed Verma module $\mathbb{M}_{\kappa,\mfrak{g}}(W^\mfrak{g}_\mfrak{b}(\lambda,\alpha))$ and
the relaxed Wakimoto module $\mathbb{W}_{\kappa,\mfrak{g}}(W^\mfrak{g}_\mfrak{b}(\lambda,\alpha))$
are isomorphic. The generic condition is given by the simplicity of the Verma module for $\widehat{\mfrak{g}}_\kappa$ with highest weight $\lambda$ and the non-criticality of the level. On the other hand, if the level is critical and $\lambda$ is dominant, then these modules are always isomorphic, which is an analogue of the corresponding result of Frenkel \cite{Frenkel2007-book} for Wakimoto modules.

Further on, in Section \ref{sec-pos-energy} we describe families of positive energy representations of the simple affine vertex algebra $\mcal{L}_\kappa(\mfrak{g})$ of an admissible level $\kappa=k\kappa_0$ associated to a simple Lie algebra $\mfrak{g}$, where $\kappa_0$ is the normalized $\mfrak{g}$-invariant bilinear form on $\mfrak{g}$ and $k \in \C$. Admissible highest weight $\mfrak{g}$-modules of level $k$ were classified in \cite{Arakawa2016}. We denote by ${\rm Pr}_{k}$  the set of admissible weights $\lambda \in \smash{\widehat{\mfrak{h}}^*}$ for which $\smash{\widehat{\Delta}(\lambda)} = \smash{y(\widehat{\Delta}(k\Lambda_0))}$ for some element $y$ from the extended affine Weyl group $\smash{\widetilde{W}}$ of $\mfrak{g}$. Then we set $\widebar{{\rm Pr}}_k = \{\widebar{\lambda};\, \lambda \in {\rm Pr}_k\}$ for the canonical projections to $\mfrak{h}^*$. Let $\mfrak{p}$ be a standard parabolic subalgebra of $\mfrak{g}$ and let $k \in \Q$ be an admissible number for $\mfrak{g}$. Let us denote by $\Omega_k(\mfrak{p})$ the set of  weights $\lambda \in \widebar{{\rm Pr}}_k \cap \Lambda^+(\mfrak{p})$ such that $\bra \lambda + \rho, \alpha^\vee \ket \notin \N$ for all $\alpha \in \Delta_+^\mfrak{u}$. For these weights the generalized Verma module $M^\mfrak{g}_\mfrak{p}(\lambda)$ is a simple $\mfrak{g}$-module, and we have the following result.

\medskip

\theorem{\label{thm-main3}
Let $\lambda \in \Omega_k(\mfrak{p})$ and $\alpha \in \Delta^\mfrak{u}_+$. Then the $\mfrak{g}$-module $W^\mfrak{g}_\mfrak{p}(\lambda,\alpha)$ is admissible of level $k$.}

The nilpotent cone of $\mfrak{g}$ is an irreducible closed algebraic subvariety of $\mfrak{g}$. It decomposes into the finite union of adjoint orbits. If $\mcal{O}$ is such an adjoint orbit, then $\mcal{O}^*$ denotes the corresponding coadjoint orbit of $\mfrak{g}^*$. A nilpotent orbit $\mcal{O}$ of $\mfrak{g}$ is the orbit of a simple $\mfrak{g}$-module $E$ if $\mcal{V}(\Ann_{U(\mfrak{g})}\!E)=\smash{\widebar{\mcal{O}^*}}$. Let $k \in \Q$ be an admissible number for $\mfrak{g}$. Nilpotent orbits of admissible $\mfrak{g}$-modules of level $k$ were described in \cite{Arakawa2015}. For a nilpotent orbit $\mcal{O}$ of $\mfrak{g}$, we define
\begin{align*}
  \widebar{{\rm Pr}}_k^\mcal{O} =\{\lambda\in \widebar{{\rm Pr}}_k;\, \mcal{V}(J_\lambda) = \widebar{\mcal{O}^*}\},
\end{align*}
where $J_\lambda$ is the annihilator of the simple $\mfrak{g}$-module with highest weight $\lambda \in \mfrak{h}^*$ and $\mcal{V}(J_\lambda)$ is the associated variety. Then for the standard Borel subalgebra $\mfrak{b}$ of $\mfrak{g}$ we have that $\Omega_k(\mfrak{b}) = \widebar{{\rm Pr}}_k^{\smash{\mcal{O}_{\rm reg}}}$, where $\mcal{O}_{\rm reg}$ is the regular nilpotent orbit (cf.\ Proposition \ref{prop-orbit-borel}). Moreover, in the case $\mfrak{g}=\mfrak{sl}_n$ we give a more convenient description of the set $\Omega_k(\mfrak{p})$ for any standard parabolic subalgebra $\mfrak{p}$ of $\mfrak{g}$ (cf.\ Theorem \ref{thm-omega-sln}).
\medskip

We denote by $\C$, $\R$, $\Z$, $\N_0$ and $\N$ the set of complex numbers, real numbers, integers, non-negative integers and positive integers, respectively. All algebras and modules are considered over the field of complex numbers.

%%%%%%%%%%%%%%%%%%%%%%%%%%%%%%%%%%%%%%%%%%%%%%%%%%%%%%%%%%%%%%%%%%%%%%%%%%%%%%%%%%%%%%%%%%
%%%%%%%%%%%%%%%%%%%%%%%%%%%%%%%%%%%%%%%%%%%%%%%%%%%%%%%%%%%%%%%%%%%%%%%%%%%%%%%%%%%%%%%%%%

\section{Preliminaries}
\label{sec:prelim}

Let $\mfrak{g}$ be a complex semisimple finite-dimensional Lie algebra and let $\mfrak{h}$ be a Cartan subalgebra of $\mfrak{g}$. We denote by $\Delta$ the root system of $\mfrak{g}$ with respect to $\mfrak{h}$, by $\Delta_+$ a positive root system in $\Delta$ and by $\Pi \subset \Delta_+$ the set of simple roots. For $\alpha \in \Delta_+$, let $h_\alpha \in \mfrak{h}$ be the corresponding coroot and let $e_\alpha$ and $f_\alpha$ be basis of root subspaces $\mfrak{g}_\alpha$ and $\mfrak{g}_{-\alpha}$, respectively, defined by the requirement $[e_\alpha, f_\alpha] = h_\alpha$. We also set
\begin{align*}
  Q = \sum_{\alpha \in \Pi} \Z \alpha \qquad \text{and} \qquad Q_+ = \sum_{\alpha \in \Pi} \N_0 \alpha
\end{align*}
together with
\begin{align*}
  P = \sum_{\alpha \in \Pi} \Z \omega_\alpha \qquad \text{and} \qquad P_+ = \sum_{\alpha \in \Pi} \N_0 \omega_\alpha,
\end{align*}
where $\omega_\alpha \in \mfrak{h}^*$ for $\alpha \in \Pi$ is the fundamental weight determined by $\omega_\alpha(h_\gamma) = \delta_{\alpha,\gamma}$ for all $\gamma \in \Pi$. We call $Q$ the root lattice and $P$ the weight lattice. Further, we define the Weyl vector $\rho \in \mfrak{h}^*$ by
\begin{align*}
 \rho =  {1 \over 2} \sum_{\alpha \in \Delta_+} \alpha.
\end{align*}
The standard Borel subalgebra $\mfrak{b}$ of $\mfrak{g}$ is defined through $\mfrak{b} = \mfrak{h} \oplus \mfrak{n}$ with the nilradical $\mfrak{n}$ and the opposite nilradical $\widebar{\mfrak{n}}$ given by
\begin{align*}
  \mfrak{n} = \bigoplus_{\alpha \in \Delta_+} \mfrak{g}_\alpha \qquad \text{and} \qquad \widebar{\mfrak{n}} = \bigoplus_{\alpha \in \Delta_+} \mfrak{g}_{-\alpha}.
\end{align*}
Besides, we have the corresponding triangular decomposition
\begin{align*}
  \mfrak{g} = \widebar{\mfrak{n}} \oplus \mfrak{h} \oplus \mfrak{n}
\end{align*}
of the Lie algebra $\mfrak{g}$.

Let $\kappa_\mfrak{g}$ be the Cartan--Killing form on $\mfrak{g}$ and $(\cdot\,,\cdot)_\mfrak{g}$ the corresponding induced bilinear form on $\mfrak{g}^*$. Whenever $\alpha \in \mfrak{h}^*$ satisfies $(\alpha,\alpha)_\mfrak{g} \neq 0$, we define $s_\alpha \in \GL(\mfrak{h}^*)$ by
\begin{align*}
  s_\alpha(\gamma) = \gamma - {2(\alpha,\gamma)_\mfrak{g} \over (\alpha,\alpha)_\mfrak{g}}\, \alpha
\end{align*}
for $\gamma \in \mfrak{h}^*$. The subgroup $W$ of $\GL(\mfrak{h}^*)$ given by
\begin{align*}
  W = \langle s_\alpha;\, \alpha \in \Pi \rangle
\end{align*}
is called the Weyl group of $\mfrak{g}$. Let us note that $W$ is a finite Coxeter group.

Moreover, if $\mfrak{g}$ is a simple Lie algebra, we denote by $\kappa_0$ the $\mfrak{g}$-invariant symmetric bilinear form on $\mfrak{g}$ normalized in such a way that $(\theta,\theta)=2$, where $\theta \in \Delta_+$ is the highest root of $\mfrak{g}$ (by definition the highest weight of the adjoint representation $\mfrak{g}$) and $(\cdot\,,\cdot)$ is the corresponding induced bilinear form on $\mfrak{g}^*$. Further, we denote by $\theta_s$ the highest short root of $\mfrak{g}$. Then we have $(\theta_s,\theta_s)= 2/r^\vee$, where $r^\vee$ is the lacing number of $\mfrak{g}$, i.e.\ the maximal number of edges in the Dynkin diagram of $\mfrak{g}$. We also define
\begin{align*}
  P^\vee = \bigoplus_{\alpha \in \Pi} \Z\omega_\alpha^\vee \qquad  \text{and} \qquad P^\vee_+ = \bigoplus_{\alpha \in \Pi} \N_0\omega_\alpha^\vee,
\end{align*}
where $\omega_\alpha^\vee \in \mfrak{h}^*$ for $\alpha \in \Pi$ is the fundamental coweight defined by $(\omega_\alpha^\vee, \gamma) = \delta_{\alpha,\gamma}$ for all $\gamma \in \Pi$. We call $P^\vee$ the coweight lattice.
\medskip

The category of all $\mfrak{g}$-modules we denote by $\mcal{M}(\mfrak{g})$. We say that a $\mfrak{g}$-module $M$ is a \emph{generalized weight} (with respect to $\mfrak{h}$) $\mfrak{g}$-module, if the action of $\mfrak{h}$ on $M$ is locally finite. If the action of $\mfrak{h}$ is semisimple on $M$, then $M$ is called a \emph{weight} (with respect to $\mfrak{h}$) $\mfrak{g}$-module. In particular, any simple generalized weight $\mfrak{g}$-module is a weight $\mfrak{g}$-module. Further, for a Lie subalgebra $\mfrak{a}$ of $\mfrak{g}$ we denote by $\mcal{I}(\mfrak{g},\mfrak{a})$ and $\mcal{I}_f(\mfrak{g},\mfrak{a})$ the full subcategories of $\mcal{M}(\mfrak{g})$ consisting of locally $\mfrak{a}$-finite weight $\mfrak{g}$-modules and finitely generated locally $\mfrak{a}$-finite weight $\mfrak{g}$-modules, respectively.
\medskip

For a commutative algebra $\Gamma$ we denote by $\Hom(\Gamma,\C)$ the set of all characters of $\Gamma$, i.e.\ algebra homomorphisms from $\Gamma$ to $\C$. Let $M$ be a $\Gamma$-module. For each character $\chi \in \Hom(\Gamma,\C)$ we set
\begin{align}
  M_\chi = \{v \in M;\, (\exists k \in \N)\,(\forall a \in \Gamma)\, (a-\chi(a))^kv=0\}.
\end{align}
When $M_\chi \neq\{0\}$, we say that $\chi$ is a $\Gamma$-weight of $M$, the vector space $M_\chi$ is called the $\Gamma$-weight subspace of $M$ with weight $\chi$ and the elements of $M_\chi$ are $\Gamma$-weight vectors with weight $\chi$. Moreover, if a $\Gamma$-module $M$ satisfies
\begin{align}
  M = \bigoplus_{\chi \in \Hom(\Gamma,\C)} M_\chi,
\end{align}
then we call $M$ a \emph{$\Gamma$-weight module}. The dimension of the $\Gamma$-weight subspace $M_\chi$ will be called the \emph{$\Gamma$-multiplicity} of $\chi$ in $M$.

Let $\mfrak{g}$ be a semisimple Lie algebra and let $\Gamma$ be a commutative subalgebra of the universal enveloping algebra $U(\mfrak{g})$ of $\mfrak{g}$. Then we denote by $\mcal{H}(\mfrak{g},\Gamma)$ the full subcategory of $\mcal{M}(\mfrak{g})$ consisting of $\Gamma$-weight $\mfrak{g}$-modules. Let us note that $\mcal{H}(\mfrak{g},\Gamma)$ is closed with respect to the operations of taking submodules and quotients. Besides, if $\Gamma$ contains the Cartan subalgebra $\mfrak{h}$, a $\Gamma$-weight $\mfrak{g}$-module $M$ is called a \emph{$\Gamma$-Gelfand--Tsetlin $\mfrak{g}$-module}.

%%%%%%%%%%%%%%%%%%%%%%%%%%%%%%%%%%%%%%%%%%%%%%%%%%%%%%%%%%%%%%%%%%%%%%%%%%%%%%%%%%%%%%%%%%
%%%%%%%%%%%%%%%%%%%%%%%%%%%%%%%%%%%%%%%%%%%%%%%%%%%%%%%%%%%%%%%%%%%%%%%%%%%%%%%%%%%%%%%%%%

\section{Affine Kac--Moody algebras and Weyl algebras}

In this section we define smooth and induced modules for affine Kac--Moody algebras and introduce a formalism for infinite-dimensional Weyl algebras.

%%%%%%%%%%%%%%%%%%%%%%%%%%%%%%%%%%%%%%%%%%%%%%%%%%%%%%%%%%%%%%%%%%%%%%%%%%%%%%%%%%%%%%%%%%

\subsection{Affine Kac--Moody algebras}

Let $\mfrak{g}$ be a semisimple (reductive) finite-dimensional Lie algebra and $\kappa$ be a $\mfrak{g}$-invariant symmetric bilinear form  on $\mfrak{g}$. The affine Kac--Moody algebra $\widehat{\mfrak{g}}_\kappa$ associated to $\mfrak{g}$ of level $\kappa$ is the $1$-dimensional central extension $\widehat{\mfrak{g}}_\kappa = \mfrak{g}(\!(t)\!) \oplus \C c$ of the formal loop algebra $\mfrak{g}(\!(t)\!)= \mfrak{g} \otimes_\C \C(\!(t)\!)$, with the commutation relations
\begin{align}
  [a \otimes f(t), b \otimes g(t)] = [a,b] \otimes f(t)g(t) - \kappa(a,b)\Res_{t=0} (f(t)dg(t))c, \label{eq:commutation relation}
\end{align}
where $c$ is the central element of $\widehat{\mfrak{g}}_\kappa$, $a, b \in \mfrak{g}$ and $f(t), g(t) \in \C(\!(t)\!)$. Let us note that Lie algebras $\widehat{\mfrak{g}}_\kappa$ and $\widehat{\mfrak{g}}_{\kappa'}$ for $\mfrak{g}$-invariant symmetric bilinear forms $\kappa$ and $\kappa'$ on $\mfrak{g}$ are isomorphic if $\kappa'=k \kappa$ for some $k \in \C^\times$. By introducing the notation $a_n = a \otimes t^n$ for $a \in \mfrak{g}$ and $n \in \Z$, the commutation relations \eqref{eq:commutation relation} can be simplified into the form
\begin{align}
  [a_m,b_n]=[a,b]_{m+n}+m \kappa(a,b) \delta_{m,-n} c \label{eq:commutation relation modes}
\end{align}
for $m,n \in \Z$ and $a,b \in \mfrak{g}$.

As $\mfrak{h}$ is a Cartan subalgebra of $\mfrak{g}$, we introduce a Cartan subalgebra $\smash{\widehat{\mfrak{h}}}$ of $\widehat{\mfrak{g}}_\kappa$ by
\begin{align*}
  \widehat{\mfrak{h}} = \mfrak{h} \otimes_\C \C 1 \oplus \C c.
\end{align*}
While any two Borel subalgebras of $\mfrak{g}$ are conjugate by an automorphism of $\mfrak{g}$, any two Borel subalgebras of $\widehat{\mfrak{g}}_\kappa$ may not be conjugate by an automorphism of $\widehat{\mfrak{g}}_\kappa$, see \cite{Futorny1997}.

Let $\mfrak{p}$ be a standard parabolic subalgebra of $\mfrak{g}$ with the nilradical $\mfrak{u}$, the opposite nilradical $\widebar{\mfrak{u}}$ and the Levi subalgebra $\mfrak{l}$. Then the \emph{standard parabolic subalgebra} $\widehat{\mfrak{p}}_{\rm st}$ of $\widehat{\mfrak{g}}_\kappa$ associated to $\mfrak{p}$ is given through
\begin{align*}
  \widehat{\mfrak{p}}_{\rm st} = \widehat{\mfrak{l}}_{\rm st} \oplus \widehat{\mfrak{u}}_{\rm st},
\end{align*}
where the Levi subalgebra $\smash{\widehat{\mfrak{l}}}_{\rm st}$ is defined by
\begin{align*}
  \widehat{\mfrak{l}}_{\rm st} = \mfrak{l} \otimes_\C \C 1 \oplus \C c
\end{align*}
and the nilradical $\widehat{\mfrak{u}}_{\rm st}$ and the opposite nilradical $\widehat{\widebar{\mfrak{u}}}_{\rm st}$ by
\begin{align*}
  \widehat{\mfrak{u}}_{\rm st} = \mfrak{u} \otimes_\C \C 1 \oplus \mfrak{g} \otimes_\C t\C[[t]] \qquad \text{and} \qquad
  \widehat{\widebar{\mfrak{u}}}_{\rm st} = \widebar{\mfrak{u}} \otimes_\C \C 1 \oplus \mfrak{g} \otimes_\C t^{-1}\C[t^{-1}].
\end{align*}
Moreover, we have the corresponding triangular decomposition
\begin{align*}
  \widehat{\mfrak{g}}_\kappa = \widehat{\widebar{\mfrak{u}}}_{\rm st} \oplus \widehat{\mfrak{l}}_{\rm st} \oplus \widehat{\mfrak{u}}_{\rm st}
\end{align*}
of the Lie algebra $\widehat{\mfrak{g}}_\kappa$. If $\mfrak{p} = \mfrak{b}$ then $\smash{\widehat{\mfrak{b}}}_{\rm st}$ is the standard Borel subalgebra, and if $\mfrak{p} = \mfrak{g}$ then $\widehat{\mfrak{g}}_{\rm st}$ is the maximal standard parabolic subalgebra.

On the other hand, the \emph{natural parabolic subalgebra} $\widehat{\mfrak{p}}_{\rm nat}$ of $\widehat{\mfrak{g}}_\kappa$ associated to $\mfrak{p}$ is given through
\begin{align*}
  \widehat{\mfrak{p}}_{\rm nat} = \widehat{\mfrak{l}}_{\rm nat} \oplus \widehat{\mfrak{u}}_{\rm nat},
\end{align*}
where the Levi subalgebra $\smash{\widehat{\mfrak{l}}}_{\rm nat}$ is defined by
\begin{align*}
  \widehat{\mfrak{l}}_{\rm nat} = \mfrak{l} \otimes_\C \C(\!(t)\!) \oplus \C c
\end{align*}
and the nilradical $\widehat{\mfrak{u}}_{\rm nat}$ and the opposite nilradical $\widehat{\widebar{\mfrak{u}}}_{\rm nat}$ by
\begin{align*}
  \widehat{\mfrak{u}}_{\rm nat} = \mfrak{u} \otimes_\C \C(\!(t)\!)
  \qquad \text{and} \qquad
  \widehat{\widebar{\mfrak{u}}}_{\rm nat} = \widebar{\mfrak{u}} \otimes_\C \C(\!(t)\!).
\end{align*}
We have again the corresponding triangular decomposition
\begin{align*}
  \widehat{\mfrak{g}}_\kappa = \widehat{\widebar{\mfrak{u}}}_{\rm nat} \oplus \widehat{\mfrak{l}}_{\rm nat} \oplus \widehat{\mfrak{u}}_{\rm nat}
\end{align*}
of the Lie algebra $\widehat{\mfrak{g}}_\kappa$. If $\mfrak{p} = \mfrak{b}$ then $\smash{\widehat{\mfrak{b}}}_{\rm nat}$ is the natural Borel subalgebra, and if $\mfrak{p} = \mfrak{g}$ then we have $\widehat{\mfrak{g}}_{\rm nat} = \widehat{\mfrak{g}}_\kappa$.

For more general parabolic subalgebras of the affine Kac--Moody algebra $\widehat{\mfrak{g}}_\kappa$ see \cite{Futorny-Kashuba2018}.
\medskip

To describe the root structure of the affine Kac--Moody algebra $\widehat{\mfrak{g}}_\kappa$ we define the extended affine Kac--Moody algebra $\widetilde{\mfrak{g}}_\kappa$ as
\begin{align*}
  \widetilde{\mfrak{g}}_\kappa = \mfrak{g}(\!(t)\!) \oplus \C c \oplus \C d,
\end{align*}
where we extend the Lie algebra structure of $\widehat{\mfrak{g}}_\kappa$ on $\widetilde{\mfrak{g}}_\kappa$ by the commutation relations
\begin{align}
  [d,a\otimes f(t)] = a\otimes t\partial_t f(t), \qquad [d,c]=0
\end{align}
for $a \in \mfrak{g}$ and $f(t) \in \C(\!(t)\!)$. The corresponding Cartan subalgebra $\smash{\widetilde{\mfrak{h}}}$ of $\widetilde{\mfrak{g}}_\kappa$ is defined by
\begin{align*}
  \widetilde{\mfrak{h}} = \mfrak{h} \otimes_\C \C 1 \oplus \C c \oplus \C d.
\end{align*}
The dual space to $\smash{\widetilde{\mfrak{h}}}$ we identify with $\mfrak{h}^* \oplus \C \Lambda_0 \oplus \C\delta$, where $\Lambda_0|_{\mfrak{h} \otimes_\C \C 1} = 0$, $\Lambda_0(c)=1$, $\Lambda_0(d)=0$ and $\delta|_{\mfrak{h} \otimes_\C \C 1} = 0$, $\delta(c)=0$, $\delta(d)=1$.
The set of roots of $\widehat{\mfrak{g}}_\kappa$ is then naturally a subset of $\mfrak{h}^* \oplus \C \Lambda_0 \oplus \C\delta$ given by
\begin{align}
  \widehat{\Delta} = \widehat{\Delta}^{\rm re} \cup \widehat{\Delta}^{\rm im},
\end{align}
where
\begin{align}
  \widehat{\Delta}^{\rm re} = \{\alpha+n\delta;\, \alpha \in \Delta,\, n \in\Z\} \qquad \text{and} \qquad \widehat{\Delta}^{\rm im} = \{n\delta;\, n \in \Z\setminus\{0\}\}
\end{align}
are the sets of real and imaginary roots, respectively. Moreover, we have the following decomposition $\widehat{\Delta}=\widehat{\Delta}_+\cup \widehat{\Delta}_-$ into positive and negative roots, where
\begin{align*}
  \widehat{\Delta}_{+}=\{\alpha+n\delta;\,  \alpha\in \Delta_+,\, n\in \N_0\} \cup
   \{-\alpha+n\delta;\, \alpha\in \Delta_+,\, n\in \N\} \cup \{n\delta;\, n \in \N\}
\end{align*}
and $\widehat{\Delta}_-=-\widehat{\Delta}_+$. We also set $\widehat{\Delta}^{\rm re}_\pm=\widehat{\Delta}_\pm \cap \widehat{\Delta}^{\rm re}$ and $\widehat{\Delta}^{\rm im}_\pm=\widehat{\Delta}_\pm \cap \widehat{\Delta}^{\rm im}$.
\medskip

Let us note that the grading element $d \in \widetilde{\mfrak{g}}_\kappa$ gives $\widetilde{\mfrak{g}}_\kappa$ and also $\widehat{\mfrak{g}}_\kappa$ the structure of $\Z$-graded topological Lie algebras with the gradation defined by $-d$, i.e.\ we have $\deg c = 0$, $\deg d =0$ and $\deg a_n = -n$ for $a \in \mfrak{g}$, $n\in \Z$.

%%%%%%%%%%%%%%%%%%%%%%%%%%%%%%%%%%%%%%%%%%%%%%%%%%%%%%%%%%%%%%%%%%%%%%%%%%%%%%%%%%%%%%%%%%

\subsection{Smooth modules over affine Kac--Moody algebras}
\label{subsec:induced}

We denote by $\mcal{M}(\widehat{\mfrak{g}}_\kappa)$ the category of $\widehat{\mfrak{g}}_\kappa$-modules. However, since the objects of this category may be very complicated, we focus our attention to some nice full subcategories of $\mcal{M}(\widehat{\mfrak{g}}_\kappa)$.
\medskip

\definition{Let $M$ be a $\widehat{\mfrak{g}}_\kappa$-module. We say that $M$ is a \emph{smooth} $\widehat{\mfrak{g}}_\kappa$-module if for each vector $v \in M$ there exists a positive integer $N_v \in \N$ such that
\begin{align*}
  (\mfrak{g} \otimes_\C t^{N_v}\C[[t]]) v = 0,
\end{align*}
or in other words that the Lie subalgebra $\mfrak{g} \otimes_\C t^{N_v}\C[[t]]$ annihilates $v$.
The category of smooth $\widehat{\mfrak{g}}_\kappa$-modules on which the central element $c$ acts as the identity we will denote by $\mcal{E}(\widehat{\mfrak{g}}_\kappa)$.}

Let us recall that by a graded $\widehat{\mfrak{g}}_\kappa$-module $M$ we mean a $\C$-graded vector space $M$ having the structure of a $\widehat{\mfrak{g}}_\kappa$-module compatible with the gradation of $\widehat{\mfrak{g}}_\kappa$. Let us note that by shifting a given gradation on $M$ by a complex number we obtain a new gradation on $M$.
\medskip

\definition{Let $M$ be a graded $\widehat{\mfrak{g}}_\kappa$-module. We say that $M$ is a \emph{positive energy} $\widehat{\mfrak{g}}_\kappa$-module if $M=\bigoplus_{n=0}^\infty M_{\lambda+n}$ and $M_\lambda \neq 0$, where $\lambda \in \C$. The category of positive energy $\widehat{\mfrak{g}}_\kappa$-modules on which the central element $c$ acts as the identity we will denote by $\mcal{E}_+\!(\widehat{\mfrak{g}}_\kappa)$.}

Let us note that if $M$ is a positive energy $\widehat{\mfrak{g}}_\kappa$-module, then it follows immediately that $M$ is also a smooth $\widehat{\mfrak{g}}_\kappa$-module. Therefore, the category $\mcal{E}_+\!(\widehat{\mfrak{g}}_\kappa)$ is a full subcategory of $\mcal{E}(\widehat{\mfrak{g}}_\kappa)$.
\medskip

Let us recall that the category $\mcal{M}(\widehat{\mfrak{g}}_\kappa)$ coincides with the category of modules over the universal enveloping algebra $U(\widehat{\mfrak{g}}_\kappa)$. There exists an analogous associative algebra for the category $\mcal{E}(\widehat{\mfrak{g}}_\kappa)$ which is constructed as follows, see \cite{Frenkel2007-book}. Since the central element $c$ acts as the identity on all $\widehat{\mfrak{g}}_\kappa$-modules from the category $\mcal{E}(\widehat{\mfrak{g}}_\kappa)$, the action of $U(\widehat{\mfrak{g}}_\kappa)$ factors through the quotient algebra
\begin{align*}
  U_c(\widehat{\mfrak{g}}_\kappa) = U(\widehat{\mfrak{g}}_\kappa)/(c-1).
\end{align*}
Further, let us introduce a linear topology on $U_c(\widehat{\mfrak{g}}_\kappa)$ in which the basis of open neighbourhoods of $0$ are the left ideals $I_N$ defined by
\begin{align*}
  I_N = U_c(\widehat{\mfrak{g}}_\kappa)(\mfrak{g} \otimes_\C t^N\C[[t]])
\end{align*}
for $N \in \N_0$. Let $\smash{\widetilde{U}_c(\widehat{\mfrak{g}}_\kappa)}$ be the completion of $U_c(\widehat{\mfrak{g}}_\kappa)$ with respect to this linear topology, i.e.\ we get
\begin{align*}
  \smash{\widetilde{U}_c(\widehat{\mfrak{g}}_\kappa)} = \lim_{\longleftarrow} U_c(\widehat{\mfrak{g}}_\kappa)/I_N.
\end{align*}
Then the structure of an associative algebra on $U_c(\widehat{\mfrak{g}}_\kappa)$ extends to the structure of an associative algebra on $\smash{\widetilde{U}_c(\widehat{\mfrak{g}}_\kappa)}$ by continuity. Hence, we obtain that $\smash{\widetilde{U}_c(\widehat{\mfrak{g}}_\kappa)}$ is a complete topological associative algebra, which we will call the \emph{completed universal enveloping algebra} of $\widehat{\mfrak{g}}_\kappa$. Moreover, the category $\mcal{E}(\widehat{\mfrak{g}}_\kappa)$ coincides with the category of discrete modules over the associative algebra $\smash{\widetilde{U}_c(\widehat{\mfrak{g}}_\kappa)}$ on which the action of $\smash{\widetilde{U}_c(\widehat{\mfrak{g}}_\kappa)}$ is pointwise continuous.
\medskip

Now, we construct a class of $\widehat{\mfrak{g}}_\kappa$-modules, the so called induced modules, which belong to the category $\mcal{E}(\widehat{\mfrak{g}}_\kappa)$. Let $\widehat{\mfrak{p}}_{\rm st} = \smash{\widehat{\mfrak{l}}}_{\rm st} \oplus \widehat{\mfrak{u}}_{\rm st}$ be the standard parabolic subalgebra of $\widehat{\mfrak{g}}_\kappa$ associated to a standard parabolic subalgebra $\mfrak{p}$ of $\mfrak{g}$ and let $E$ be an $\mfrak{l}$-module. Then the induced $\widehat{\mfrak{g}}_\kappa$-module
\begin{align*}
\mathbb{M}_{\kappa,\mfrak{p}}(E) = U(\widehat{\mfrak{g}}_\kappa)\otimes_{U(\widehat{\mfrak{p}}_{\rm st})}\! E,
\end{align*}
where $E$ is considered as the $\widehat{\mfrak{p}}_{\rm st}$-module on which $\widehat{\mfrak{u}}_{\rm st}$ acts trivially and $c$ acts as the identity, has a unique maximal $\widehat{\mfrak{g}}_\kappa$-submodule $\mathbb{K}_{\kappa,\mfrak{p}}(E)$ having zero intersection with the $\mfrak{l}$-submodule $E$ of $\mathbb{M}_{\kappa,\mfrak{p}}(E)$. Therefore, we may set
\begin{align*}
  \mathbb{L}_{\kappa,\mfrak{p}}(E) = \mathbb{M}_{\kappa,\mfrak{p}}(E)/\mathbb{K}_{\kappa,\mfrak{p}}(E)
\end{align*}
for an $\mfrak{l}$-module $E$. Moreover, it is easy to see that if $E$ is a simple $\mfrak{l}$-module, then $\mathbb{L}_{\kappa,\mfrak{p}}(E)$ is also a simple $\widehat{\mfrak{g}}_\kappa$-module. The $\widehat{\mfrak{g}}_\kappa$-module $\mathbb{M}_{\kappa,\mfrak{p}}(E)$ is called the \emph{generalized Verma module} induced from $E$ for the standard parabolic subalgebra $\widehat{\mfrak{p}}_{\rm st}$.

Further, if we set $F=U(\widebar{\mfrak{u}})E \subset \mathbb{M}_{\kappa,\mfrak{p}}(E)$, then $F$ is a $\widehat{\mfrak{g}}_{\rm st}$-module such that $\mfrak{g} \otimes_\C t\C[[t]]$ acts trivially and $c$ acts as the identity. The induced $\widehat{\mfrak{g}}_\kappa$-module $\mathbb{M}_{\kappa,\mfrak{g}}(F)$ is isomorphic to $\mathbb{M}_{\kappa,\mfrak{p}}(E)$ and $\mathbb{L}_{\kappa,\mfrak{g}}(F) \simeq \mathbb{L}_{\kappa,\mfrak{p}}(E)$ if $F$ is a simple $\mfrak{g}$-module. We will always consider induced $\widehat{\mfrak{g}}_\kappa$-modules as generalized Verma modules for the maximal standard parabolic subalgebra $\widehat{\mfrak{g}}_{\rm st}$.

Therefore, we have the induction functor
\begin{align}
  \mathbb{M}_{\kappa,\mfrak{g}} \colon \mcal{M}(\mfrak{g}) \rarr \mcal{E}_+\!(\widehat{\mfrak{g}}_\kappa)
\end{align}
and the functor
\begin{align}
  \mathbb{L}_{\kappa,\mfrak{g}} \colon \mcal{M}(\mfrak{g}) \rarr \mcal{E}_+\!(\widehat{\mfrak{g}}_\kappa).
\end{align}
Let us recall that if $E$ is a simple finite-dimensional $\mfrak{g}$-module, then $\mathbb{M}_{\kappa,\mfrak{g}}(E)$ is usually called the \emph{Weyl module}.
\medskip

In the same way as $\mcal{E}(\widehat{\mfrak{g}}_\kappa)$ we may define the category $\mcal{E}(\widetilde{\mfrak{g}}_\kappa)$. Moreover, as $\widehat{\mfrak{g}}_\kappa$ is a Lie subalgebra of $\widetilde{\mfrak{g}}_\kappa$, we have also a natural forgetful functor $\mcal{E}(\widetilde{\mfrak{g}}_\kappa) \rarr \mcal{E}(\widehat{\mfrak{g}}_\kappa)$. On the other hand, for a non-critical level $\kappa$ we can view on $\mcal{E}(\widehat{\mfrak{g}}_\kappa)$ as a full subcategory of $\mcal{E}(\widetilde{\mfrak{g}}_\kappa)$. Let $M$ be a smooth $\widehat{\mfrak{g}}_\kappa$-module on which the central element $c$ acts as the identity. Since $\kappa$ is a non-critical level, any smooth $\widehat{\mfrak{g}}_\kappa$-module carries an action of the Virasoro algebra obtained by the Segal--Sugawara construction, and so in particular an action of $L_0$ (the nontrivial semisimple generator). Hence, the action of the grading element $d$ of $\widetilde{\mfrak{g}}_\kappa$ is then defined as the action of $-L_0$.
However, if $\kappa$ is not a non-critical level, then general smooth $\widehat{\mfrak{g}}_\kappa$-modules do not necessarily carry an action of $L_0$.
\medskip

In the rest of this section we described non-critical levels for a semisimple Lie algebra $\mfrak{g}$. Let us consider a $\mfrak{b}$-invariant symmetric bilinear form on $\mfrak{b}$ defined by
\begin{align}
  \kappa_c^\mfrak{b}(a,b) = - \tr_{\mfrak{g}/\mfrak{b}}(\ad(a)\ad(b))
\end{align}
for $a, b \in \mfrak{b}$.
\medskip

\definition{Let $\kappa$ be a $\mfrak{g}$-invariant symmetric bilinear form on $\mfrak{g}$. We say that $\kappa$ is \emph{non-critical} if $\kappa-\kappa_c^\mfrak{b}$ is non-degenerate on $\mfrak{h}$, \emph{partially critical} if $\kappa-\kappa_c^\mfrak{b}$ is degenerate on $\mfrak{h}$, and \emph{critical} if $\kappa-\kappa_c^\mfrak{b}$ is zero on $\mfrak{h}$. The critical $\mfrak{g}$-invariant symmetric bilinear form on $\mfrak{g}$ we will denote by $\kappa_c$.}

Since $\mfrak{g}$ is a semisimple Lie algebra, we have the direct sum decomposition
\begin{align*}
  \mfrak{g} = \bigoplus_{i=1}^r \mfrak{g}_i
\end{align*}
of $\mfrak{g}$ into the direct sum of simple Lie algebras $\mfrak{g}_i$ for $i=1,2,\dots,r$ such that these direct summands are mutually orthogonal with respect to the Cartan--Killing form $\kappa_\mfrak{g}$ on $\mfrak{g}$. We denote by $\kappa_0^{\mfrak{g}_i}$ the normalized $\mfrak{g}_i$-invariant symmetric bilinear form on $\mfrak{g}_i$, i.e.\ we have $\kappa_{\mfrak{g}_i} = 2h^\vee_i \kappa_0^{\mfrak{g}_i}$, where $h_i^\vee$ is the dual Coxeter number of $\mfrak{g}_i$, for $i=1,2,\dots,r$. We have the following criterion.
\medskip

\lemma{Let $\kappa$ be a $\mfrak{g}$-invariant symmetric bilinear form on $\mfrak{g}$. Then $\kappa$ is partially critical if $\kappa_{|\mfrak{g}_i} = -h^\vee_i \kappa_0^{\mfrak{g}_i}$ for some $i=1,2,\dots,r$, and $\kappa$ is critical if $\kappa_{|\mfrak{g}_i} = -h^\vee_i \kappa_0^{\mfrak{g}_i}$ for all $i=1,2,\dots,r$.}

\proof{For $a,b \in \mfrak{h}$, we have
\begin{align*}
  \kappa_c^\mfrak{b}(a,b) = -\tr_{\mfrak{g}/\mfrak{b}}(\ad(a)\ad(b)) = -\sum_{\alpha \in \Delta_+} \alpha(a)\alpha(b).
\end{align*}
Further, for each $i=1,2,\dots,r$ there exists $k_i \in \C$ such that $\smash{\kappa_{|\mfrak{g}_i}} =k_i\kappa_0^{\mfrak{g}_i}$, since $\mfrak{g}_i$ is a simple Lie algebra. Hence, we may write
\begin{align*}
  (\kappa-\kappa_c^\mfrak{b})(a,b) &= (k_i\kappa_0^{\mfrak{g}_i}-\kappa_c^\mfrak{b})(a,b) = {k_i \over 2h_i^\vee}\, \kappa_{\mfrak{g}_i}(a,b) - \kappa_c^\mfrak{b}(a,b) = {k_i + h_i^\vee \over h_i^\vee} \sum_{\alpha \in \Delta_+^{\mfrak{g}_i}} \alpha(a)\alpha(b)
\end{align*}
for $a,b \in \mfrak{h} \cap \mfrak{g}_i$, where $\Delta_+^{\mfrak{g}_i} \subset \Delta_+$ is the set of positive roots of $\mfrak{g}_i$. The required statement then follows immediately.}

\vspace{-2mm}

%%%%%%%%%%%%%%%%%%%%%%%%%%%%%%%%%%%%%%%%%%%%%%%%%%%%%%%%%%%%%%%%%%%%%%%%%%%%%%%%%%%%%%%%%%

\subsection{Weyl algebras}

Let us consider the commutative algebra $\mcal{K}=\C(\!(t)\!)$ with the subalgebra $\mcal{O}=\C[[t]]$. Further, let $\Omega_\mcal{K} = \C(\!(t)\!)dt$ and $\Omega_\mcal{O} = \C[[t]]dt$ be the corresponding modules of K\"ahler differentials. Then for a finite-dimensional complex vector space $V$ we introduce the infinite-dimensional vector spaces $\mcal{K}(V) = V \otimes_\C \mcal{K}$ and $\Omega_\mcal{K}(V) = V \otimes_\C \Omega_\mcal{K}$. Using the pairing $(\,\cdot\,,\cdot\,) \colon \Omega_\mcal{K}(V^*) \otimes \mcal{K}(V) \rarr \C$ defined by
\begin{align}
  (\alpha \otimes f(t)dt, v \otimes g(t)) = \alpha(v) \Res_{t=0}(g(t)f(t)dt), \label{eq:non-degenerate form}
\end{align}
where $\alpha \in V^*$, $v \in V$ and $f(t),g(t) \in \mcal{K}$, we identify the restricted dual space to $\mcal{K}(V)$ with the vector space $\Omega_\mcal{K}(V^*)$, and vice versa. Moreover, the pairing \eqref{eq:non-degenerate form} gives us a skew-symmetric non-degenerate bilinear form $\langle \cdot\,,\cdot \rangle$ on $\Omega_\mcal{K}(V^*) \oplus \mcal{K}(V)$ defined by
\begin{align*}
\langle \alpha \otimes f(t)dt,v \otimes g(t) \rangle = - \langle v \otimes g(t), \alpha \otimes f(t)dt \rangle = (\alpha \otimes f(t)dt, v \otimes g(t))
\end{align*}
for $\alpha \in V^*$, $v \in V$, $f(t), g(t) \in \mcal{K}$, and
\begin{align*}
  \langle v \otimes f(t), w \otimes g(t) \rangle = \langle \alpha \otimes f(t)dt, \beta \otimes g(t)dt \rangle = 0
\end{align*}
for $\alpha,\beta \in V^*$, $v,w \in V$, $f(t),g(t) \in \mcal{K}$. The Weyl algebras $\eus{A}_{\mcal{K}(V)}$ and $\eus{A}_{\Omega_\mcal{K}(V^*)}$ are given by
\begin{align*}
  \eus{A}_{\mcal{K}(V)} = T(\Omega_\mcal{K}(V^*) \oplus \mcal{K}(V))/I_{\mcal{K}(V)} \qquad \text{and} \qquad \eus{A}_{\Omega_\mcal{K}(V^*)} = T(\Omega_\mcal{K}(V^*) \oplus \mcal{K}(V))/I_{\Omega_\mcal{K}(V^*)},
\end{align*}
where $I_{\mcal{K}(V)}$ and $I_{\Omega_\mcal{K}(V^*)}$ denote the two-sided ideals of the tensor algebra $T(\Omega_\mcal{K}(V^*) \oplus \mcal{K}(V))$ generated by $a \otimes b - b \otimes a +\langle a,b \rangle \cdot 1$ for $a, b  \in \Omega_\mcal{K}(V^*) \oplus \mcal{K}(V)$ and by $a \otimes b - b \otimes a -\langle a,b \rangle \cdot 1$ for $a,b \in \Omega_\mcal{K}(V^*) \oplus \mcal{K}(V)$, respectively. The algebras of polynomials on $\mcal{K}(V)$ and $\Omega_\mcal{K}(V^*)$ are defined by $\Pol \mcal{K}(V) = S(\Omega_\mcal{K}(V^*))$ and $\Pol \Omega_\mcal{K}(V^*) = S(\mcal{K}(V))$.

Let $\mcal{L}$ and $\mcal{L}^{\rm c}$ be complementary Lagrangian (maximal isotropic) subspaces of $\Omega_\mcal{K}(V^*) \oplus \mcal{K}(V)$, i.e.\ we have $\Omega_\mcal{K}(V^*) \oplus \mcal{K}(V) = \mcal{L} \oplus \mcal{L}^{\rm c}$. Then the symmetric algebra $S(\mcal{L})$ is a subalgebra of $\eus{A}_{\mcal{K}(V)}$ since the elements of $\mcal{L}$ commute in $\eus{A}_{\mcal{K}(V)}$. Moreover, it is a maximal commutative subalgebra of $\eus{A}_{\mcal{K}(V)}$. Further, let us consider the induced $\eus{A}_{\mcal{K}(V)}$-module
\begin{align*}
  \Ind_{S(\mcal{L})}^{\eus{A}_{\mcal{K}(V)}} \!\C \simeq S(\mcal{L}^{\rm c}),
\end{align*}
where $\C$ is the trivial $1$-dimensional $S(\mcal{L})$-module. It follows immediately that the induced $\eus{A}_{\mcal{K}(V)}$-module has the natural structure of a commutative algebra. We denote by $M_\mcal{L}$ the commutative algebra which is the completion of $S(\mcal{L}^{\rm c})$ with
respect to the linear topology on $S(\mcal{L}^{\rm c})$ in which the basis of open neighbourhoods of $0$ are the subspaces $\mcal{I}_{n,m}$ for $n,m\in \Z$, where $\mcal{I}_{n,m}$ is the ideal of $S(\mcal{L}^{\rm c})$ generated by $\mcal{L}^{\rm c}\cap (V^*\!\otimes_\C t^n\Omega_\mcal{O} \oplus V\! \otimes_\C t^m\mcal{O})$. In addition, we may extend the action of the Weyl algebra $\eus{A}_{\mcal{K}(V)}$ to $M_\mcal{L}$.
\medskip

Our next step is to pass to a completion of the Weyl algebra $\smash{\eus{A}_{\mcal{K}(V)}}$,
because $\smash{\eus{A}_{\mcal{K}(V)}}$ is not sufficiently large for our considerations.
Let us denote by $\Fun\mcal{K}(V)$ the completion of the commutative algebra $\Pol \mcal{K}(V)$ with respect to the linear topology on $\Pol \mcal{K}(V)$ in which the basis of open neighbourhoods of $0$ are the subspaces $\mcal{J}_n$ for $n \in \Z$, where $\mcal{J}_n$ is the ideal of $\Pol \mcal{K}(V)$ generated by $V^*\! \otimes_\C t^n\Omega_\mcal{O}$. Consequently, we have $\Fun \mcal{K}(V)= M_{\mcal{K}(V)}$. Then a vector field on $\mcal{K}(V)$ is by definition a continuous linear endomorphism $\xi$ of $\Fun \mcal{K}(V)$ which satisfies the Leibniz rule
\begin{align*}
  \xi(fg)=\xi(f)g+f\xi(g)
\end{align*}
for all $f,g \in \Fun \mcal{K}(V)$. In other words, a vector field on $\mcal{K}(V)$ is a linear endomorphism $\xi$ of $\Fun \mcal{K}(V)$ such that for any $m \in \Z$ there exists $n \in \Z$, $m \leq n$ and a derivation
\begin{align*}
  \xi_{n,m} \colon \Pol \mcal{K}(V)/\mcal{J}_n \rarr \Pol \mcal{K}(V)/\mcal{J}_m
\end{align*}
satisfying
\begin{align*}
  s_m(\xi(f))=\xi_{n,m}(s_n(f))
\end{align*}
for all $f \in \Fun \mcal{K}(V)$, where
\begin{align*}
  s_n \colon \Fun \mcal{K}(V) \rarr \Pol \mcal{K}(V)/\mcal{J}_n
\end{align*}
is the canonical homomorphism of algebras. The vector space of all vector fields is naturally a topological Lie algebra, which we denote by $\Vect \mcal{K}(V)$. There is a non-split short exact sequence
\begin{align*}
  0 \rarr \Fun \mcal{K}(V) \rarr \widetilde{\eus{A}}_{\mcal{K}(V),\leq 1} \rarr \Vect \mcal{K}(V) \rarr 0
\end{align*}
of topological Lie algebras, see \cite{Frenkel-Ben-Zvi2004} for more details. This extension of the topological Lie algebra $\Vect \mcal{K}(V)$ by its module $\Fun \mcal{K}(V)$ is however different from the standard split extension
\begin{align*}
  0 \rarr \Fun \mcal{K}(V) \rarr \eus{A}^\sharp_{\mcal{K}(V),\leq 1} \rarr \Vect \mcal{K}(V) \rarr 0
\end{align*}
of $\Vect \mcal{K}(V)$ by $\Fun \mcal{K}(V)$. The completed Weyl algebra $\smash{\widetilde{\eus{A}}_{\mcal{K}(V)}}$ is then defined as the associative algebra over $\Fun \mcal{K}(V)$  generated by the images of the homomorphisms
 $i \colon \Fun\mcal{K}(V) \rarr \smash{\widetilde{\eus{A}}_{\mcal{K}(V)}}$ (as associative algebras) and  $j \colon  \smash{\widetilde{\eus{A}}_{\mcal{K}(V),\leq 1}} \rarr \smash{\widetilde{\eus{A}}_{\mcal{K}(V)}}$ (as Lie algebras),
  with the relation
\begin{align*}
  [j(P),i(f)]=i(a(P)(f))
\end{align*}
for $f \in \Fun \mcal{K}(V)$ and $P \in \smash{\widetilde{\eus{A}}_{\mcal{K}(V),\leq 1}}$, where $a \colon \smash{\widetilde{\eus{A}}_{\mcal{K}(V),\leq 1}} \rarr \Vect \mcal{K}(V)$ is the homomorphism of Lie algebras originated from the corresponding short exact sequence.
\medskip

Now, we define a class of $\eus{A}_{\mcal{K}(V)}$-modules called induced modules. Let us consider the vector subspaces $\mcal{L}_+$, $\mcal{L}_-$ and $\mcal{L}_0$ of $\Omega_\mcal{K}(V^*) \oplus \mcal{K}(V)$ given by
\begin{gather*}
  \mcal{L}_- = V^*\! \otimes_\C t^{-2}\C[t^{-1}]dt \oplus V\! \otimes_\C t^{-1}\C[t^{-1}], \qquad  \mcal{L}_+ = V^*\! \otimes_\C \C[[t]]dt \oplus V\! \otimes_\C t\C[[t]], \\
  \mcal{L}_0 = V^*\! \otimes_\C \C t^{-1}dt \oplus V\! \otimes_\C \C 1.
\end{gather*}
Then we have the direct sum decomposition
\begin{align*}
  \Omega_\mcal{K}(V^*) \oplus \mcal{K}(V) = \mcal{L}_- \oplus \mcal{L}_0 \oplus \mcal{L}_+
\end{align*}
of $\Omega_\mcal{K}(V^*) \oplus \mcal{K}(V)$, which induces the triangular decomposition
\begin{align*}
  \eus{A}_{\mcal{K}(V)} \simeq \eus{A}_{\mcal{K}(V),-} \otimes_\C \eus{A}_{\mcal{K}(V),0} \otimes_\C \eus{A}_{\mcal{K}(V),+}
\end{align*}
of the Weyl algebra $\eus{A}_{\mcal{K}(V)}$, where
\begin{align*}
  \eus{A}_{\mcal{K}(V),-} \simeq S(\mcal{L}_-), \qquad \eus{A}_{\mcal{K}(V),0} \simeq \eus{A}_V, \qquad \eus{A}_{\mcal{K}(V),+} \simeq S(\mcal{L}_+),
\end{align*}
and $\eus{A}_V$ is the Weyl algebra of the finite-dimensional vector space $V$. Moreover, the Weyl algebra $\eus{A}_{\mcal{K}(V)}$ is a $\Z$-graded algebra with the gradation determined by
\begin{align*}
  \deg(v \otimes t^n) = -n, \qquad \deg 1 = 0, \qquad \deg(\alpha \otimes t^{-n-1}dt) = n
\end{align*}
for $v \in V$, $\alpha \in V^*$ and $n \in \Z$.
\medskip

\definition{Let $M$ be an $\eus{A}_{\mcal{K}(V)}$-module. Then we say that $M$ is a \emph{smooth} $\eus{A}_{\mcal{K}(V)}$-module if for each vector $v \in M$ there exists a positive integer $N_v \in \N$ such that
\begin{align*}
  (V^*\! \otimes_\C t^{N_v}\C[[t]]dt \oplus V\! \otimes_\C t^{N_v}\C[[t]]) v = 0.
\end{align*}
The category of smooth $\eus{A}_{\mcal{K}(V)}$-modules we will denote by $\mcal{E}(\eus{A}_{\mcal{K}(V)})$.}

Completely analogously as for $\widehat{\mfrak{g}}_\kappa$-modules, we may introduce graded $\eus{A}_{\mcal{K}(V)}$-modules and positive energy $\eus{A}_{\mcal{K}(V)}$-modules. The category of positive energy $\eus{A}_{\mcal{K}(V)}$-modules we will denote by $\mcal{E}_+\!(\eus{A}_{\mcal{K}(V)})$.
\medskip

Let $E$ be an $\eus{A}_V$-module. Then the induced $\eus{A}_{\mcal{K}(V)}$-module
\begin{align*}
  \mathbb{M}_{\mcal{K}(V)}(E) = \eus{A}_{\mcal{K}(V)} \otimes_{\eus{A}_{\mcal{K}(V),0} \otimes_\C \eus{A}_{\mcal{K}(V),+}}\! E,
\end{align*}
where $E$ is considered as the $\eus{A}_{\mcal{K}(V),0} \otimes_\C \eus{A}_{\mcal{K}(V),+}$-module on which the Weyl algebra $\eus{A}_{\mcal{K}(V),0}$ acts via the canonical isomorphism $\eus{A}_{\mcal{K}(V),0} \simeq \eus{A}_V$ and $\eus{A}_{\mcal{K}(V),+}$ acts trivially, has a unique maximal $\eus{A}_{\mcal{K}(V)}$-submodule $\mathbb{K}_{\eus{A}_{\mcal{K}(V)}}(E)$ having zero intersection with the $\eus{A}_V$-submodule $E$ of $\mathbb{M}_{\mcal{K}(V)}(E)$. Therefore, we may set
\begin{align*}
  \mathbb{L}_{\mcal{K}(V)}(E) = \mathbb{M}_{\mcal{K}(V)}(E)/\mathbb{K}_{\mcal{K}(V)}(E)
\end{align*}
for an $\eus{A}_V$-module $E$. Moreover, it is easy to see that if $E$ is a simple $\eus{A}_V$-module, then $\mathbb{L}_{\mcal{K}(V)}(E)$ is also a simple $\eus{A}_{\mcal{K}(V)}$-module.

Therefore, we have the induction functor
\begin{align}
  \mathbb{M}_{\mcal{K}(V)} \colon \mcal{M}(\eus{A}_V) \rarr \mcal{E}_+\!(\eus{A}_{\mcal{K}(V)})
\end{align}
and the functor
\begin{align}
  \mathbb{L}_{\mcal{K}(V)} \colon \mcal{M}(\eus{A}_V) \rarr \mcal{E}_+\!(\eus{A}_{\mcal{K}(V)})
\end{align}
Besides, it follows immediately that $\mathbb{M}_{\mcal{K}(V)}(E)$ and $\mathbb{L}_{\mcal{K}(V)}(E)$ are also $\smash{\widetilde{\eus{A}}}_{\mcal{K}(V)}$-modules.
\medskip

Let us consider a Borel subalgebra $\mfrak{b}$ of a semisimple Lie algebra $\mfrak{g}$ with the nilradical $\mfrak{n}$, the opposite nilradical $\widebar{\mfrak{n}}$ and the Cartan subalgebra $\mfrak{h}$. Further, let $\{f_\alpha;\, \alpha \in \Delta_+\}$ be a root basis of the opposite nilradical $\widebar{\mfrak{n}}$. We denote by $\{x_\alpha;\, \alpha \in \Delta_+\}$ the linear coordinate functions on $\widebar{\mfrak{n}}$ with respect to the given basis of $\widebar{\mfrak{n}}$.
Then the set $\{f_\alpha \otimes t^n;\, \alpha \in \Delta_+,\, n \in \Z\}$ forms a topological basis of $\mcal{K}(\widebar{\mfrak{n}}) =\widebar{\mfrak{n}}_{{\rm nat}}$, and the set $\{x_\alpha \otimes t^{-n-1}dt;\, \alpha \in \Delta_+,\, n \in \Z\}$ forms the dual topological basis of $\Omega_\mcal{K}(\widebar{\mfrak{n}}^*) \simeq (\widebar{\mfrak{n}}_{{\rm nat}})^*$ with respect to the pairing \eqref{eq:non-degenerate form}, i.e.\ we have
\begin{align*}
  (x_\alpha \otimes t^{-n-1}dt,f_\beta \otimes t^m)=x_\alpha(f_\beta)\Res_{t=0} t^{m-n-1}dt= \delta_{\alpha,\beta} \delta_{n,m}
\end{align*}
for $\alpha,\beta \in \Delta_+$ and $m,n \in \Z$. If we denote $x_{\alpha,n}=x_\alpha \otimes t^{-n-1}dt$ and $\partial_{x_{\alpha,n}}=f_\alpha \otimes t^n$ for $\alpha \in \Delta_+$ and $n \in \Z$, then the two-sided ideal $I_{\mcal{K}(\widebar{\mfrak{n}})}$ is generated by elements
\begin{align*}
\bigg(\sum_{n \in \Z} a_nx_{\alpha,n}\bigg)\! \otimes \!\bigg(\sum_{m\in \Z} b_m\partial_{x_{\beta,m}}\! \bigg) - \bigg(\sum_{m\in \Z} b_m\partial_{x_{\beta,m}}\! \bigg) \!\otimes\! \bigg(\sum_{n \in \Z} a_nx_{\alpha,n}\bigg)+\delta_{\alpha,\beta}\bigg(\sum_{n\in \Z} a_nb_n\bigg)\!\cdot 1
\end{align*}
and it coincides with the canonical commutation relations
\begin{align*}
[x_{\alpha,n},\partial_{x_{\beta,m}}]=-\delta_{\alpha,\beta}\delta_{n,m}
\end{align*}
for $\alpha,\beta \in \Delta_+$ and $m,n \in \Z$. Therefore, we obtain that the Weyl algebra $\eus{A}_{\mcal{K}(\widebar{\mfrak{n}})}$ is topologically generated by the set $\{x_{\alpha,n}, \partial_{x_{\alpha,n}};\, \alpha \in \Delta_+,\, n \in \Z\}$ with the canonical commutation relations.

%%%%%%%%%%%%%%%%%%%%%%%%%%%%%%%%%%%%%%%%%%%%%%%%%%%%%%%%%%%%%%%%%%%%%%%%%%%%%%%%%%%%%%%%%%
%%%%%%%%%%%%%%%%%%%%%%%%%%%%%%%%%%%%%%%%%%%%%%%%%%%%%%%%%%%%%%%%%%%%%%%%%%%%%%%%%%%%%%%%%%

\section{Representations of vertex algebras}

%%%%%%%%%%%%%%%%%%%%%%%%%%%%%%%%%%%%%%%%%%%%%%%%%%%%%%%%%%%%%%%%%%%%%%%%%%%%%%%%%%%%%%%%%%

\subsection{Vertex algebras}

In this section we recall some notions and basic facts on vertex algebra, for more details see \cite{Borcherds1986}, \cite{Kac1998}, \cite{Dong-Li-Mason1998}, \cite{Frenkel-Ben-Zvi2004}, \cite{Frenkel2007-book}.
\medskip

Let $R$ be an algebra over $\C$. Then an $R$-valued formal power series (or formal distribution) in the variables $z_1,\dots,z_n$ is a series
\begin{align*}
  a(z_1,\dots,z_n)=\sum_{m_1,\dots,m_n \in \Z} a_{m_1,\dots,m_n} z_1^{m_1} \dotsm z_n^{m_n},
\end{align*}
where $a_{m_1,\dots,m_n} \in R$. The complex vector space of all $R$-valued formal power series is denoted by $R[[z_1^{\pm 1},\dots,z_n^{\pm 1}]]$. For a formal power series $a(z) = \sum_{m \in \Z} a_m z^m$, the residue is defined by
\begin{align*}
  \Res_{z=0} a(z) = a_{-1}.
\end{align*}
A particulary important example of a $\C$-valued formal power series in two variables $z$, $w$ is the formal delta function $\delta(z-w)$ given by
\begin{align*}
  \delta(z-w) = \sum_{m \in \Z} z^m w^{-m-1}.
\end{align*}

Let $V$ be a complex vector space, so $\End V$ is an algebra over $\C$. We say that a formal power series $a(z) \in \End V[[z^{\pm 1}]]$ is a field, if $a(z)v \in V(\!(z)\!)$ for all $v \in V$.
We shall write the field $a(z)$ as
\begin{align}
  a(z) = \sum_{n \in \Z} a_{(n)} z^{-n-1}.
\end{align}
The complex vector space of all fields on $V$ in the variable $z$ we will be denote by $\eus{F}(V)(z)$.
\medskip

\definition{A vertex algebra consists of the following data:
\begin{enumerate}[topsep=3pt,itemsep=0pt]
  \item[1)] a complex vector space $\mcal{V}$ (the space of states);
  \item[2)] a vector $\vac \in \mcal{V}$ (the vacuum vector);
  \item[3)] an endomorphism $T \colon \mcal{V} \rarr \mcal{V}$ (the translation operator);
  \item[4)] a linear mapping $Y(\,\cdot\,,z) \colon \mcal{V} \rarr \End \mcal{V}[[z^{\pm 1}]]$ sending
  \begin{align*}
    a \in \mcal{V} \mapsto Y(a,z) = \sum_{n \in \Z} a_{(n)} z^{-n-1} \in \eus{F}(\mcal{V})(z)
  \end{align*}
  (the state-field correspondence)
\end{enumerate}
satisfying the subsequent axioms:
\begin{enumerate}[topsep=3pt,itemsep=0pt]
  \item[1)] $Y(\vac,z)= \id_\mcal{V}$, $Y(a,z)\vac|_{z=0}=a$ (the vacuum axiom);
  \item[2)] $T\vac=0$, $[T,Y(a,z)] = \partial_z Y(a,z)$ (the translation axiom);
  \item[3)] for all $a,b \in \mcal{V}$, there is a non-negative integer $N_{a,b} \in \N_0$ such that
  \begin{align*}
    (z-w)^{N_{a,b}}[Y(a,z),Y(b,w)]=0
  \end{align*}
  (the locality axiom).
\end{enumerate}

A vertex algebra $\mcal{V}$ is called $\Z$-graded if $\mcal{V}$ is a $\Z$-graded vector space, $\vac$ is a vector of degree $0$, $T$ is an endomorphism of degree $1$, and for $a \in \mcal{V}_m$ the field $Y(a,z)$ has conformal dimension $m$, i.e.\ we have
\begin{align*}
  \deg a_{(n)} = -n+m-1
\end{align*}
for all $n \in \Z$.}

Let us note that according to the translation axiom, the action of $T$ on the space of states $\mcal{V}$ is completely determined by $Y$, since we have $T(a) = a_{(-2)}\vac$. Moreover, we have $a = a_{(-1)}\vac$.

%%%%%%%%%%%%%%%%%%%%%%%%%%%%%%%%%%%%%%%%%%%%%%%%%%%%%%%%%%%%%%%%%%%%%%%%%%%%%%%%%%%%%%%%%%

\subsection{Positive energy representations}

Let us consider a $\Z$-graded vertex algebra $\mcal{V}$. Then a $\mcal{V}$-module $M$ is called graded if $M$ is a $\C$-graded vector space and for $a \in \mcal{V}_m$ the field $Y_M(a,z)$ has conformal dimension $m$, i.e.\ we have
\begin{align*}
  \deg a_{(n)}^M = -n+m-1
\end{align*}
for all $n \in \Z$. Let us note that by shifting a given gradation on $M$ by a complex number we obtain a new gradation on $M$.
\medskip

\definition{Let $\mcal{V}$ be a $\Z$-graded vertex algebra. We say that a graded $\mcal{V}$-module $M$ is a \emph{positive energy $\mcal{V}$-module} if $M=\bigoplus_{n=0}^\infty M_{\lambda+n}$ and $M_\lambda \neq 0$, where $\lambda \in \C$. Moreover, we denote by $M_{\rm top}$ the top degree component $M_\lambda$ of $M$. The category of positive energy $\mcal{V}$-modules we will denote by $\mcal{E}_+\!(\mcal{V})$.}

In \cite{Zhu1996}, Zhu introduced a functorial construction which assigns to a $\Z$-graded vertex algebra an associative algebra known as the Zhu's algebra. Let $\mcal{V}$ be a $\Z$-graded vertex algebra. We define a bilinear mapping on $\mcal{V}$ by
\begin{align}
  a * b = \Res_{z=0}\!\bigg(Y(a,z){(1+z)^{\deg a} \over z}\,b\bigg) = \sum_{i=0}^{\deg a} \binom{\deg a}{i} a_{(i-1)}b \label{eq:Zhu algebra mult}
\end{align}
for homogeneous elements $a, b \in \mcal{V}$ and extend linearly. The Zhu's algebra $A(\mcal{V})$ is defined as
\begin{align}
  A(\mcal{V}) = \mcal{V}/O(\mcal{V}),
\end{align}
where $O(\mcal{V})$ is the vector subspace of $\mcal{V}$ spanned by
\begin{align}
  \Res_{z=0}\!\bigg(Y(a,z){(1+z)^{\deg a} \over z^2}\,b\bigg) = \sum_{i=0}^{\deg a} \binom{\deg a}{i} a_{(i-2)}b
\end{align}
for homogeneous elements $a,b \in \mcal{V}$. We denote by $\pi_{\rm Zhu}$ the canonical projection from $\mcal{V}$ to $A(\mcal{V})$. The bilinear mapping \eqref{eq:Zhu algebra mult} induces an associative multiplication on the quotient $A(\mcal{V})$. Further, we define
\begin{align*}
  o(a) =a_{(\deg a -1)}
\end{align*}
for a homogeneous element $a \in \mcal{V}$. Then it easily follows that for a homogeneous element $a \in \mcal{V}$ the operator $o_M(a) = a^M_{(\deg a -1)}$ preserves the homogeneous components of any graded $\mcal{V}$-module $M$.

As the following theorem proved in \cite{Zhu1996} shows, the Zhu's algebra $A(\mcal{V})$ plays a prominent role in the representation theory of vertex algebras.
\medskip

\theorem{\label{thm:Zhu correspondence}
Let $\mcal{V}$ be a $\Z$-graded vertex algebra and let $M$ be a positive energy $\mcal{V}$-module. Then the top degree component $M_{\rm top}$ is an $A(\mcal{V})$-module, where the action of $\pi_{\rm Zhu}(a) \in A(\mcal{V})$ for $a \in \mcal{V}$ is given by $o_M(a)$. In addition, the correspondence $M \mapsto M_{\rm top}$ gives a bijection between the set of simple positive energy $\mcal{V}$-modules and that of simple $A(\mcal{V})$-modules.}

To a $\Z$-graded vertex algebra $\mcal{V}$ we may associate a complete topological Lie algebra $U(\mcal{V})$, first introduced by Borcherds \cite{Borcherds1986}, by
\begin{align*}
  U(\mcal{V}) = (\mcal{V} \otimes_\C \C(\!(t)\!)) / \im \partial,
\end{align*}
where
\begin{align*}
  \partial = T \otimes \id + \id \otimes \partial_t.
\end{align*}
If we denote by $a_{[n]}$ for $a \in \mcal{V}$ and $n \in \Z$ the projection of $a \otimes t^n \in \mcal{V} \otimes_\C \C(\!(t)\!)$ onto $U(\mcal{V})$, then the Lie bracket on $U(\mcal{V})$ is determined by
\begin{align*}
  [a_{[m]}, b_{[n]}] = \sum_{k=0}^\infty \binom{m}{k} (a_{(k)}b)_{[m+n-k]}
\end{align*}
for $a, b \in \mcal{V}$ and $m,n \in \Z$.

Further, for a homogeneous element $a \in \mcal{V}$, we set $\deg a_{[n]} = -n+\deg a -1$. Then the degree assignment to elements of $U(\mcal{V})$ gives us a triangular decomposition
\begin{align}
  U(\mcal{V}) = U(\mcal{V})_- \oplus U(\mcal{V})_0 \oplus U(\mcal{V})_+ \label{eq:triangular docomposition U(V)}
\end{align}
of the Lie algebra $U(\mcal{V})$ together with a canonical surjective homomorphism
\begin{align*}
  U(\mcal{V})_0  \rarr A(\mcal{V})
\end{align*}
of Lie algebras defined by
\begin{align*}
  a_{[\deg a -1 ]} \mapsto \pi_{\rm Zhu}(a)
\end{align*}
for a homogeneous element $a \in \mcal{V}$.

Let us consider a $\mcal{V}$-module $M$. Then it has also a natural structure of a $U(\mcal{V})$-module since we have a canonical homomorphism
\begin{align*}
  U(\mcal{V}) \rarr \End \mcal{V}
\end{align*}
of Lie algebras defined through
\begin{align*}
  a_{[n]} \mapsto a_{(n)}
\end{align*}
for $a \in \mcal{V}$ and $n \in \Z$. We denote by $\Omega_\mcal{V}(M)$ the vector subspace of $M$ consisting of lowest weight vectors, i.e.\ we have
\begin{align}
  \Omega_\mcal{V}(M) = \{v\in M;\, U(\mcal{V})_-v=0\}.
\end{align}
It follows immediately using the triangular decomposition of $U(\mcal{V})$ that $\Omega_\mcal{V}(M)$ is a $U(\mcal{V})_0$-module. Moreover, by \cite{Dong-Li-Mason1998} we have that $\Omega_\mcal{V}(M)$ is an $A(\mcal{V})$-module, where the action of $\pi_{\rm Zhu}(a) \in A(\mcal{V})$ for $a \in \mcal{V}$ is given by $o_M(a)$. It is clear that
\begin{align}
  \Omega_\mcal{V} \colon \mcal{E}(\mcal{V}) \rarr \mcal{M}(A(\mcal{V}))
\end{align}
is a functor. Let us note that if $M$ is a positive energy $\mcal{V}$-module, then we have $\Omega_\mcal{V}(M) \supset M_{\rm top}$ and $\Omega_\mcal{V}(M) = M_{\rm top}$ provided $M$ is a simple $\mcal{V}$-module.

Therefore, we may consider a functor $\Omega_\mcal{V} \colon \mcal{E}_+\!(\mcal{V}) \rarr \mcal{M}(A(\mcal{V}))$. On the other hand, there exists also an induction functor
\begin{align}
  \mathbb{M}_\mcal{V} \colon \mcal{M}(A(\mcal{V})) \rarr \mcal{E}_+\!(\mcal{V})
\end{align}
which is a left adjoint functor to $\Omega_\mcal{V}$ and has the following universal property. For a $\mcal{V}$-module $M$ and a morphism $\varphi \colon  E \rarr \Omega_\mcal{V}(M)$ of $A(\mcal{V})$-modules, there is a unique morphism $\widetilde{\varphi} \colon \mathbb{M}_\mcal{V}(E) \rarr M$ of $\mcal{V}$-modules which extends $\varphi$, see \cite{Dong-Li-Mason1998}. Moreover, for an $A(\mcal{V})$-module $E$ we have $\mathbb{M}_\mcal{V}(E)_{\rm top} \simeq E$ as modules over $A(\mcal{V})$. Besides, since the $\mcal{V}$-module $\mathbb{M}_\mcal{V}(E)$ has a unique maximal $\mcal{V}$-submodule $\mathbb{K}_\mcal{V}(E)$ having zero intersection with the $A(\mcal{V})$-submodule $E$ of $\mathbb{M}_\mcal{V}(E)$, we may set
\begin{align}
  \mathbb{L}_\mcal{V}(E) = \mathbb{M}_\mcal{V}(E)/\mathbb{K}_\mcal{V}(E)
\end{align}
for an $A(\mcal{V})$-module $E$.

%%%%%%%%%%%%%%%%%%%%%%%%%%%%%%%%%%%%%%%%%%%%%%%%%%%%%%%%%%%%%%%%%%%%%%%%%%%%%%%%%%%%%%%%%%

\subsection{Affine vertex algebras}
\label{subsec:affine vertex}

Let $\mfrak{g}$ be a semisimple (reductive) finite-dimensional Lie algebra and $\kappa$ be a $\mfrak{g}$-invariant symmetric bilinear form on $\mfrak{g}$. The induced $\widehat{\mfrak{g}}_\kappa$-module $\mathbb{M}_{\kappa,\mfrak{g}}(\C)$, where $\C$ is the trivial $1$-dimensional $\mfrak{g}$-module, is of a special importance in the theory of vertex algebras since it is equipped with the natural structure of an $\N_0$-graded vertex algebra, called the \emph{universal affine vertex algebra}, see \cite{Kac1998}, which we will denote by $\mcal{V}_\kappa(\mfrak{g})$. For an element $a \in \mfrak{g}$, we denote by $a(z) \in \widehat{\mfrak{g}}_\kappa[[z^{\pm 1}]]$ the formal distribution defined by
\begin{align}
a(z) = \sum_{n \in \Z} a_n z^{-n-1}.  \label{eq:field a}
\end{align}
By using this formal power series we may rewrite the commutation relations \eqref{eq:commutation relation modes} for $\widehat{\mfrak{g}}_\kappa$ into the form
\begin{align}
  [a(z),b(w)] = [a,b](w)\delta(z-w) + \kappa(a,b)c\partial_w \delta(z-w) \label{eq:comm relations vertex algebra}
\end{align}
for $a, b \in \mfrak{g}$. The state field correspondence $Y \colon \mcal{V}_\kappa(\mfrak{g}) \rarr \End \mcal{V}_\kappa(\mfrak{g})[[z^{\pm 1}]]$ is given by
\begin{align*}
  Y(a_{1,-n_1-1} \dots a_{k,-n_k-1}\vac,z) = {1 \over n_1! \dots n_k!}\, \normOrd{\partial_z^{n_1} a_1(z) \dots \partial_z^{n_k}a_k(z)}
\end{align*}
for $k \in \N$, $n_1,n_2,\dots,n_k \in \N_0$ and $a_1,a_2,\dots,a_k \in \mfrak{g}$, where $\vac \in \mcal{V}_\kappa(\mfrak{g})$ is the vacuum vector (a highest weight vector of $\mathbb{M}_{\kappa,\mfrak{g}}(\C)$). The translation operator $T \colon \mcal{V}_\kappa(\mfrak{g}) \rarr \mcal{V}_\kappa(\mfrak{g})$ is defined by $T\vac = 0$ and $[T,a_n] = -na_{n-1}$ for $a \in \mfrak{g}$ and $n \in \Z$.
\medskip

To describe positive energy representations of $\mcal{V}_\kappa(\mfrak{g})$, we need to know its Zhu's algebra. It is easy to see that for $\mcal{V}_\kappa(\mfrak{g})$ we have a canonical isomorphism
\begin{align}
  A(\mcal{V}_\kappa(\mfrak{g})) \simeq U(\mfrak{g})
\end{align}
of associative algebras determined by
\begin{align}
  a_{1,-n_1-1}a_{2,-n_2-1}\dots a_{k,-n_k-1}\vac \mapsto
  (-1)^{n_1+n_2+\dots +n_k} a_k\dots a_2a_1
\end{align}
for $k \in \N$, $n_1,n_2,\dots,n_k \in \N_0$ and $a_1,a_2,\dots,a_k \in \mfrak{g}$.

Let us note that for the universal affine vertex algebra $\mcal{V}_\kappa(\mfrak{g})$ the functors $\mathbb{M}_{\mcal{V}_\kappa(\mfrak{g})}$ and $\mathbb{L}_{\mcal{V}_\kappa(\mfrak{g})}$ coincide with the functors $\mathbb{M}_{\kappa,\mfrak{g}}$ and $\mathbb{L}_{\kappa,\mfrak{g}}$, respectively. Hence, according to a theorem of Zhu the assignment $E \mapsto \mathbb{L}_{\kappa,\mfrak{g}}(E)$ gives a one-to-one correspondence between the isomorphism classes of simple $\mfrak{g}$-modules and simple positive energy $\mcal{V}_\kappa(\mfrak{g})$-modules. Therefore, the study of positive energy $\mcal{V}_\kappa(\mfrak{g})$-modules reduces to the study of $\mfrak{g}$-modules.
\medskip

In addition, the unique simple quotient $\mathbb{L}_{\kappa,\mfrak{g}}(\C)$ of $\mathbb{M}_{\kappa,\mfrak{g}}(\C)$ has also the natural structure of an $\N_0$-graded vertex algebra, called the \emph{simple affine vertex algebra}, which we will denote by $\mcal{L}_\kappa(\mfrak{g})$. The Zhu's algebra $A(\mcal{L}_\kappa(\mfrak{g}))$ is a homomorphic image of $U(\mfrak{g})$, hence we have
\begin{align}
  A(\mcal{L}_\kappa(\mfrak{g})) \simeq U(\mfrak{g})/I_\kappa
\end{align}
for some two-sided ideal $I_\kappa$ of $U(\mfrak{g})$. Moreover, the assignment $E \mapsto \mathbb{L}_{\kappa,\mfrak{g}}(E)$ gives a one-to-one correspondence between isomorphism classes of simple modules over $U(\mfrak{g})/I_\kappa$ and simple positive energy $\mcal{L}_\kappa(\mfrak{g})$-modules.

%%%%%%%%%%%%%%%%%%%%%%%%%%%%%%%%%%%%%%%%%%%%%%%%%%%%%%%%%%%%%%%%%%%%%%%%%%%%%%%%%%%%%%%%%%

\subsection{Weyl vertex algebras}
\label{subsec:Weyl vertex}

Let $V$ be a finite-dimensional vector space. Then the induced $\eus{A}_{\mcal{K}(V)}$-module $\mathbb{M}_{\mcal{K}(V)}(S(V^*))$ carries the natural structure of an $\N_0$-graded vertex algebra, called the \emph{Weyl vertex algebra}, which we will denote by $\mcal{M}_V$. Let $\{x_1,x_2,\dots,x_m\}$, where $\dim V = m$, be linear coordinate functions on $V$ and let $a_i(z), a_i^*(z) \in \eus{A}_{\mcal{K}(V)}[[z^{\pm 1}]]$ for $i=1,2,\dots,m$ be the formal distributions defined by
\begin{align}
  a_i(z)= \sum_{n\in \Z} a_{i,n}z^{-n-1} \qquad \text{and} \qquad a_i^*(z) = \sum_{n\in \Z} a^*_{i,n}z^{-n}, \label{eq:field a,a*}
\end{align}
where $a_{i,n}=\partial_{x_{i,n}}$ and $a^*_{i,n}=x_{i,-n}$ for $n \in \Z$ and $i=1,2,\dots,m$. By using these formal power series, the canonical commutation relations for $\eus{A}_{\mcal{K}(V)}$ we may write in the form
\begin{align}
  [a_i(z), a_j(w)] = 0, \qquad [a_i(z), a_j^*(w)] = \delta_{i,j} \delta(z-w), \qquad [a_i^*(z), a_j^*(w)] = 0
\end{align}
for $i,j = 1,2,\dots,m$. The state field correspondence $Y \colon \mcal{M}_V \rarr \End \mcal{M}_V[[z^{\pm 1}]]$ is given by
\begin{multline*}
  Y(a_{i_1,-n_1-1}\dots a_{i_r,-n_r-1}a^*_{j_1,-m_1}\dots a^*_{j_s,-m_s}\vac,z) = \\ {1 \over n_1! \dots n_r!} {1 \over m_1!\dots m_s!} \, \normOrd{\partial_z^{n_1}a_{i_1}(z) \dots \partial_z^{n_r}a_{i_r}(z) \partial_z^{m_1}a_{j_1}^*(z) \dots \partial_z^{m_s}a_{j_s}^*(z)}
\end{multline*}
for $r,s \in \N$, $n_1,\dots,n_r,m_1,\dots,m_s \in \N_0$, where $\vac \in \mcal{M}_V$ is the vacuum vector (a highest weight vector of $\mathbb{M}_{\mcal{K}(V)}(S(V^*))$). The translation operator $T \colon \mcal{M}_V \rarr \mcal{M}_V$ is defined by $T\vac = 0$, $[T,a_{i,n}] = -na_{i,n-1}$ and $[T,a^*_{i,n}] = -(n-1)a^*_{i,n-1}$ for $n \in \Z$ and $i=1,2,\dots,m$. Moreover, we have a canonical isomorphism
\begin{align}
  A(\mcal{M}_V) \simeq \eus{A}_V
\end{align}
of associative algebras determined by
\begin{multline}
  a_{i_1,-n_1-1}\dots a_{i_r,-n_r-1}a^*_{j_1,-m_1}\dots a^*_{j_s,-m_s}\vac \mapsto \\
  \delta_{m_1,0} \dots \delta_{m_s,0}(-1)^{n_1+\dots + n_r}x_{j_1}\dots x_{j_s}\partial_{x_{i_1}}\dots\partial_{x_{i_r}}
\end{multline}
for $r,s \in \N$, $n_1,\dots,n_r,m_1,\dots,m_s \in \N_0$.

Let us note that for the Weyl vertex algebra $\mcal{M}_V$ the functors $\mathbb{M}_{\mcal{M}_V}$ and $\mathbb{L}_{\mcal{M}_V}$ coincide with the functors $\mathbb{M}_{\mcal{K}(V)}$ and $\mathbb{L}_{\mcal{K}(V)}$, respectively.

%%%%%%%%%%%%%%%%%%%%%%%%%%%%%%%%%%%%%%%%%%%%%%%%%%%%%%%%%%%%%%%%%%%%%%%%%%%%%%%%%%%%%%%%%%
%%%%%%%%%%%%%%%%%%%%%%%%%%%%%%%%%%%%%%%%%%%%%%%%%%%%%%%%%%%%%%%%%%%%%%%%%%%%%%%%%%%%%%%%%%

\section{Twisting functors and Gelfand--Tsetlin modules}
\label{sec:GT modules}

We recall the definition of Gelfand--Tsetlin modules for a complex semisimple finite-dimensional Lie algebra $\mfrak{g}$ and describe the construction of the twisting functor $T_\alpha$ assigned to a positive root $\alpha$ of $\mfrak{g}$.

%%%%%%%%%%%%%%%%%%%%%%%%%%%%%%%%%%%%%%%%%%%%%%%%%%%%%%%%%%%%%%%%%%%%%%%%%%%%%%%%%%%%%%%%%%

\subsection{Twisting functors for semisimple Lie algebras}

Let us consider a complex semisimple Lie algebra $\mfrak{g}$. For a positive root $\alpha \in \Delta_+$ of $\mfrak{g}$ we denote by $\mfrak{s}_\alpha$ the Lie subalgebra of $\mfrak{g}$ generated by the $\mfrak{sl}_2$-triple $(e_\alpha,h_\alpha,f_\alpha)$, where $e_\alpha \in \mfrak{g}_\alpha$ and $f_\alpha \in \mfrak{g}_{-\alpha}$ satisfy $[e_\alpha,f_\alpha]=h_\alpha$.
Further, we define the Lie subalgebras $\mfrak{s}_\alpha^+= \mfrak{s}_\alpha \cap \mfrak{n} = \C e_\alpha$ and $\mfrak{s}_\alpha^-= \mfrak{s}_\alpha \cap \widebar{\mfrak{n}} = \C f_\alpha$ of $\mfrak{s}_\alpha$. Let $c_\alpha$ be the quadratic Casimir element given through
\begin{align*}
   c_\alpha = e_\alpha f_\alpha + f_\alpha e_\alpha + {\textstyle {1 \over 2}} h_\alpha^2,
\end{align*}
which is a free generator of the center $Z(\mfrak{s}_\alpha)$ of $U(\mfrak{s}_\alpha)$.
Then we denote by $\Gamma_\alpha$ the commutative subalgebra of $U(\mfrak{g})$ generated by the Cartan subalgebra $\mfrak{h}$ and by the center $Z(\mfrak{s}_\alpha)$. Therefore, we have the category $\mcal{H}(\mfrak{g},\Gamma_\alpha)$ of $\Gamma_\alpha$-Gelfand--Tsetlin modules, or simply $\alpha$-Gelfand--Tsetlin modules, which was studied in \cite{Futorny-Krizka2019} and \cite{Futorny-Krizka2019b}. Let us note that $\mcal{H}(\mfrak{g}, \Gamma_\alpha)$ contains $\mcal{I}(\mfrak{g},\mfrak{s}_\alpha^+)$ and $\mcal{I}(\mfrak{g},\mfrak{s}_\alpha^-)$ as full subcategories.

Since the multiplicative set $\{f_\alpha^n;\, n \in \N_0\}$ in the universal enveloping algebra $U(\mfrak{g})$ is a left (right) denominator set, thanks to the fact that $f_\alpha$ is a locally $\ad$-nilpotent regular element in $U(\mfrak{g})$, we define the \emph{twisting functor}
\begin{align*}
  T_\alpha=T_{f_\alpha} \colon \mcal{M}(\mfrak{g}) \rarr \mcal{M}(\mfrak{g})
\end{align*}
by
\begin{align*}
  T_\alpha(M) = (U(\mfrak{g})_{(f_\alpha)}/U(\mfrak{g})) \otimes_{U(\mfrak{g})} M  \simeq M_{(f_\alpha)}/M
\end{align*}
for $M \in \mcal{M}(\mfrak{g})$. Let us note that the quotient $M_{(f_\alpha)}/M$ means the quotient of $M_{(f_\alpha)}$ by the image of the canonical homomorphism $M \rarr M_{(f_\alpha)}$ of $\mfrak{g}$-modules. As $f_\alpha$ is a locally $\ad$-nilpotent regular element in $U(\mfrak{g})$, the functor $T_\alpha$ is right exact. The twisting functor $T_\alpha$ for a simple root $\alpha \in \Pi$ is well studied (see e.g.\ \cite{Andersen-Stroppel2003}). In this case the functor $T_\alpha$ preserves the category $\mcal{O}(\mfrak{g})$ up to a conjugation of the action of $\mfrak{g}$.
\medskip

\theorem{\cite[Theorem 3.3]{Futorny-Krizka2019}\label{thm:finitely generated finite}
For $\alpha \in \Delta_+$ the functor $T_\alpha$ induces the restricted functor
\begin{align*}
  T_\alpha \colon \mcal{I}_f(\mfrak{g},\mfrak{s}_\alpha^+) \rarr \mcal{I}_f(\mfrak{g},\mfrak{s}_\alpha^-),
\end{align*}
where $\mcal{I}_f(\mfrak{g},\mfrak{s}_\alpha^\pm)$ is the category of finitely generated locally $\mfrak{s}_\alpha^\pm$-finite weight $\mfrak{g}$-modules.}

It was shown in \cite{Futorny-Krizka2019b} how the twisting functor $T_\alpha$ can be used to construct $\mfrak{g}$-modules in $\mcal{H}(\mfrak{g},\Gamma_\alpha)$ with finite $\Gamma_\alpha$-multiplicities. Namely, we have the following theorem.
\medskip

\theorem{\cite[Theorem 3.4]{Futorny-Krizka2019}\label{thm:homol finite}
Let $\alpha \in \Delta_+$.
\begin{itemize}[topsep=3pt,itemsep=0pt]
\item[i)] For a weight  $\mfrak{g}$-module $M$, the $\mfrak{g}$-module $T_\alpha(M)$ belongs to $\mcal{H}(\mfrak{g},\Gamma_\alpha)$. In addition, $T_\alpha(M)$ has finite $\Gamma_\alpha$-multiplicities if and only if the first cohomology group $H^1(\mfrak{s}_\alpha^-;M)$ is a weight $\mfrak{h}$-module with finite-dimensional weight spaces.
\item[ii)] If $M$ is a highest weight $\mfrak{g}$-module, then $T_\alpha(M)$ is a cyclic weight $\mfrak{g}$-module with finite $\Gamma_\alpha$-multiplicities.
\end{itemize}
}

For a $\mfrak{g}$-module $M$, we denote by $\Ann_{U(\mfrak{g})}\!M$ the annihilator of $M$ in $U(\mfrak{g})$. Then we get the following important proposition.
\medskip

\proposition{\label{prop-annihilator}
Let $M$ be a $U(\mfrak{s}_\alpha^-)$-free $\mfrak{g}$-module for $\alpha \in \Delta_+$. Then we have $\Ann_{U(\mfrak{g})}\!T_\alpha(M) = \Ann_{U(\mfrak{g})}\!M$.}

\proof{First, we observe that $\Ann_{U(\mfrak{g})}\! M_{(f_\alpha)}=\Ann_{U(\mfrak{g})}\!M$. Indeed, the inclusion $M\subset M_{(f_\alpha)}$ gives us
\begin{align*}
 \Ann_{U(\mfrak{g})}\! M_{(f_\alpha)}\subset \Ann_{U(\mfrak{g})}\!M.
\end{align*}
The opposite embedding is obvious since $f_\alpha$ is a locally ad-nilpotent element in $U(\mfrak{g})$. Further, from $T_{\alpha}(M)=M_{(f_\alpha)}/M$ we get that $\Ann_{U(\mfrak{g})}\!M \subset \Ann_{U(\mfrak{g})}\!T_{\alpha}(M)$. Now, we prove the opposite inclusion. Let us consider a vector $v \in M$ and let us assume that $a \in \Ann_{U(\mfrak{g})}\!T_{\alpha}(M)$. Then we have $af_\alpha^{-n}v \in M$ for all $n \in \N_0$. By using the fact $f_\alpha^naf_\alpha^{-n-m}v \in M$ and the formula
\begin{align*}
  f_\alpha^na f_\alpha^{-n-m}v = f_\alpha^{-m} \sum_{k=0}^\infty \binom{n+m+k-1}{k} f_\alpha^{-k} \ad(f_\alpha)^k(a)v
\end{align*}
for all $n,m \in \N_0$, we obtain immediately $f_\alpha^{-k-m}\ad(f_\alpha)^k(a)v \in M$ for all $k,m \in \N_0$ which gives us $\ad(f_\alpha)^k(a)v=0$ for all $k\in \N_0$. We may write
\begin{align*}
  af_\alpha^{-n}v = f_\alpha^{-n}\sum_{k=0}^\infty \binom{n+k-1}{k} f_\alpha^{-k} \ad(f_\alpha)^k(a)v = 0
\end{align*}
for all $n \in \N_0$. Therefore, we have $a \in \Ann_{U(\mfrak{g})}\! M_{(f_\alpha)} = \Ann_{U(\mfrak{g})}\! M$ and we are done.}

If $M$ is a $\mfrak{g}$-module, then $T_\alpha(M)$ is a locally $\mfrak{s}_\alpha^-$-finite $\mfrak{g}$-module for any $\alpha \in \Delta_+$. On the other hand, if $M$ is a locally $\mfrak{s}_\alpha^-$-finite $\mfrak{g}$-module, then $T_\alpha(M)=0$.
Let us also note that for any $\alpha \in \Delta_+$ the twisting functor $T_\alpha$ commutes with the translation functors by \cite[Theorem 3.6]{Futorny-Krizka2019b}.
\medskip

Let $\mfrak{p}=\mfrak{l} \oplus \mfrak{u}$ be the standard parabolic subalgebra of $\mfrak{g}$ associated to a subset $\Sigma$ of $\Pi$ with the nilradical $\mfrak{u}$, the opposite nilradical $\widebar{\mfrak{u}}$ and the Levi subalgebra $\mfrak{l}$. We define the subsets
\begin{align*}
  \Delta_+^\mfrak{u} = \{\alpha \in \Delta_+;\, \mfrak{g}_\alpha \subset \mfrak{u}\}, \qquad \Delta_+^\mfrak{l} = \{\alpha \in \Delta_+;\, \mfrak{g}_\alpha \subset \mfrak{l}\}
\end{align*}
of $\Delta_+$ and we set
\begin{align*}
  \Lambda^+(\mfrak{p}) = \{\lambda \in \mfrak{h}^*;\, (\forall \alpha \in \Sigma)\, \lambda(h_\alpha) \in \N_0\}.
\end{align*}
For a weight $\lambda \in \Lambda^+(\mfrak{p})$, we denote by $\mathbb{F}_\lambda$ the simple finite-dimensional $\mfrak{p}$-module with highest weight $\lambda$ and by $M^\mfrak{g}_\mfrak{p}(\lambda)$ the generalized Verma $\mfrak{g}$-module with highest weight $\lambda$ defined by
\begin{align*}
  M^\mfrak{g}_\mfrak{p}(\lambda) = U(\mfrak{g}) \otimes_{U(\mfrak{p})} \!\mathbb{F}_\lambda.
\end{align*}
Then for $\alpha \in \Delta_+^\mfrak{u}$ we define the $\alpha$-Gelfand--Tsetlin $\mfrak{g}$-module $W^\mfrak{g}_\mfrak{p}(\lambda,\alpha)$ by
\begin{align*}
  W^\mfrak{g}_\mfrak{p}(\lambda,\alpha) = T_\alpha(M^\mfrak{g}_\mfrak{p}(\lambda)).
\end{align*}
As a consequence of the Poincar\'e--Birkhoff--Witt theorem we have isomorphisms
\begin{align*}
  M^\mfrak{g}_\mfrak{p}(\lambda) \simeq U(\widebar{\mfrak{u}}) \otimes_\C \mathbb{F}_\lambda \qquad \text{and} \qquad W^\mfrak{g}_\mfrak{p}(\lambda,\alpha) \simeq (U(\widebar{\mfrak{u}})_{(f_\alpha)}/U(\widebar{\mfrak{u}})) \otimes_\C \mathbb{F}_\lambda
\end{align*}
of $U(\widebar{\mfrak{u}})$-modules. Moreover, since the twisting functor $T_\alpha$ is right exact, by applying of $T_\alpha$ on a surjective homomorphism
\begin{align*}
  M^\mfrak{g}_\mfrak{b}(\lambda) \rarr M^\mfrak{g}_\mfrak{p}(\lambda)
\end{align*}
of generalized Verma modules, we get a surjective homomorphism
\begin{align*}
  W^\mfrak{g}_\mfrak{b}(\lambda,\alpha) \rarr W^\mfrak{g}_\mfrak{p}(\lambda,\alpha)
\end{align*}
of $\alpha$-Gelfand--Tsetlin modules.
\medskip

By applying Theorem \ref{thm:homol finite} and Proposition \ref{prop-annihilator} on $W^\mfrak{g}_\mfrak{p}(\lambda,\alpha)$ for $\lambda \in \Lambda^+(\mfrak{p})$ and $\alpha \in \Delta_+^\mfrak{u}$ we obtain the following immediate corollaries.
\medskip

\corollary{Let $\lambda \in \Lambda^+(\mfrak{p})$ and $\alpha \in \Delta_+^\mfrak{u}$. Then the $\mfrak{g}$-module $W^\mfrak{g}_\mfrak{p}(\lambda,\alpha)\in \mcal{H}(\mfrak{g}, \Gamma_\alpha)$ is cyclic weight with central character and finite $\Gamma_\alpha$-multiplicities. Moreover, $W^\mfrak{g}_\mfrak{p}(\lambda,\alpha)$ belongs to $\mcal{H}(\mfrak{g}, \Gamma)$ and has finite $\Gamma$-multiplicities for any commutative subalgebra $\Gamma$ of $U(\mfrak{g})$ containing $\Gamma_\alpha$.}

\corollary{\label{cor-annihilator}
Let $\lambda \in \Lambda^+(\mfrak{p})$ and $\alpha \in \Delta_+^\mfrak{u}$. Then we have
$\Ann_{U(\mfrak{g})}\! W^\mfrak{g}_\mfrak{p}(\lambda,\alpha) = \Ann_{U(\mfrak{g})}\! M^\mfrak{g}_\mfrak{p}(\lambda)$.}

\vspace{-2mm}

%%%%%%%%%%%%%%%%%%%%%%%%%%%%%%%%%%%%%%%%%%%%%%%%%%%%%%%%%%%%%%%%%%%%%%%%%%%%%%%%%%%%%%%%%%

\subsection{Affine $\alpha$-Gelfand--Tsetlin modules}
\label{subsec:affine Gelfand--Tsetlin modules}

Let $\widetilde{\mfrak{g}}_\kappa$ be the extended affine Kac--Moody algebra associated to a complex semisimple Lie algebra $\mfrak{g}$ of level $\kappa$. Then for a commutative subalgebra $\Gamma$ of the completed universal enveloping algebra $\smash{\widetilde{U}_c(\widetilde{\mfrak{g}}_\kappa)}$ of $\widetilde{\mfrak{g}}_\kappa$, we denote by $\mcal{H}(\widetilde{\mfrak{g}}_\kappa,\Gamma)$ the full subcategory of $\mcal{E}(\widetilde{\mfrak{g}}_\kappa)$ consisting of smooth $\Gamma$-weight $\widetilde{\mfrak{g}}_\kappa$-modules on which the central element $c$ acts as the identity. Let us note that $\mcal{H}(\widetilde{\mfrak{g}}_\kappa,\Gamma)$ is closed with respect to the operations of taking submodules and quotients. In addition, if $\Gamma$ contains the Cartan subalgebra $\smash{\widetilde{\mfrak{h}}}$, the $\widetilde{\mfrak{g}}_\kappa$-modules from $\mcal{H}(\smash{\widetilde{\mfrak{g}}}_\kappa, \Gamma)$ are called \emph{$\Gamma$-Gelfand--Tsetlin $\widetilde{\mfrak{g}}_\kappa$-modules}.
\medskip

For a real positive root $\alpha \in \smash{\widehat{\Delta}^{\rm re}_+}$, we denote by $\mfrak{s}_\alpha$ the Lie subalgebra of $\widetilde{\mfrak{g}}_\kappa$ generated by the $\mfrak{sl}_2$-triple $(e_\alpha,h_\alpha,f_\alpha)$, where
\begin{align*}
  e_\alpha = e_{\gamma,n}, \qquad h_\alpha = h_{\gamma,0} + n\kappa(e_\gamma,f_\gamma)c, \qquad  f_\alpha = f_{\gamma,-n}
\end{align*}
if $\alpha = \gamma+n\delta$ with $\gamma \in \Delta_+$, $n \in \N_0$ and
\begin{align*}
  e_\alpha = f_{\gamma,n}, \qquad h_\alpha = -h_{\gamma,0} + n\kappa(e_\gamma,f_\gamma)c, \qquad  f_\alpha = e_{\gamma,-n}
\end{align*}
provided $\alpha = -\gamma+n\delta$ with $\gamma \in \Delta_+$, $n \in \N$. The quadratic Casimir element $c_\alpha$ given by
\begin{align*}
   c_\alpha = e_\alpha f_\alpha + f_\alpha e_\alpha + {\textstyle {1 \over 2}} h_\alpha^2
\end{align*}
is a free generator of the center $Z(\mfrak{s}_\alpha)$ of $U(\mfrak{s}_\alpha)$. Furthermore, the commutative subalgebra of $\smash{\widetilde{U}_c(\widetilde{\mfrak{g}}_\kappa)}$ generated by the Cartan subalgebra $\smash{\widetilde{\mfrak{h}}}$ and by the center $Z(\mfrak{s}_\alpha)$ we denote by $\Gamma_\alpha$. We also define the Lie subalgebras $\mfrak{s}_\alpha^+= \mfrak{s}_\alpha \cap \widehat{\mfrak{n}}_{\rm st} = \C e_\alpha$ and $\mfrak{s}_\alpha^-= \mfrak{s}_\alpha \cap \smash{\widehat{\widebar{\mfrak{n}}}}_{\rm st} = \C f_\alpha$ of $\widetilde{\mfrak{g}}_\kappa$. The objects of $\mcal{H}(\widetilde{\mfrak{g}}_\kappa,\Gamma_\alpha)$ will be simply called \emph{$\alpha$-Gelfand--Tsetlin $\widetilde{\mfrak{g}}_\kappa$-modules}.
\medskip

Let $\Gamma$ be a commutative subalgebra of $\smash{\widetilde{U}_c(\widetilde{\mfrak{g}}_\kappa)}$ containing the Cartan subalgebra $\smash{\widetilde{\mfrak{h}}}$. Besides, let us assume that $[a,\mfrak{s}_\alpha^\pm] = \pm\chi(a) \mfrak{s}_\alpha^\pm$ for all $a \in \Gamma$ and some character $\chi \in \Gamma$. Then for a $\Gamma$-Gelfand--Tsetlin $\widetilde{\mfrak{g}}_\kappa$-module $M$, the Lie algebra cohomology groups $H^n(\mfrak{s}_\alpha^\pm;M)$ are $\Gamma$-weight modules for all $n \in \N_0$.
By general results we get $H^0(\mfrak{s}_\alpha^\pm;M)\simeq \smash{M^{\mfrak{s}_\alpha^\pm}}$, $H^1(\mfrak{s}_\alpha^\pm;M) \simeq M/\mfrak{s}_\alpha^\pm M$ and $H^n(\mfrak{s}_\alpha^\pm;M)=0$ for $n>1$. Moreover, we have that $Z(\mfrak{s}_\alpha)$ commutes with $\Gamma$, which enables us to extend $\Gamma$ by $Z(\mfrak{s}_\alpha)$. Therefore, we denote by $\Gamma_{\alpha,{\rm ext}}$ the commutative subalgebra of $\smash{\widetilde{U}_c(\widetilde{\mfrak{g}}_\kappa)}$ generated by $\Gamma$ and by $Z(\mfrak{s}_\alpha)$. As a mild and straightforward generalization of Theorem 2.3 in \cite{Futorny-Krizka2019b} we obtain the following important statement.
\medskip

\theorem{\label{thm:GT modules general}
Let $\alpha \in \smash{\widehat{\Delta}}^{\rm re}_+$ and let $M$ be a $\Gamma$-Gelfand--Tsetlin $\widetilde{\mfrak{g}}_\kappa$-module which is locally $\mfrak{s}_\alpha^-$-finite. Then $M$ is a $\Gamma_{\alpha,{\rm ext}}$-Gelfand--Tsetlin $\widetilde{\mfrak{g}}_\kappa$-module with finite $\Gamma_{\alpha,{\rm ext}}$-multiplicities if and only if the zeroth cohomology group $H^0(\mfrak{s}_\alpha^-;M)$ is a $\Gamma$-weight module with finite $\Gamma$-multiplicities.}

The following result is a consequence of Theorem \ref{thm:GT modules general} and the fact that the extension of $\smash{U_c(\widetilde{\mfrak{h}})}$ by $Z(\mfrak{s}_\alpha)$ is the commutative algebra $\Gamma_\alpha$.
\medskip

\corollary{\label{cor:highest weight GT modules}
Let $\alpha \in \smash{\widehat{\Delta}}^{\rm re}_+$ and let $M$ be a locally $\mfrak{s}_\alpha^-$-finite smooth weight $\widetilde{\mfrak{g}}_\kappa$-module on which the central element $c$ acts as the identity. Then $M$ is a $\Gamma_\alpha$-Gelfand--Tsetlin $\widetilde{\mfrak{g}}_\kappa$-module with finite $\Gamma_\alpha$-multiplicities if and only if the zeroth cohomology group $H^0(\mfrak{s}_\alpha^-;M)$ is a weight $\smash{\widetilde{\mfrak{h}}}$-module with finite-dimensional weight spaces.}

By duality we can analogously prove a similar statement if we replace $\mfrak{s}_\alpha^-$ by $\mfrak{s}_\alpha^+$ and $H^0(\mfrak{s}_\alpha^-;M)$ by $H^0(\mfrak{s}_\alpha^+;M)$.
\medskip

For a Lie subalgebra $\mfrak{a}$ of $\widetilde{\mfrak{g}}_\kappa$ we denote by $\mcal{I}(\widetilde{\mfrak{g}}_\kappa,\mfrak{a})$ and $\mcal{I}_f(\widetilde{\mfrak{g}}_\kappa,\mfrak{a})$ the full subcategories  of
$\mcal{E}(\widetilde{\mfrak{g}}_\kappa)$ consisting of locally $\mfrak{a}$-finite and finitely generated locally $\mfrak{a}$-finite, respectively, smooth weight $\widetilde{\mfrak{g}}_\kappa$-modules on which the central element $c$ acts as the identity. Therefore, for a real positive root $\alpha \in \smash{\widehat{\Delta}}^{\rm re}_+$ we have the full subcategories $\mcal{I}(\widetilde{\mfrak{g}}_\kappa,\mfrak{s}_\alpha^+)$ and $\mcal{I}(\widetilde{\mfrak{g}}_\kappa,\mfrak{s}_\alpha^-)$ of $\mcal{E}(\widetilde{\mfrak{g}}_\kappa)$ assigned to the Lie subalgebras $\mfrak{s}_\alpha^+$ and $\mfrak{s}_\alpha^-$ of $\widetilde{\mfrak{g}}_\kappa$.
\medskip

The following statement is obvious, see Proposition 1.3 in \cite{Futorny-Krizka2019} for details.
\medskip

\proposition{The categories $\mcal{I}(\widetilde{\mfrak{g}}_\kappa,\mfrak{s}_\alpha^+)$ and $\mcal{I}(\widetilde{\mfrak{g}}_\kappa,\mfrak{s}_\alpha^-)$ are full subcategories of $\mcal{H}(\widetilde{\mfrak{g}}_\kappa, \Gamma_\alpha)$ for $\alpha \in \smash{\widehat{\Delta}}^{\rm re}_+$.}

\vspace{-2mm}

%%%%%%%%%%%%%%%%%%%%%%%%%%%%%%%%%%%%%%%%%%%%%%%%%%%%%%%%%%%%%%%%%%%%%%%%%%%%%%%%%%%%%%%%%%

\subsection{Twisting functors for affine Kac--Moody algebras}

For a real positive root $\alpha \in \smash{\widehat{\Delta}}^{\rm re}_+$ the multiplicative sets $\{f_\alpha^n;\, n \in \N_0\}$ and $\{e_\alpha^n;\, n \in \N_0\}$ in the universal enveloping algebra $U(\widehat{\mfrak{g}}_\kappa)$ of $\widehat{\mfrak{g}}_\kappa$ are left (right) denominator sets, since $f_\alpha$ and $e_\alpha$ are locally $\ad$-nilpotent regular elements. Therefore, based on the general construction we define the \emph{twisting functor}
\begin{align*}
  T_\alpha=T_{f_\alpha} \colon \mcal{M}(\widehat{\mfrak{g}}_\kappa) \rarr \mcal{M}(\widehat{\mfrak{g}}_\kappa) \qquad \text{and} \qquad T_{-\alpha} =T_{e_\alpha} \colon \mcal{M}(\widehat{\mfrak{g}}_\kappa) \rarr \mcal{M}(\widehat{\mfrak{g}}_\kappa)
\end{align*}
by
\begin{align*}
  T_\alpha(M) = U(\widehat{\mfrak{g}}_\kappa)_{(f_\alpha)}/U(\widehat{\mfrak{g}}_\kappa) \otimes_{U(\widehat{\mfrak{g}}_\kappa)}\! M, \qquad
  T_{-\alpha}(M) = U(\widehat{\mfrak{g}}_\kappa)_{(e_\alpha)}/U(\widehat{\mfrak{g}}_\kappa) \otimes_{U(\widehat{\mfrak{g}}_\kappa)}\! M
\end{align*}
for $\alpha \in \smash{\widehat{\Delta}}^{\rm re}_+$ and $M \in \mcal{M}(\widehat{\mfrak{g}}_\kappa)$. In addition, it is easy to see that if $M$ is a $\widetilde{\mfrak{g}}_\kappa$-module, then $T_\alpha(M)$ has also a natural structure of a $\widetilde{\mfrak{g}}_\kappa$-module for $\alpha \in \smash{\widehat{\Delta}}^{\rm re}$.
\medskip

In the next, we prove some basic characteristics of $T_\alpha$ for $\alpha \in \smash{\widehat{\Delta}}^{\rm re}$. Let us note that the twisting functor $T_\alpha$ for $\alpha \in \Delta \subset \smash{\widehat{\Delta}}^{\rm re}$ is of a special importance and will be discussed later in details.
\medskip

Let us consider the subset
\begin{align*}
  \Phi_\alpha = \{\gamma \in \widehat{\Delta} \setminus \{\pm\alpha\};\, \widehat{\mfrak{g}}_{\kappa,\gamma,\alpha} \subset \widehat{\mfrak{n}}_{\rm st}\}
\end{align*}
of $\smash{\widehat{\Delta}}$ for $\alpha \in \smash{\widehat{\Delta}}^{\rm re}_+$, where $\widehat{\mfrak{g}}_{\kappa,\gamma,\alpha}$ for $\gamma \in \smash{\widehat{\Delta}}$ is the finite-dimensional $\mfrak{s}_\alpha$-module given by
\begin{align*}
  \widehat{\mfrak{g}}_{\kappa,\gamma,\alpha} = \bigoplus_{j \in \Z} \widehat{\mfrak{g}}_{\kappa,\gamma+j\alpha}.
\end{align*}
It easily follows that $\Phi_\alpha \subset \smash{\widehat{\Delta}}_+$ and that $\Phi_\alpha$ is closed, i.e.\ if $\gamma_1, \gamma_2 \in \Phi_\alpha$ and $\gamma_1+\gamma_2 \in \smash{\widehat{\Delta}}$, then we have $\gamma_1+\gamma_2 \in \Phi_\alpha$. Therefore, the subset $\Phi_\alpha$ of $\smash{\widehat{\Delta}}_+$ gives rise to the Lie subalgebras $\mfrak{t}_\alpha^+$, $\mfrak{t}_\alpha^-$ and $\mfrak{t}_\alpha$ of $\widehat{\mfrak{g}}_\kappa$ defined through
\begin{align*}
  \mfrak{t}_\alpha^+ = \widehat{\mfrak{g}}_{\kappa,\alpha} \oplus \bigoplus_{\gamma \in \Phi_\alpha} \widehat{\mfrak{g}}_{\kappa,\gamma}, \qquad \mfrak{t}_\alpha = \bigoplus_{\gamma \in \Phi_\alpha} \widehat{\mfrak{g}}_{\kappa,\gamma}, \qquad  \mfrak{t}_\alpha^- = \widehat{\mfrak{g}}_{\kappa,-\alpha} \oplus \bigoplus_{\gamma \in \Phi_\alpha} \widehat{\mfrak{g}}_{\kappa,\gamma},
\end{align*}
where the direct sum over the set $\Phi_\alpha$ is the topological direct sum. Moreover, since we have the inclusions $\mfrak{s}_\alpha^+ \subset \mfrak{t}_\alpha^+ \subset \widehat{\mfrak{n}}_{\rm st}$ of Lie algebras, we get immediately the embeddings of categories
\begin{align*}
  \mcal{O}(\widetilde{\mfrak{g}}_\kappa) \subset \mcal{I}_f(\widetilde{\mfrak{g}}_\kappa, \mfrak{t}_\alpha^+) \subset \mcal{I}_f(\widetilde{\mfrak{g}}_\kappa,\mfrak{s}_\alpha^+)
\end{align*}
for $\alpha \in \smash{\widehat{\Delta}}^{\rm re}_+$. There is also an automorphism $\Ad(\dot{r}_\alpha)$ of $\widetilde{\mfrak{g}}_\kappa$ given by
\begin{align*}
  \Ad(\dot{r}_\alpha) = \exp(\ad f_\alpha) \exp(- \ad e_\alpha) \exp(\ad f_\alpha),
\end{align*}
which satisfies
\begin{align*}
\Ad(\dot{r}_\alpha)(\mfrak{s}_\alpha^\pm) = \mfrak{s}_\alpha^\mp \qquad \text{and} \qquad \Ad(\dot{r}_\alpha)(\mfrak{t}_\alpha^\pm) = \mfrak{t}_\alpha^\mp
\end{align*}
for $\alpha \in \smash{\widehat{\Delta}}^{\rm re}_+$.
\medskip

\theorem{\label{thm:twisting functor cat E}
For $\alpha \in \smash{\widehat{\Delta}}^{\rm re}$, the twisting functor $T_\alpha$ preserves the category $\mcal{E}(\widehat{\mfrak{g}}_\kappa)$. Moreover, if $\alpha \in \Delta \subset \smash{\widehat{\Delta}}^{\rm re}$, then the twisting functor $T_\alpha$ preserves also the category $\mcal{E}_+\!(\widehat{\mfrak{g}}_\kappa)$.}

\proof{Let $M$ be a smooth $\widehat{\mfrak{g}}_\kappa$-module on which the central element $c$ acts as the identity. Then for each $v \in M$ there exists $N_v \in \N$ such that $(\mfrak{g} \otimes_\C t^{N_v}\C[[t]])v=0$. Further, for $\alpha \in \smash{\widehat{\Delta}}^{\rm re}_+$ there exists $k_\alpha \in \N$ satisfying $\ad(f_\alpha)^{k_\alpha}(a)=0$ for all $a \in \widehat{\mfrak{g}}_\kappa$. Hence, we may write
\begin{align*}
  af_\alpha^{-n}v = \sum_{k=0}^{k_\alpha-1} \binom{n+k-1}{k}f_\alpha^{-n-k} \ad(f_\alpha)^k(a) v
\end{align*}
for $a \in \widehat{\mfrak{g}}_\kappa$. Moreover, since $\ad(f_\alpha)^k(a_m) \in \mfrak{g} \otimes_\C t^{m-k\alpha(d)}\C[[t]]$ for $a \in \mfrak{g}$ and $m\in\Z$, we get that $(\mfrak{g} \otimes_\C t^{N_{n,v}}\C[[t]])f_\alpha^{-n}v=0$ for $N_{n,v}=\max\{N_v+k\alpha(d);\, k\in\{0,1,\dots,k_\alpha-1\}\}$. Hence, we have that $T_\alpha(M)$ for $\alpha \in \smash{\widehat{\Delta}}^{\rm re}_+$ is a smooth $\widehat{\mfrak{g}}_\kappa$-module on which the central element $c$ acts as the identity. For $\alpha \in \smash{\widehat{\Delta}}^{\rm re}_-$ the proof goes along the same lines.

Further, let $M$ belong to the category $\mcal{E}_+\!(\widehat{\mfrak{g}}_\kappa)$. Then from the previous part we immediately get that $T_\alpha(M)$ for $\alpha \in \Delta$ is from the category $\mcal{E}(\widehat{\mfrak{g}}_\kappa)$. As $M$ is a positive energy $\widehat{\mfrak{g}}_\kappa$-module we have $M = \bigoplus_{n = 0}^\infty M_{\lambda+n}$ with $M_\lambda \ne 0$ for some $\lambda \in \C$. By setting $\deg f_\alpha^{-n}v = \deg v$ for $n \in \N_0$ and a homogeneous vector $v \in M$, it follows that $T_\alpha(M)$ is a positive energy $\widehat{\mfrak{g}}_\kappa$-module.}

Let us note that for $\alpha \in \smash{\widehat{\Delta}}^{\rm re}_- \setminus \Delta_-$ the twisting functor $T_\alpha$ also preserves the category $\mcal{E}_+\!(\widehat{\mfrak{g}}_\kappa)$. However, we have $T_\alpha(M) = 0$ for $\alpha \in \smash{\widehat{\Delta}}^{\rm re}_- \setminus \Delta_-$ and $M \in \mcal{E}_+\!(\widehat{\mfrak{g}}_\kappa)$ which follows from Proposition \ref{prop:locally nilpotent} and the fact that $M$ is a locally $\mfrak{s}_\alpha^+$-finite $\widehat{\mfrak{g}}_\kappa$-module.
\medskip

\proposition{\label{prop:locally nilpotent}
Let $M$ be a $\widehat{\mfrak{g}}_\kappa$-module. Then $T_\alpha(M)$ is a locally $\mfrak{s}_\alpha^-$-finite $\widehat{\mfrak{g}}_\kappa$-module for any $\alpha \in \smash{\widehat{\Delta}}^{\rm re}_+$. Moreover, if $M$ is a locally $\mfrak{s}_\alpha^-$-finite $\widehat{\mfrak{g}}_\kappa$-module, then we obtain $T_\alpha(M)=0$.}

\proof{Let $M$ be a $\widehat{\mfrak{g}}_\kappa$-module. Then by definition we have
\begin{align*}
  T_\alpha(M) = (U(\widehat{\mfrak{g}}_\kappa)_{(f_\alpha)}/U(\widehat{\mfrak{g}}_\kappa)) \otimes_{U(\widehat{\mfrak{g}}_\kappa)}\! M \simeq M_{(f_\alpha)}/M
\end{align*}
for $\alpha \in \smash{\widehat{\Delta}}^{\rm re}_+$. Since every element of $M_{(f_\alpha)}$ can be written in the form $f_\alpha^{-n} v$ for $n \in \N_0$ and $v \in M$, we obtain immediately that $T_\alpha(M)$ is locally $\mfrak{s}_\alpha^-$-finite. Further, let us assume that $M$ is a locally $\mfrak{s}_\alpha^-$-finite $\widehat{\mfrak{g}}_\kappa$-module. Since for each $v \in M$ there exists $n_v \in \N_0$ such that $f_\alpha^{n_v}v=0$, we may write $f_\alpha^{-n}v = f_\alpha^{-n-n_v}f_\alpha^{n_v}v = 0$ for $n\in \N_0$. This implies the required statement.}

\lemma{\label{lem:sl2 action}
Let $\alpha \in \smash{\widehat{\Delta}}^{\rm re}_+$. Then we have
\begin{align*}
\begin{aligned}
  e_\alpha f_\alpha^{-n} &= f_\alpha^{-n} e_\alpha -nf_\alpha^{-n-1} h_\alpha - n(n+1)f_\alpha^{-n-1}, \\
  h f_\alpha^{-n} &= f_\alpha^{-n} h + n\alpha(h) f_\alpha^{-n}, \\
  f_\alpha f_\alpha^{-n} &= f_\alpha^{-n} f_\alpha
\end{aligned}
\end{align*}
in $U(\widehat{\mfrak{g}}_\kappa)_{(f_\alpha)}$ for $n \in \Z$ and $h \in \smash{\widehat{\mfrak{h}}}$.}

\proof{It follows immediately from the formula
\begin{align*}
  a f_\alpha^{-n} = f_\alpha^{-n} \sum_{k=0}^\infty \binom{n+k-1}{k} f_\alpha^{-k} \ad(f_\alpha)^k(a)
\end{align*}
in $U(\widehat{\mfrak{g}}_\kappa)_{(f_\alpha)}$ for all $a \in U(\widehat{\mfrak{g}}_\kappa)$ and $n \in \N_0$.}

\theorem{\label{thm:finitely generated}
We have
\begin{enumerate}[topsep=3pt,itemsep=0pt]
  \item[i)] if $M \in \mcal{I}_f(\widetilde{\mfrak{g}}_\kappa,\mfrak{s}_\alpha^+)$, then $T_\alpha(M) \in \mcal{I}_f(\widetilde{\mfrak{g}}_\kappa,\mfrak{s}_\alpha^-)$;
  \item[ii)] if $M \in \mcal{I}_f(\widetilde{\mfrak{g}}_\kappa,\mfrak{t}_\alpha^+)$, then $T_\alpha(M) \in \mcal{I}_f(\widetilde{\mfrak{g}}_\kappa,\mfrak{t}_\alpha^-)$
\end{enumerate}
for $\alpha \in \smash{\widehat{\Delta}}^{\rm re}_+$. Therefore, we have the restricted functors
\begin{align*}
  T_\alpha \colon \mcal{I}_f(\widetilde{\mfrak{g}}_\kappa,\mfrak{s}_\alpha^+) \rarr \mcal{I}_f(\widetilde{\mfrak{g}}_\kappa,\mfrak{s}_\alpha^-) \qquad \text{and} \qquad T_\alpha \colon \mcal{I}_f(\widetilde{\mfrak{g}}_\kappa,\mfrak{t}_\alpha^+) \rarr \mcal{I}_f(\widetilde{\mfrak{g}}_\kappa,\mfrak{t}_\alpha^-)
\end{align*}
for $\alpha \in \smash{\widehat{\Delta}}^{\rm re}_+$.}

\proof{i) If $M$ is a weight $\widetilde{\mfrak{g}}_\kappa$-module, then it easily follows from definition that $T_\alpha(M)$ is also a weight $\widetilde{\mfrak{g}}_\kappa$-module. Moreover, by Theorem \ref{thm:twisting functor cat E} and Proposition \ref{prop:locally nilpotent} we have that $T_\alpha(M)$ belongs to the category $\mcal{I}(\widetilde{\mfrak{g}}_\kappa, \mfrak{s}_\alpha^-)$. Hence, the rest of the proof is to show that $T_\alpha(M)$ is finitely generated provided $M$ is finitely generated and locally $\mfrak{s}_\alpha^+$-finite.

Let $R \subset M$ be a finite set of generators of $M$. Then the vector subspace $V = U(\mfrak{s}_\alpha^+)\langle R \rangle$ of $M$ is finite dimensional. Further, let us introduce a filtration $\{F_k V\}_{k \in \N_0}$ on $V$ by
\begin{align*}
  F_k V = \{v \in V;\, e_\alpha^kv=0\}
\end{align*}
for $k \in \N_0$ and let $n_0 \in \N$ be the smallest positive integer satisfying
\begin{align*}
  n_0 \geq \max\{-\mu(h_\alpha);\, \mu \in \smash{\widetilde{\mfrak{h}}}^*,\, V_\mu \neq \{0\},\, \mu(h_\alpha) \in \R\}.
\end{align*}
Let us consider a vector $v \in F_kV \cap V_\mu$ for $\mu \in \smash{\widetilde{\mfrak{h}}}^*$ and $k \in \N_0$. Then by Lemma \ref{lem:sl2 action} we get
\begin{align*}
   e_\alpha f_\alpha^{-(n_0+n)}v = f_\alpha^{-(n_0+n)}e_\alpha v - (n_0+n)(\mu(h_\alpha)+n_0+n+1) f_\alpha^{-(n_0+n+1)}v
\end{align*}
for $n \in \N_0$, which together with $f_\alpha f_\alpha^{-n-1}v = f_\alpha^{-n}v$ for $n\in \N_0$ gives us
\begin{align*}
  U(\mfrak{s}_\alpha)f_\alpha^{-n_0}F_kV/\C[f_\alpha^{-1}]F_{k-1}V = \C[f_\alpha^{-1}]F_kV/\C[f_\alpha^{-1}]F_{k-1}V
\end{align*}
for all $k \in \N$. As $F_0V = \{0\}$ and $V$ is a finite-dimensional vector space, we immediately obtain $U(\mfrak{s}_\alpha)f_\alpha^{-n_0}V = \C[f_\alpha^{-1}]V$.

Further, since for any $a \in U(\widehat{\mfrak{g}}_\kappa)$ and $n \in \N_0$ there exist $b \in U(\widehat{\mfrak{g}}_\kappa)$ and $m \in \N_0$ satisfying $f_\alpha^{-n}a = bf_\alpha^{-m}$, we obtain $f_\alpha^{-n}U(\widehat{\mfrak{g}}_\kappa) \subset U(\widehat{\mfrak{g}}_\kappa)\C[f_\alpha^{-1}] \subset U(\widehat{\mfrak{g}}_\kappa)_{(f_\alpha)}$ for $n \in \N_0$. Hence, we may write
\begin{align*}
  f_\alpha^{-n}M = f_\alpha^{-n}U(\widehat{\mfrak{g}}_\kappa)V \subset U(\widehat{\mfrak{g}}_\kappa)\C[f_\alpha^{-1}]V,
\end{align*}
which implies $M_{(f_\alpha)} =U(\widehat{\mfrak{g}}_\kappa)\C[f_\alpha^{-1}]V = U(\widehat{\mfrak{g}}_\kappa)f_\alpha^{-n_0}V$. In other words, this means that $T_\alpha(M) \simeq M_{(f_\alpha)}/M$ is finitely generated and the number of generators is bounded by $\dim V$.
\smallskip

ii) Since $\mcal{I}_f(\widetilde{\mfrak{g}}_\kappa,\mfrak{t}_\alpha^+)$ is a full subcategory of $\mcal{I}_f(\widetilde{\mfrak{g}}_\kappa,\mfrak{s}_\alpha^+)$, by item (i) we only need to show that $T_\alpha(M)$ is locally $\mfrak{t}_\alpha^-$-finite if $M$ is locally $\mfrak{t}_\alpha^+$-finite. As we have $U(\mfrak{t}_\alpha^-) = U(\mfrak{s}_\alpha^-) U(\mfrak{t}_\alpha)$ and $T_\alpha(M)$ is locally $\mfrak{s}_\alpha^-$-finite, the condition that $T_\alpha(M)$ is locally $\mfrak{t}_\alpha^-$-finite is equivalent to saying that $T_\alpha(M)$ is locally $\mfrak{t}_\alpha$-finite. Further, since $T_\alpha(M)$ is a smooth $\widehat{\mfrak{g}}_\kappa$-module and the quotient $\mfrak{t}_\alpha/ (\mfrak{t}_\alpha \cap (\mfrak{g} \otimes_\C t^N\C[[t]]))$ for $N \in \N$ is a finite-dimensional nilpotent Lie algebra, we get that $T_\alpha(M)$ is locally $\mfrak{t}_\alpha$-finite if the commutative Lie algebra $\widehat{\mfrak{g}}_{\kappa,\gamma}$ acts locally nilpotently on $T_\alpha(M)$ for all $\gamma \in \Phi_\alpha$.

Let $\htt \colon \Z\smash{\widehat{\Delta}} \rarr \Z$ be the $\Z$-linear height function with $\htt(\alpha)=1$ if $\alpha$ is a simple root. Then we get an $\N_0$-grading $U(\mfrak{t}_\alpha) = \bigoplus_{n \in \N_0} U(\mfrak{t}_\alpha)_n$, where
\begin{align*}
  U(\mfrak{t}_\alpha)_n = \bigoplus_{\mu \in \smash{\widetilde{\mfrak{h}}}^*\!,\, \htt(\mu)=n} U(\mfrak{t}_\alpha)_\mu.
\end{align*}
Let $\gamma \in \Phi_\alpha$ and let $r \in \N_0$ be the smallest nonnegative integer such that $\gamma-(r+1)\alpha \notin \Phi_\alpha$. Let us recall that if $\gamma - k\alpha \in \smash{\widehat{\Delta}}$ for some $k \in \Z$ then $\gamma - k\alpha \in \Phi_\alpha$. Let us consider a vector $v \in M$. Then from the formula
\begin{align*}
  a f_\alpha^{-n} = f_\alpha^{-n} \sum_{k=0}^\infty \binom{n+k-1}{k} f_\alpha^{-k} \ad(f_\alpha)^k(a)
\end{align*}
in $U(\widehat{\mfrak{g}}_\kappa)_{(f_\alpha)}$ for $a \in U(\widehat{\mfrak{g}}_\kappa)$ and $n \in \N_0$, we obtain
\begin{align*}
  e_\gamma^t f_\alpha^{-n}v = f_\alpha^{-n} \sum_{k=0}^{tr} \binom{n+k-1}{k} f_\alpha^{-k} \ad(f_\alpha)^k(e_\gamma^t)v
\end{align*}
for $e_\gamma \in \widehat{\mfrak{g}}_{\kappa,\gamma}$, $t \in \N$ and $n \in \N_0$, where we used the fact that $\ad(f_\alpha)^k(e_\gamma) \neq 0$ only for $k=0,1,\dots,r$. Moreover, we have $\ad(f_\alpha)^k(e_\gamma^t) \in U(\mfrak{t}_\alpha)_{t\gamma-k\alpha} \subset U(\mfrak{t}_\alpha)_{\htt(t\gamma-k\alpha)}$ for $k=0,1,\dots,tr$. Since $M$ is locally $\mfrak{t}_\alpha^+$-finite, there exists an integer $n_v \in \N_0$ such that $U(\mfrak{t}_\alpha)_nv=\{0\}$ for $n > n_v$. Therefore, it is enough to show that $\htt(t\gamma-k\alpha) > n_v$ for $k=0,1,\dots,tr$.

As we may write
\begin{align*}
  \htt(t\gamma-k\alpha) \geq \htt(t\gamma-tr\alpha) = t\htt(\gamma-r\alpha) \geq t
\end{align*}
for $k = 0,1,\dots,tr$, since $\gamma - r \alpha \in \Phi_\alpha$ and hence $\htt(\gamma -r\alpha) \geq 1$, we obtain that $e_\gamma^t f_\alpha^{-n}v=0$ for $n \in \N_0$ provided $t > n_v$. Hence, the element $e_\gamma$ acts locally nilpotently on $T_\alpha(M)$ for $\gamma \in \Phi_\alpha$.}

The following shows how the twisting functor $T_\alpha$ can be used to construct $\alpha$-Gelfand--Tsetlin modules with finite $\Gamma_\alpha$-multiplicities for $\alpha \in \smash{\widehat{\Delta}}^{\rm re}_+$.
\medskip

\theorem{\label{thm:highest weight modules condition}
Let $\alpha \in \smash{\widehat{\Delta}}^{\rm re}_+$ and let $M$ be a smooth weight $\widetilde{\mfrak{g}}_\kappa$-module on which the central element $c$ acts as the identity.
\begin{itemize}[topsep=3pt,itemsep=0pt]
\item[i)] The $\widetilde{\mfrak{g}}_\kappa$-module $T_\alpha(M)$ is an $\alpha$-Gelfand--Tsetlin module with finite $\Gamma_\alpha$-multiplicities if and only if the first cohomology group $H^1(\mfrak{s}_\alpha^-;M)$ is a weight $\smash{\widetilde{\mfrak{h}}}$-module with finite-dimensional weight spaces.
\item[ii)] If $M$ is a highest weight $\widetilde{\mfrak{g}}_\kappa$-module, then $T_\alpha(M)$ is a locally $\mfrak{t}_\alpha^-$-finite smooth cyclic weight $\widetilde{\mfrak{g}}_\kappa$-module with finite $\Gamma_\alpha$-multiplicities.
\end{itemize}
}

\proof{i) Let $M$ be a smooth weight $\widetilde{\mfrak{g}}_\kappa$-module. Then by Proposition \ref{prop:locally nilpotent} and Theorem \ref{thm:twisting functor cat E} we obtain that $T_\alpha(M)$ is a locally $\mfrak{s}_\alpha^-$-finite smooth weight $\widetilde{\mfrak{g}}_\kappa$-module for $\alpha \in \smash{\widehat{\Delta}}^{\rm re}_+$. Hence, we may apply Corollary \ref{cor:highest weight GT modules} on $T_\alpha(M)$ and we get that $T_\alpha(M)$ is an $\alpha$-Gelfand--Tsetlin module with finite $\Gamma_\alpha$-multiplicities if and only if $H^0(s_\alpha^-;T_\alpha(M))$ is a weight $\smash{\widetilde{\mfrak{h}}}$-module with finite-dimensional weight spaces. Further, the linear mapping
$\varphi_\alpha \colon M \rarr T_\alpha(M)$ defined by
\begin{align*}
  \varphi_\alpha(v) = f_\alpha^{-1}v
\end{align*}
for $v \in M$ gives rise to the linear mapping
\begin{align*}
 \widetilde{\varphi}_\alpha \colon H^1(\mfrak{s}_\alpha^-;M) \rarr H^0(\mfrak{s}_\alpha^-;T_\alpha(M))
\end{align*}
for $\alpha \in \smash{\widehat{\Delta}}^{\rm re}_+$, which is in fact an isomorphism. Therefore, we obtain an isomorphism $H^1(\mfrak{s}_\alpha^-;M) \simeq H^0(\mfrak{s}_\alpha^-;T_\alpha(M)) \otimes_\C \C_{-\alpha}$ of $\smash{\widetilde{\mfrak{h}}}$-modules, where $\C_{-\alpha}$ is the $1$-dimensional $\smash{\widetilde{\mfrak{h}}}$-module determined by the character $-\alpha$ of $\smash{\widetilde{\mfrak{h}}}$, which implies the first statement.
\smallskip

ii) If $M$ is a highest weight $\widetilde{\mfrak{g}}_\kappa$-module, then it belongs to the category $\mcal{O}(\widetilde{\mfrak{g}}_\kappa)$. Hence, using Theorem \ref{thm:finitely generated} we obtain that $T_\alpha(M)$ is a locally $\mfrak{t}_\alpha^-$-finite smooth weight $\widetilde{\mfrak{g}}_\kappa$-module. In fact, from the proof of Theorem \ref{thm:finitely generated} (i) it follows that $T_\alpha(M)$ is not only finitely generated by also cyclic. To finish the proof, we need to show by Theorem \ref{thm:highest weight modules condition} that $H^1(\mfrak{s}_\alpha^-;M)$ is a weight $\smash{\widetilde{\mfrak{h}}}$-module with finite-dimensional weight spaces. Since $H^1(\mfrak{s}_\alpha^-;M) \simeq M/\mfrak{s}_\alpha^-M$ and $M$ is a weight $\smash{\widetilde{\mfrak{h}}}$-module with finite-dimensional weight spaces, we immediately obtain that also $M/\mfrak{s}_\alpha^-M$ is a weight $\smash{\widetilde{\mfrak{h}}}$-module with finite-dimensional weight spaces. This gives us the required statement.}

If we denote by $\Theta_\alpha \colon \mcal{M}(\widetilde{\mfrak{g}}_\kappa) \rarr \mcal{M}(\widetilde{\mfrak{g}}_\kappa)$ the functor sending a $\widetilde{\mfrak{g}}_\kappa$-module to the same $\widetilde{\mfrak{g}}_\kappa$-module with the action twisted by the automorphism $\Ad(\dot{r}_\alpha) \colon \widetilde{\mfrak{g}}_\kappa \rarr \widetilde{\mfrak{g}}_\kappa$, then we obtain the endofunctor
\begin{align*}
 \Theta_\alpha \circ T_\alpha \colon \mcal{I}_f(\widetilde{\mfrak{g}}_\kappa,\mfrak{t}_\alpha^+) \rarr \mcal{I}_f(\widetilde{\mfrak{g}}_\kappa,\mfrak{t}_\alpha^+).
\end{align*}
Besides, for $\alpha \in \smash{\widehat{\Pi}}$ we have $\mfrak{t}_\alpha^+ = \widehat{\mfrak{n}}_{\rm st}$ which implies $\mcal{I}_f(\widetilde{\mfrak{g}}_\kappa,\mfrak{t}_\alpha^+) = \mcal{O}(\widetilde{\mfrak{g}}_\kappa)$ and in this case the functor coincides with the Arkhipov's twisting functor, see \cite{Arkhipov1997}, \cite{Arkhipov2004}.

%%%%%%%%%%%%%%%%%%%%%%%%%%%%%%%%%%%%%%%%%%%%%%%%%%%%%%%%%%%%%%%%%%%%%%%%%%%%%%%%%%%%%%%%%%

\subsection{Tensoring with Weyl modules}

In this subsection we show that the twisting functors behave well with respect to tensoring with certain $\widehat{\mfrak{g}}_\kappa$-modules. In the finite-dimensional setting this was considered in \cite{Futorny-Krizka2019b} and \cite{Andersen-Stroppel2003}.
\medskip

Let us recall that the universal enveloping algebra $U(\widehat{\mfrak{g}}_\kappa)$ is a Hopf algebra with the comultiplication $\Delta \colon U(\widehat{\mfrak{g}}_\kappa) \rarr U(\widehat{\mfrak{g}}_\kappa) \otimes_\C U(\widehat{\mfrak{g}}_\kappa)$, the counit $\veps \colon U(\widehat{\mfrak{g}}_\kappa) \rarr \C$ and the antipode $S \colon U(\widehat{\mfrak{g}}_\kappa) \rarr U(\widehat{\mfrak{g}}_\kappa)$ given by
\begin{align*}
  \Delta(a) = a \otimes 1 + 1 \otimes a, \qquad \veps(a) = 0, \qquad S(a) = -a
\end{align*}
for $a \in \widehat{\mfrak{g}}_\kappa$. For $\alpha \in \smash{\widehat{\Delta}}^{\rm re}_+$, the localization $U(\widehat{\mfrak{g}}_\kappa)_{(f_\alpha)}$ has the structure of a left $\C[f_\alpha^{-1}]$-module, hence also $U(\widehat{\mfrak{g}}_{\kappa})_{(f_\alpha)} \otimes_\C U(\widehat{\mfrak{g}}_\kappa)_{(f_\alpha)}$ is a left $\C[f_\alpha^{-1}]$-module and we denote by $U(\widehat{\mfrak{g}}_\kappa)_{(f_\alpha)} \,\smash{\widehat{\otimes}_\C}\, U(\widehat{\mfrak{g}}_\kappa)_{(f_\alpha)}$ its extension to a left $\C[[f_\alpha^{-1}]]$-module, i.e.\ we set
\begin{align*}
  U(\widehat{\mfrak{g}}_\kappa)_{(f_\alpha)} \,\widehat{\otimes}_\C\, U(\widehat{\mfrak{g}}_\kappa)_{(f_\alpha)} = \C[[f_\alpha^{-1}]] \otimes_{\C[f_\alpha^{-1}]}  U(\widehat{\mfrak{g}}_\kappa)_{(f_\alpha)} \otimes_\C U(\widehat{\mfrak{g}}_\kappa)_{(f_\alpha)}.
\end{align*}
There is an obvious extension of the algebra structure on $U(\widehat{\mfrak{g}}_\kappa)_{(f_\alpha)} \otimes_\C U(\widehat{\mfrak{g}}_\kappa)_{(f_\alpha)}$ to the completion $U(\widehat{\mfrak{g}}_\kappa)_{(f_\alpha)} \,\smash{\widehat{\otimes}_\C}\, U(\widehat{\mfrak{g}}_\kappa)_{(f_\alpha)}$. Then the linear mapping
\begin{align*}
  \widetilde{\Delta} \colon U(\widehat{\mfrak{g}}_\kappa)_{(f_\alpha)} &\rarr U(\widehat{\mfrak{g}}_\kappa)_{(f_\alpha)}\, \widehat{\otimes}_\C \, U(\widehat{\mfrak{g}}_\kappa)_{(f_\alpha)}
\end{align*}
given through
\begin{align*}
  \widetilde{\Delta}(f_\alpha^{-n}u) = \bigg(\sum_{k=0}^\infty (-1)^k \binom{n+k-1}{k} f_\alpha^{-n-k} \otimes f_\alpha^k \bigg) \Delta(u)
\end{align*}
for $n \in \N_0$ and $u \in U(\widehat{\mfrak{g}}_\kappa)$ defines an algebra homomorphism, see \cite{Andersen-Stroppel2003}. The following theorem is analogous to \cite[Theorem 3.2]{Andersen-Stroppel2003}.
\medskip

\theorem{\label{thm-And}
Let $\alpha \in \smash{\widehat{\Delta}}^{\rm re}_\pm$. Then there exists a family $\{\eta_E\}_{E \in \mcal{M}(\widehat{\mfrak{g}}_\kappa, \mfrak{s}_\alpha^\mp)}$ of natural isomorphisms
\begin{align*}
  \eta_E \colon T_\alpha \circ (\,\bullet \otimes_\C E) \rarr (\,\bullet \otimes_\C E) \circ T_\alpha
\end{align*}
of functors, where $\mcal{M}(\widehat{\mfrak{g}}_\kappa, \mfrak{s}_\alpha^\pm)$ is the category of locally $\mfrak{s}_\alpha^\pm$-finite $\widehat{\mfrak{g}}_\kappa$-modules.}

\remark{Theorem \ref{thm-And} implies that for a root $\alpha \in \Delta \subset \smash{\widehat{\Delta}}^{\rm re}$ the twisting functor $T_\alpha$ commutes with tensoring by Weyl modules.}

%%%%%%%%%%%%%%%%%%%%%%%%%%%%%%%%%%%%%%%%%%%%%%%%%%%%%%%%%%%%%%%%%%%%%%%%%%%%%%%%%%%%%%%%%%

\subsection{Relaxed Verma modules}

Let us recall that we have the induction functor
\begin{align}
  \mathbb{M}_{\kappa,\mfrak{g}} \colon \mcal{M}(\mfrak{g}) \rarr \mcal{E}_+\!(\widehat{\mfrak{g}}_\kappa) \label{eq:induction functor}
\end{align}
introduced in Section \ref{subsec:induced} by $\mathbb{M}_{\kappa,\mfrak{g}}(E) = U(\widehat{\mfrak{g}}_\kappa) \otimes_{U(\widehat{\mfrak{g}}_{\rm st})}\! E$ for a $\mfrak{g}$-module $E$, where $E$ is considered as $\widehat{\mfrak{g}}_{\rm st}$-module on which $\mfrak{g} \otimes_\C t\C[[t]]$ acts trivially and $c$ acts as the identity.
\medskip

The link between the twisting functor $T_\alpha$ for $\alpha \in \Delta \subset \smash{\widehat{\Delta}}^{\rm re}$ and the induction functor $\mathbb{M}_{\kappa,\mfrak{g}}$ is given in the following theorem.
\medskip

\theorem{\label{thm:twisting functor intertwining}
Let $\alpha \in \Delta \subset \smash{\widehat{\Delta}}^{\rm re}$. Then there exists a natural isomorphism
\begin{align*}
  \eta_\alpha \colon T_\alpha \circ \mathbb{M}_{\kappa,\mfrak{g}} \rarr \mathbb{M}_{\kappa,\mfrak{g}} \circ\, T_\alpha^\mfrak{g}
\end{align*}
of functors, where $T_\alpha^\mfrak{g} \colon \mcal{M}(\mfrak{g}) \rarr \mcal{M}(\mfrak{g})$ is the twisting functor for $\mfrak{g}$ assigned to $\alpha$. In particular, we have
\begin{align*}
  T_\alpha(\mathbb{M}_{\kappa,\mfrak{g}}(M^\mfrak{g}_\mfrak{p}(\lambda))) \simeq \mathbb{M}_{\kappa,\mfrak{g}}(W^\mfrak{g}_\mfrak{p}(\lambda,\alpha))
\end{align*}
for $\lambda \in \Lambda^+(\mfrak{p})$ and $\alpha \in \Delta_+^\mfrak{u}$.}

\proof{We prove the statement only for $\alpha \in \Delta_+$, since for $\alpha \in \Delta_-$ the proof goes along the same lines. By using the triangular decomposition $\widehat{\mfrak{g}}_\kappa = \widehat{\mfrak{g}}_{\kappa,-} \oplus \widehat{\mfrak{g}}_{\kappa,0} \oplus \widehat{\mfrak{g}}_{\kappa,+}$, where
\begin{align*}
  \widehat{\mfrak{g}}_{\kappa,-} = \mfrak{g} \otimes_\C t^{-1}\C[t^{-1}], \qquad \widehat{\mfrak{g}}_{\kappa,0} = \mfrak{g} \otimes_\C \C 1 \oplus \C c, \qquad \widehat{\mfrak{g}}_{\kappa,+} = \mfrak{g} \otimes_\C t\C[[t]],
\end{align*}
and the Poincar\'e--Birkhoff--Witt theorem we get an isomorphism
\begin{align*}
  U(\widehat{\mfrak{g}}_\kappa)_{(f_\alpha)} \simeq U(\widehat{\mfrak{g}}_{\kappa,-})_{(f_\alpha)} \otimes_\C U(\widehat{\mfrak{g}}_{\rm st})
\end{align*}
of $U(\widehat{\mfrak{g}}_{\kappa,-})$-modules for $\alpha \in \Delta_+$. Hence, we may write
\begin{align*}
  U(\widehat{\mfrak{g}}_\kappa)_{(f_\alpha)} \otimes_{U(\widehat{\mfrak{g}}_\kappa)} \mathbb{M}_{\kappa,\mfrak{g}}(E) \simeq  U(\widehat{\mfrak{g}}_\kappa)_{(f_\alpha)} \otimes_{U(\widehat{\mfrak{g}}_\kappa)} U(\widehat{\mfrak{g}}_\kappa) \otimes_{U(\widehat{\mfrak{g}}_{\rm st})} \!E \simeq U(\widehat{\mfrak{g}}_{\kappa,-})_{(f_\alpha)} \otimes_\C E
\end{align*}
for a $\mfrak{g}$-module $E$ and $\alpha \in \Delta_+$, which gives us an isomorphism
\begin{align*}
  (T_\alpha \circ \mathbb{M}_{\kappa,\mfrak{g}})(E) \simeq (U(\widehat{\mfrak{g}}_{\kappa,-})_{(f_\alpha)}/U(\widehat{\mfrak{g}}_{\kappa,-})) \otimes_\C E
\end{align*}
of $U(\widehat{\mfrak{g}}_{\kappa,-})$-modules. On the other hand, for $\alpha \in \Delta_+$ we have an isomorphism of $U(\widehat{\mfrak{g}}_{\kappa,-})$-modules
\begin{align*}
  (\mathbb{M}_{\kappa,\mfrak{g}} \circ T^\mfrak{g}_\alpha)(E) \simeq U(\widehat{\mfrak{g}}_{\kappa,-}) \otimes_\C T^\mfrak{g}_\alpha(E).
\end{align*}
Therefore, we define the isomorphism $\eta_{\alpha,E}$ of $U(\widehat{\mfrak{g}}_{\kappa,-})$-modules by
\begin{align*}
  f_\alpha^{-n}u \otimes v \mapsto \sum_{k=0}^\infty (-1)^k \binom{n+k-1}{k} \ad(f_\alpha)^k(u) \otimes f_\alpha^{-n-k}v
\end{align*}
with the inverse
\begin{align*}
  u \otimes f_\alpha^{-n}v \mapsto \sum_{k=0}^\infty \binom{n+k-1}{k} f_\alpha^{-n-k}\ad(f_\alpha)^k(u) \otimes v
\end{align*}
for $u \in U(\widehat{\mfrak{g}}_{\kappa,-})$, $v \in E$ and $n \in \N$. It is straightforward to check that $\eta_{\alpha,E}$ is an isomorphism of $U(\widehat{\mfrak{g}}_\kappa)$-modules.

The rest of the statement follows from the fact that $W^\mfrak{g}_\mfrak{p}(\lambda,\alpha) = T^\mfrak{g}_\alpha(M^\mfrak{g}_\mfrak{p}(\lambda))$ for $\lambda \in \Lambda^+(\mfrak{p})$ and $\alpha \in \Delta^\mfrak{u}_+$ by definition.}

\vspace{-2mm}

%%%%%%%%%%%%%%%%%%%%%%%%%%%%%%%%%%%%%%%%%%%%%%%%%%%%%%%%%%%%%%%%%%%%%%%%%%%%%%%%%%%%%%%%%%
%%%%%%%%%%%%%%%%%%%%%%%%%%%%%%%%%%%%%%%%%%%%%%%%%%%%%%%%%%%%%%%%%%%%%%%%%%%%%%%%%%%%%%%%%%

\section{Relaxed Wakimoto modules}
\label{sec-rel-wak}

We introduce a class of positive energy $\widehat{\mfrak{g}}_\kappa$-modules which we will call relaxed Wakimoto modules and also give a free field realization of relaxed Verma modules.

%%%%%%%%%%%%%%%%%%%%%%%%%%%%%%%%%%%%%%%%%%%%%%%%%%%%%%%%%%%%%%%%%%%%%%%%%%%%%%%%%%%%%%%%%%

\subsection{Feigin--Frenkel homomorphism}

Let $\mfrak{g}$ be a semisimple finite-dimensional Lie algebra. Let us consider a Borel subalgebra $\mfrak{b}$ of $\mfrak{g}$ with the nilradical $\mfrak{n}$, the opposite nilradical $\widebar{\mfrak{n}}$ and the Cartan subalgebra $\mfrak{h}$. Let $\{f_\alpha;\, \alpha \in \Delta_+\}$ be a root basis of the opposite nilradical $\widebar{\mfrak{n}}$. We denote by $\{x_\alpha;\, \alpha \in \Delta_+\}$ the linear coordinate functions on $\widebar{\mfrak{n}}$ with respect to the given basis of $\widebar{\mfrak{n}}$. Then the Weyl algebra $\eus{A}_{\widebar{\mfrak{n}}}$ of the vector space $\widebar{\mfrak{n}}$ is generated by $\{x_\alpha, \partial_{x_\alpha};\, \alpha \in \Delta_+\}$ together with the canonical commutation relations. Further, by \cite{Krizka-Somberg2017} there exists a homomorphism
\begin{align}
  \pi_\mfrak{g} \colon U(\mfrak{g}) \rarr \eus{A}_{\widebar{\mfrak{n}}} \otimes_\C U(\mfrak{h}) \label{eq:pi action}
\end{align}
of associative algebras uniquely determined by
\begin{align}
\pi_\mfrak{g}(a)= -\sum_{\alpha \in \Delta_+}\bigg[{\ad(u(x))e^{\ad(u(x))} \over e^{\ad(u(x))}-\id}\,(e^{-\ad(u(x))}a)_{\widebar{\mfrak{n}}}\bigg]_\alpha \partial_{x_\alpha} + (e^{-\ad(u(x))}a)_\mfrak{h} \label{eq:pi action general}
\end{align}
for $a \in \mfrak{g}$, where $[a]_\alpha$ denotes the $\alpha$-th coordinate of $a \in \widebar{\mfrak{n}}$ with respect to the basis $\{f_\alpha;\, \alpha \in \Delta_+\}$ of $\widebar{\mfrak{n}}$, $a_{\widebar{\mfrak{n}}}$ and $a_\mfrak{h}$ are the $\widebar{\mfrak{n}}$-part and $\mfrak{h}$-part of $a \in \mfrak{g}$ with respect to the triangular decomposition $\mfrak{g}=\widebar{\mfrak{n}} \oplus \mfrak{h} \oplus \mfrak{n}$, and the element $u(x) \in \C[\widebar{\mfrak{n}}] \otimes_\C \mfrak{g}$ is given by
\begin{align*}
u(x)=\sum_{\alpha \in \Delta_+} x_\alpha f_\alpha.
\end{align*}
Let us note that $\C[\widebar{\mfrak{n}}] \otimes_\C \mfrak{g}$ has the natural structure of a Lie algebra. Hence, we have a well-defined linear mapping $\ad(u(x)) \colon \C[\widebar{\mfrak{n}}] \otimes_\C \mfrak{g} \rarr \C[\widebar{\mfrak{n}}] \otimes_\C \mfrak{g}$.
\medskip

Let us recall that by Section \ref{subsec:Weyl vertex} we may assign to the vector space $\widebar{\mfrak{n}}$ the Weyl vertex algebra $\mcal{M}_{\widebar{\mfrak{n}}}$ generated by the fields
\begin{align*}a_\alpha(z)=\sum_{n\in \Z}\partial_{x_{\alpha,n}}z^{-n-1} \qquad \text{and}  \qquad
a^*_\alpha(z)=\sum_{n\in \Z}x_{\alpha, -n}z^{-n}
\end{align*}
for $\alpha\in \Delta_+$.
\medskip

The following theorem can be deduced from \cite[Theorem 5.1]{Frenkel2005}.
\medskip

\theorem{\label{thm:Wakimoto realization}
Let $\kappa$ be a $\mfrak{g}$-invariant symmetric bilinear form on $\mfrak{g}$. Then there exists a homomorphism
\begin{align*}
  w_{\kappa,\mfrak{g}} \colon \mcal{V}_\kappa(\mfrak{g}) \rarr \mcal{M}_{\widebar{\mfrak{n}}} \otimes_\C \mcal{V}_{\kappa - \kappa_c}(\mfrak{h})
\end{align*}
of $\N_0$-graded vertex algebras such that
\begin{align*}
\begin{aligned}
  w_{\kappa,\mfrak{g}}(e_\gamma(z)) & = - \sum_{\alpha \in \Delta_+} \normOrd{q^\gamma_\alpha(a_\beta^*(z))a_\alpha(z)} - (c_\gamma + (\kappa-\kappa_c) (e_\gamma,f_\gamma)) \partial_z a^*_\gamma(z) + a^*_\gamma(z)h_\gamma(z), \\
  w_{\kappa,\mfrak{g}}(h_\gamma(z)) & = \sum_{\alpha \in \Delta_+} \alpha(h_\gamma) \normOrd{a^*_\alpha(z)a_\alpha(z)} + h_\gamma(z), \\
  w_{\kappa,\mfrak{g}}(f_\gamma(z)) & = -a_\gamma(z) - \sum_{\alpha \in \Delta_+} \normOrd{p^\gamma_\alpha(a_\beta^*(z))a_\alpha(z)}
\end{aligned}
\end{align*}
for $\gamma \in \Pi$, where $c_\gamma \in \C$ are constants and the polynomials $p^\gamma_\alpha, q^\gamma_\alpha \in \C[\widebar{\mfrak{n}}]$ are given by
\begin{align*}
  p^\gamma_\alpha(x_\beta) &= \bigg[\bigg({\ad(u(x)) \over e^{\ad(u(x))}-\id} - \id\bigg)f_\gamma\bigg]_\alpha, \qquad q^\gamma_\alpha(x_\beta) = \bigg[{\ad(u(x))e^{\ad(u(x))} \over e^{\ad(u(x))}-\id}\,(e^{-\ad(u(x))}e_\gamma)_{\widebar{\mfrak{n}}}\bigg]_\alpha
\end{align*}
for $\gamma \in \Pi$ and $\alpha \in \Delta_+$.}

The link between the Feigin--Frenkel  homomorphism $ w_{\kappa,\mfrak{g}}$ of vertex algebras and the homomorphism
$\pi_\mfrak{g}$ of associative algebras is given by the following theorem.
\medskip

\theorem{Let $\kappa$ be a $\mfrak{g}$-invariant symmetric bilinear form on $\mfrak{g}$. Then the diagram
\begin{align} \label{eq:comm diag Zhu}
\bfig
\square(0,0)|alrb|<900,500>[\mcal{V}_\kappa(\mfrak{g})`\mcal{M}_{\widebar{\mfrak{n}}} \otimes_\C \mcal{V}_{\kappa-\kappa_c}(\mfrak{h})`U(\mfrak{g})`\eus{A}_{\widebar{\mfrak{n}}} \otimes_\C U(\mfrak{h});w_{\kappa,\mfrak{g}}`\pi_{\rm Zhu}`\pi_{\rm Zhu}`\pi_\mfrak{g}]
\efig
\end{align}
is commutative.}

\proof{As $w_{\kappa,\mfrak{g}} \colon  \mcal{V}_\kappa(\mfrak{g}) \rarr \mcal{M}_{\widebar{\mfrak{n}}} \otimes_\C \mcal{V}_{\kappa-\kappa_c}(\mfrak{h})$ is a homomorphism of $\N_0$-graded vertex algebras by Theorem \ref{thm:Wakimoto realization}, we obtain a homomorphism $\widebar{w}_{\kappa,\mfrak{g}} \colon U(\mfrak{g}) \rarr \eus{A}_{\widebar{\mfrak{n}}} \otimes_\C U(\mfrak{h})$ of the corresponding Zhu's algebras, since we have $A(\mcal{V}_\kappa(\mfrak{g})) \simeq U(\mfrak{g})$ and $A(\mcal{M}_{\widebar{\mfrak{n}}} \otimes_\C \mcal{V}_{\kappa-\kappa_c}(\mfrak{h})) \simeq \eus{A}_{\widebar{\mfrak{n}}} \otimes_\C U(\mfrak{h})$. Moreover, we have $\pi_{\rm Zhu} \circ w_{\kappa,\mfrak{g}} = \widebar{w}_{\kappa,\mfrak{g}} \circ \pi_{\rm Zhu}$. Therefore, we need only to show that $\widebar{w}_{\kappa,\mfrak{g}} = \pi_\mfrak{g}$. We may write
\begin{align*}
  \widebar{w}_{\kappa,\mfrak{g}}(e_\gamma) &= \widebar{w}_{\kappa,\mfrak{g}}(\pi_{\rm Zhu}(e_{\gamma,-1}\vac)) = \pi_{\rm Zhu}(w_{\kappa,\mfrak{g}}(e_{\gamma,-1}\vac)) \\
  &= \pi_{\rm Zhu}\!\bigg(\!-\sum_{\alpha \in \Delta_+} q^\gamma_\alpha(x_{\beta,0}) \partial_{x_{\alpha,-1}} - (c_\gamma + (\kappa-\kappa_c)(e_\gamma,f_\gamma))x_{\gamma,1} + x_{\gamma,0} h_{\gamma,-1}\vac\!\bigg) \\
  & = -\sum_{\alpha \in \Delta_+} q^\gamma_\alpha(x_\beta) \partial_{x_\alpha} + x_\gamma h_\gamma = \pi_\mfrak{g}(e_\gamma)
\intertext{and}
  \widebar{w}_{\kappa,\mfrak{g}}(f_\gamma) &= \widebar{w}_{\kappa,\mfrak{g}}(\pi_{\rm Zhu}(f_{\gamma,-1}\vac)) = \pi_{\rm Zhu}(w_{\kappa,\mfrak{g}}(f_{\gamma,-1}\vac)) \\
  &= \pi_{\rm Zhu}\!\bigg(\!-\partial_{x_{\gamma,-1}} - \sum_{\alpha \in \Delta_+} p^\gamma_\alpha(x_{\beta,0})\partial_{x_{\alpha,-1}}\!\bigg) = -\partial_{x_\gamma} - \sum_{\alpha \in \Delta_+} p^\gamma_\alpha(x_\beta) \partial_{x_\alpha} = \pi_\mfrak{g}(f_\gamma)
\end{align*}
for $\gamma \in \Pi$, which immediately implies that $\widebar{w}_{\kappa,\mfrak{g}} = \pi_\mfrak{g}$.}

Let us note that an explicit form of the homomorphism $w_{\kappa,\mfrak{g}}$ in Theorem \ref{thm:Wakimoto realization} is given only for the generators of the vertex algebra.
An alternative approach based on the Hamiltonian reduction of the WZNW model was considered in \cite{deBoer-Feher1997}. However, we will need explicit formulas for all elements of the opposite nilradical $\widebar{\mfrak{n}}$ which will be established in the next theorem.
\medskip

Let us denote by $\eus{P}^{\mfrak{g},\mfrak{b}}_{{\rm loc}}(z)$ the vector space of all polynomials in $a^*_\alpha(z)$ for $\alpha \in \Delta_+$ and by $\eus{F}^{\mfrak{g},\mfrak{b}}_{{\rm loc}}(z)$ the vector space of all differential polynomials in $a^*_\alpha(z)$ for $\alpha \in \Delta_+$.
We define a formal power series $u(z) \in \mfrak{g} \otimes_\C \smash{\eus{P}^{\mfrak{g},\mfrak{b}}_{{\rm loc}}(z)}$ by
\begin{align}
  u(z) = \sum_{\alpha \in \Delta_+} a_\alpha^*(z)f_\alpha. \label{eq:u(z) def}
\end{align}
Besides, the vector space $\mfrak{g} \otimes_\C \smash{\eus{F}^{\mfrak{g},\mfrak{b}}_{{\rm loc}}(z)}$ has the natural structure of a Lie algebra.
\medskip

We recall the following statement from \cite{Futorny-Krizka-Somberg2019}.
\medskip

\proposition{\label{prop:relation} \cite[Proposition 3.8]{Futorny-Krizka-Somberg2019}
We have the  identities
\begin{enumerate}[topsep=3pt,itemsep=0pt]
\item[1)]
\begin{align*}
  \bigg({{\rm d} \over {\rm d}t}_{|t=0} e^{\ad(u(z)+tx(z))}\!\bigg)e^{-\ad(u(z))}=\ad\!\bigg({e^{\ad(u(z))}-\id \over \ad(u(z))}\,x(z)\!\!\bigg)
\end{align*}
for $x(z) \in \mfrak{g} \otimes_\C \smash{\eus{F}^{\mfrak{g},\mfrak{b}}_{{\rm loc}}(z)}$,
\item[2)]
\begin{multline*}
   \bigg[{e^{\ad(u(z))}- \id \over \ad(u(z))}\,x(z),{e^{\ad(u(z))}- \id \over \ad(u(z))}\,y(z)\bigg] \\
   = {{\rm d} \over {\rm d}t}_{|t=0} {e^{\ad(u(z)+tx(z))}-\id \over \ad(u(z)+tx(z))}\,y(z) - {{\rm d} \over {\rm d}t}_{|t=0} {e^{\ad(u(z)+ty(z))}-\id \over \ad(u(z)+ty(z))}\,x(z)
\end{multline*}
for $x(z), y(z) \in \mfrak{g} \otimes_\C \smash{\eus{F}^{\mfrak{g},\mfrak{b}}_{{\rm loc}}(z)}$.
\end{enumerate}}

Let us introduce an element of $\widebar{\mfrak{n}} \otimes \smash{\eus{P}^{\mfrak{g},\mfrak{b}}_{{\rm loc}}(z)}$ by
\begin{align*}
  T(a,z) = {\ad(u(z)) \over e^{\ad(u(z))} - \id}\,a
\end{align*}
for $a \in \widebar{\mfrak{n}}$. Then we may write
\begin{align*}
  T(a,z) = \sum_{\alpha \in \Delta_+}  T_\alpha(a,z) f_\alpha,
\end{align*}
where $T_\alpha(a,z)=[T(a,z)]_\alpha$ for $\alpha \in \Delta_+$. Further, we define the linear mapping
\begin{align*}
  dT(a,z) \colon \widebar{\mfrak{n}} \otimes_\C \eus{F}^{\mfrak{g},\mfrak{b}}_{{\rm loc}}(z) \rarr \widebar{\mfrak{n}} \otimes_\C \eus{F}^{\mfrak{g},\mfrak{b}}_{{\rm loc}}(z)
\end{align*}
by
\begin{align*}
  dT(a,z)(x(z)) = {{\rm d} \over {\rm d}t}_{|t=0} {\ad(u(z)+tx(z)) \over e^{\ad(u(z)+tx(z))} - \id}\, a
\end{align*}
for $x(z) \in \widebar{\mfrak{n}} \otimes_\C \smash{\eus{F}^{\mfrak{g},\mfrak{b}}_{{\rm loc}}(z)}$. We have also
\begin{align*}
  dT(a,z) = \sum_{\alpha \in \Delta_+} dT_\alpha(a,z)f_\alpha,
\end{align*}
where $dT_\alpha(a,z) = [dT(a,z)]_\alpha$ for $\alpha \in \Delta_+$.
\medskip

\lemma{\label{lem:commutators}
We have
\begin{align*}
  [a_\alpha(z),u(w)]=f_\alpha \delta(z-w)
\end{align*}
and
\begin{align*}
  [a_\alpha(z),T(a,w)]=  dT(a,w)(f_\alpha) \delta(z-w)
\end{align*}
for $\alpha \in \Delta_+$ and $a \in \widebar{\mfrak{n}}$.}

\proof{By definition of $u(w)$ we have
\begin{align*}
  [a_\alpha(z),u(w)]= \sum_{\beta \in \Delta_+} [a_\alpha(z),a^*_\beta(w)]f_\beta=  \sum_{\beta \in \Delta_+} \delta_{\alpha,\beta} f_\beta \delta(z-w)= f_\alpha \delta(z-w).
\end{align*}
for $\alpha \in \Delta_+$. Further, we may write
\begin{align*}
  [a_\alpha(z), T(a,w)] = {{\rm d} \over {\rm d}t}_{|t=0} {\ad(u(w) + tf_\alpha) \over e^{\ad(u(w)+tf_\alpha)} - \id}\,a\, \delta(z-w) = dT(a,w)(f_\alpha) \delta(z-w)
\end{align*}
for $\alpha \in \Delta_+$ and $a \in \widebar{\mfrak{n}}$.}

\theorem{\label{thm-FF-hom}
Let $\kappa$ be a $\mfrak{g}$-invariant symmetric bilinear form on $\mfrak{g}$. Then we have
\begin{align}
  w_{\kappa,\mfrak{g}}(a(z))=-\sum_{\alpha\in \Delta_+} \normOrd{\!\bigg[{\ad(u(z)) \over e^{\ad(u(z))}-\id}\,a\bigg]_\alpha a_\alpha(z)} \label{eq:wakimoto homomorphism nilradical op}
\end{align}
for $a \in \widebar{\mfrak{n}}$ and
\begin{align}
   w_{\kappa,\mfrak{g}}(a(z))=  \sum_{\alpha\in \Delta_+} \normOrd{\,[\ad(u(z))(a)]_\alpha a_\alpha(z)} + a(z) \label{eq:wakimoto homomorphism cartan}
\end{align}
for $a \in \mfrak{h}$.}

\proof{From Theorem \ref{thm:Wakimoto realization} we know that the statement holds for all $a \in \mfrak{h}$ and for the root vectors $f_\gamma \in \widebar{\mfrak{n}}$ for $\gamma \in \Pi$. Hence, it is enough to show that $w_{\kappa,\mfrak{g}}$ given by \eqref{eq:wakimoto homomorphism nilradical op} gives rise to a homomorphism of vertex algebras from $\mcal{V}_\kappa(\widebar{\mfrak{n}})$ to $\mcal{M}_{\widebar{\mfrak{n}}} \otimes_\C \mcal{V}_{\kappa-\kappa_c}\!(\mfrak{h})$.

By using the relation \eqref{eq:comm relations vertex algebra}, we have
\begin{align*}
  w_{\kappa,\mfrak{g}}([a(z),b(w)]) &= w_{\kappa,\mfrak{g}}([a,b](w))\delta(z-w) + \kappa(a,b)\partial_w\delta(z-w) = w_{\kappa,\mfrak{g}}([a,b](w))\delta(z-w) \\
  & = - \sum_{\alpha \in \Delta_+} \normOrd{T_\alpha([a,b],w)a_\alpha(w)}\delta(z-w)
\end{align*}
for $a,b \in \widebar{\mfrak{n}}$. On the other hand, we get
\begin{align*}
  [w_{\kappa,\mfrak{g}}(a(z)), w_{\kappa,\mfrak{g}}(b(w))] &= \sum_{\alpha,\beta \in \Delta_+} [\normOrd{T_\alpha(a,z)a_\alpha(z)}, \normOrd{T_\beta(b,w)a_\beta(w)}] \\
  & = \sum_{\alpha,\beta \in \Delta_+} \normOrd{dT_\beta(b,w)(f_\alpha)T_\alpha(a,z)a_\beta(w)}\delta(z-w) \\
  & \quad - \sum_{\alpha, \beta \in \Delta_+} \normOrd{dT_\alpha(a,z)(f_\beta)T_\beta(b,w)a_\alpha(z)}\delta(z-w) \\
  & \quad - \sum_{\alpha,\beta \in \Delta_+} dT_\beta(b,w)(f_\alpha) dT_\alpha(a,z)(f_\beta) \partial_w \delta(z-w) \\
  & = \sum_{\alpha \in \Delta_+} \normOrd{\,(dT_\alpha(b,w)(T(a,w))- dT_\alpha(a,w)(T(b,w))) a_\alpha(w)}\delta(z-w) \\
  & \quad - \sum_{\alpha,\beta \in \Delta_+} dT_\beta(b,w)(f_\alpha) dT_\alpha(a,w)(f_\beta) \partial_w \delta(z-w) \\
  & \quad - \sum_{\alpha,\beta \in \Delta_+} dT_\beta(b,w)(f_\alpha) \partial_w dT_\alpha(a,w)(f_\beta) \delta(z-w),
\end{align*}
where we used the Wick theorem and Lemma \ref{lem:commutators}. Further, for $a, b \in \widebar{\mfrak{n}}$ we may write
\begin{align*}
  \sum_{\alpha,\beta \in \Delta_+} dT_\beta(b,w)(f_\alpha) dT_\alpha(a,w)(f_\beta) &= \sum_{\beta \in \Delta_+} dT_\beta(b,w)(dT(a,w)(f_\beta)) \\
  &= \tr_{\widebar{\mfrak{n}}} (dT(b,w) \circ dT(a,w))
\end{align*}
and
\begin{align*}
  \sum_{\alpha,\beta \in \Delta_+} dT_\beta(b,z)(f_\alpha) \partial_w dT_\alpha(a,z)(f_\beta) &= \sum_{\beta \in \Delta_+} dT_\beta(b,w)(\partial_w dT(a,w)(f_\beta)) \\
  &= \tr_{\widebar{\mfrak{n}}}(dT(b,w) \circ \partial_w dT(a,w)).
\end{align*}
Hence, we obtain
\begin{align*}
   [w_{\kappa,\mfrak{g}}(a(z)), w_{\kappa,\mfrak{g}}(b(w))] &= \sum_{\alpha \in \Delta_+} \normOrd{\,(dT_\alpha(b,w)(T(a,w))- dT_\alpha(a,w)(T(b,w))) a_\alpha(w)}\delta(z-w) \\
  & \quad - \tr_{\widebar{\mfrak{n}}} (dT(b,w) \circ dT(a,w))\partial_w\delta(z-w) \\
  & \quad - \tr_{\widebar{\mfrak{n}}}(dT(b,w) \circ \partial_w dT(a,w)) \delta(z-w)
\end{align*}
for $a,b \in \widebar{\mfrak{n}}$. Therefore, it is enough to show that
\begin{align*}
  T_\alpha([a,b],w) = \sum_{\alpha \in \Delta_+} (dT_\alpha(a,w)(T(b,w)) - dT_\alpha(b,w)(T(a,w)))
\end{align*}
for $\alpha \in \Delta_+$ and $a,b \in \widebar{\mfrak{n}}$, or equivalently
\begin{align*}
  T([a,b],w) = dT(a,w)(T(b,w)) - dT(b,w)(T(a,w))
\end{align*}
for $a,b \in \widebar{\mfrak{n}}$, and that
\begin{align*}
  \tr_{\widebar{\mfrak{n}}} (dT(b,w) \circ dT(a,w)) = 0, \qquad \tr_{\widebar{\mfrak{n}}}(dT(b,w) \circ \partial_w dT(a,w)) = 0
\end{align*}
for $a,b \in \widebar{\mfrak{n}}$.

A proof of the previous system of equations is a subject of the following lemmas, which then completes the proof of the present theorem.}

\lemma{We have
\begin{align*}
  \tr_{\widebar{\mfrak{n}}} (dT(b,z) \circ dT(a,z)) = 0 \qquad \text{and} \qquad \tr_{\widebar{\mfrak{n}}}(dT(b,z) \circ \partial_z dT(a,z)) = 0
\end{align*}
for $a,b \in \widebar{\mfrak{n}}$.}

\proof{Since $\widebar{\mfrak{n}}$ is a nilpotent Lie algebra, we have a canonical filtration
\begin{align*}
  \widebar{\mfrak{n}} = \widebar{\mfrak{n}}_0 \supset \widebar{\mfrak{n}}_1 \supset \dots \supset \widebar{\mfrak{n}}_n \supset \widebar{\mfrak{n}}_{n+1} = 0
\end{align*}
on $\widebar{\mfrak{n}}$ given by the lower central series of $\widebar{\mfrak{n}}$, i.e.\ $\widebar{\mfrak{n}}_k = [\widebar{\mfrak{n}}, \widebar{\mfrak{n}}_{k-1}]$ for $k \in \N$ with $\widebar{\mfrak{n}}_0 = \widebar{\mfrak{n}}$. Moreover, we have
\begin{align*}
  dT(a,z) \colon \widebar{\mfrak{n}}_k \otimes_\C \eus{F}^{\mfrak{g},\mfrak{b}}_{\rm loc}(z) \rarr \widebar{\mfrak{n}}_{k+1} \otimes_\C \eus{F}^{\mfrak{g},\mfrak{b}}_{\rm loc}(z)
\end{align*}
and
\begin{align*}
  \partial_zdT(a,z) \colon \widebar{\mfrak{n}}_k \otimes_\C \eus{F}^{\mfrak{g},\mfrak{b}}_{\rm loc}(z) \rarr \widebar{\mfrak{n}}_{k+1} \otimes_\C \eus{F}^{\mfrak{g},\mfrak{b}}_{\rm loc}(z)
\end{align*}
for $a \in \widebar{\mfrak{n}}$, which implies immediately that
\begin{align*}
  \tr_{\widebar{\mfrak{n}}} (dT(b,z) \circ dT(a,z)) = 0 \qquad \text{and} \qquad \tr_{\widebar{\mfrak{n}}}(dT(b,z) \circ \partial_z dT(a,z)) = 0
\end{align*}
for all $a, b \in \widebar{\mfrak{n}}$.}

\lemma{We have
\begin{align*}
  T([a,b],z) = dT(a,z)(T(b,z)) - dT(b,z)(T(a,z))
\end{align*}
for $a,b \in \widebar{\mfrak{n}}$.}

\proof{For $a,b \in \widebar{\mfrak{n}}$, we may write
\begin{align*}
  dT(a,z)(T(b,z)) &= {{\rm d} \over {\rm d}t}_{|t=0} {\ad(u(z)+tT(b,z)) \over e^{\ad(u(z)+tT(b,z))} - \id}\,a \\
  &= {{\rm d} \over {\rm d}t}_{|t=0} {\ad(u(z)+tT(b,z)) \over e^{\ad(u(z)+tT(b,z))} - \id} {e^{\ad(u(z))} - \id \over \ad(u(z))} {\ad(u(z)) \over e^{\ad(u(z))} - \id}\,a \\
  & = - {\ad(u(z)) \over e^{\ad(u(z))} - \id} {{\rm d} \over {\rm d}t}_{|t=0} {e^{\ad(u(z)+tT(b,z))} - \id \over \ad(u(z)+tT(b,z))}\, T(a,z),
\end{align*}
 which gives us
\begin{multline*}
  dT(a,z)(T(b,z)) - dT(b,z)(T(a,z)) \\ = - {\ad(u(z)) \over e^{\ad(u(z))} - \id} \bigg({{\rm d} \over {\rm d}t}_{|t=0} {e^{\ad(u(z)+tT(b,z))} - \id \over \ad(u(z)+tT(b,z))}\, T(a,z) - {{\rm d} \over {\rm d}t}_{|t=0} {e^{\ad(u(z)+tT(a,z))} - \id \over \ad(u(z)+tT(b,z))}\, T(b,z)\!\!\bigg).
\end{multline*}
By using Proposition \ref{prop:relation}, we obtain immediately
\begin{align*}
  dT(a,z)(T(b,z)) - dT(b,z)(T(a,z)) &= {\ad(u(z)) \over e^{\ad(u(z))} - \id} \bigg[{e^{\ad(u(z))}- \id \over \ad(u(z))}\,T(a,z),{e^{\ad(u(z))}- \id \over \ad(u(z))}\,T(b,z)\bigg] \\
  & =  {\ad(u(z)) \over e^{\ad(u(z))} - \id}\,[a,b] = T([a,b],z)
\end{align*}
for $a,b \in \widebar{\mfrak{n}}$.}

\vspace{-2mm}

%%%%%%%%%%%%%%%%%%%%%%%%%%%%%%%%%%%%%%%%%%%%%%%%%%%%%%%%%%%%%%%%%%%%%%%%%%%%%%%%%%%%%%%%%%

\subsection{Wakimoto functor}

We introduce a functor $\mathbb{W}_{\kappa,\mfrak{g}}$ from a certain subcategory of $\mcal{M}(\mfrak{g})$ to the category $\mcal{E}_+\!(\widehat{\mfrak{g}}_\kappa)$. For a generic level $\kappa$ this functor provides a free field realization of the induction functor $\mathbb{M}_{\kappa,\mfrak{g}}$.
\medskip

A straightforward reformulation of Theorem \ref{thm:Wakimoto realization} and Theorem \ref{thm-FF-hom} by using the completed Weyl algebra $\smash{\widetilde{\eus{A}}}_{\mcal{K}(\widebar{\mfrak{n}})}$ and the universal enveloping algebras $U(\widehat{\mfrak{g}}_\kappa)$, $U(\smash{\widehat{\mfrak{h}}}_{\kappa-\kappa_c})$ gives us the following statement.
\medskip

\theorem{\label{thm:Wakimoto realization operator}
Let $\kappa$ be a $\mfrak{g}$-invariant symmetric bilinear form on $\mfrak{g}$. Then there exists a homomorphism
\begin{align*}
  \pi_{\kappa,\mfrak{g}} \colon U(\widehat{\mfrak{g}}_\kappa) \rarr \smash{\widetilde{\eus{A}}}_{\mcal{K}(\widebar{\mfrak{n}})} \widehat{\otimes}_\C\, U(\widehat{\mfrak{h}}_{\kappa-\kappa_c})
\end{align*}
of associative $\Z$-graded algebras such that
\begin{align}
  \pi_{\kappa,\mfrak{g}}(e_\gamma(z)) = - \sum_{\alpha \in \Delta_+} \normOrd{q^\gamma_\alpha(a_\beta^*(z))a_\alpha(z)} - (c_\gamma + (\kappa-\kappa_c) (e_\gamma,f_\gamma)) \partial_z a^*_\gamma(z) + a^*_\gamma(z)h_\gamma(z)
\end{align}
for $\gamma \in \Pi$, where $c_\gamma \in \C$ are constants and the polynomials $q^\gamma_\alpha \in \C[\widebar{\mfrak{n}}]$ are given by
\begin{align*}
  q^\gamma_\alpha(x_\beta) = \bigg[{\ad(u(x))e^{\ad(u(x))} \over e^{\ad(u(x))}-\id}\,(e^{-\ad(u(x))}e_\gamma)_{\widebar{\mfrak{n}}}\bigg]_\alpha
\end{align*}
for $\gamma \in \Pi$ and $\alpha \in \Delta_+$. Further, we have
\begin{align}
  \pi_{\kappa,\mfrak{g}}(a(z))=-\sum_{\alpha\in \Delta_+} \normOrd{\!\bigg[{\ad(u(z)) \over e^{\ad(u(z))}-\id}\,a\bigg]_\alpha a_\alpha(z)} \label{eq:wakimoto homomorphism nilradical op operator}
\end{align}
for $a \in \widebar{\mfrak{n}}$ and
\begin{align}
   \pi_{\kappa,\mfrak{g}}(a(z))=  \sum_{\alpha\in \Delta_+} \normOrd{\,[\ad(u(z))(a)]_\alpha a_\alpha(z)} + a(z) \label{eq:wakimoto homomorphism cartan operator}
\end{align}
for $a \in \mfrak{h}$.}

\definition{Let $N$ be an $\eus{A}_{\widebar{\mfrak{n}}}$-module and $E$ be an $\mfrak{h}$-module. Then the smooth $\widehat{\mfrak{g}}_\kappa$-module
\begin{align*}
  \mathbb{W}_{\kappa,\mfrak{g}}(N \otimes_\C E) =  \mathbb{M}_{\mcal{K}(\widebar{\mfrak{n}})}(N) \otimes_\C \mathbb{M}_{\kappa-\kappa_c,\mfrak{h}}(E)
\end{align*}
is called the \emph{relaxed Wakimoto module} induced from the $\mfrak{g}$-module $N \otimes_\C E$.}

Let us note that it is easy to see that $\mathbb{M}_{\mcal{K}(\widebar{\mfrak{n}})}(N) \simeq \mathbb{L}_{\mcal{K}(\widebar{\mfrak{n}})}(N)$ for any $\eus{A}_{\widebar{\mfrak{n}}}$-module $N$ and that $\mathbb{M}_{\kappa-\kappa_c,\mfrak{h}}(E) \simeq \mathbb{L}_{\kappa-\kappa_c,\mfrak{h}}(E)$ for any $\mfrak{h}$-module $E$ provided $\kappa$ is a non-critical level.
\medskip

\theorem{\label{thm:Verma-Wakimoto homomorphism}
Let $N$ be an $\eus{A}_{\widebar{\mfrak{n}}}$-module and $E$ be an $\mfrak{h}$-module. Then $\mathbb{W}_{\kappa,\mfrak{g}}(N \otimes_\C E)$ is a positive energy $\widehat{\mfrak{g}}_\kappa$-module and the top degree component of $\mathbb{W}_{\kappa,\mfrak{g}}(N \otimes_\C E)$ is isomorphic to $N \otimes_\C E$ as a $\mfrak{g}$-module. Moreover, there is a non-trivial homomorphism
\begin{align*}
  \mathbb{M}_{\kappa,\mfrak{g}}(N \otimes_\C E) \rarr \mathbb{W}_{\kappa,\mfrak{g}}(N \otimes_\C E)
\end{align*}
of $\widehat{\mfrak{g}}_\kappa$-modules. If both $\mathbb{M}_{\kappa,\mfrak{g}}(N \otimes_\C E)$ and $\mathbb{W}_{\kappa,\mfrak{g}}(N \otimes_\C E)$ are simple $\widehat{\mfrak{g}}_\kappa$-modules, then the latter homomorphism is in fact an isomorphism.}

\proof{By definition we have that $\mathbb{M}_{\mcal{K}(\widebar{\mfrak{n}})}(N)$ and $\mathbb{M}_{\kappa-\kappa_c,\mfrak{h}}(E)$ are a positive energy $\eus{A}_{\mcal{K}(\widebar{\mfrak{n}})}$-module and a positive energy $\smash{\widehat{\mfrak{h}}}_{\kappa-\kappa_c}$-module, respectively. Hence, we have
\begin{align*}
  \mathbb{M}_{\mcal{K}(\widebar{\mfrak{n}})}(N) = \bigoplus_{n = 0}^\infty \mathbb{M}_{\mcal{K}(\widebar{\mfrak{n}})}(N)_n \qquad \text{and} \qquad \mathbb{M}_{\kappa-\kappa_c,\mfrak{h}}(E) = \bigoplus_{n=0}^\infty \mathbb{M}_{\kappa-\kappa_c,\mfrak{h}}(E)_n
\end{align*}
with $\mathbb{M}_{\mcal{K}(\widebar{\mfrak{n}})}(N)_0 \simeq N$ as $\eus{A}_{\widebar{\mfrak{n}}}$-modules and $\mathbb{M}_{\kappa-\kappa_c,\mfrak{h}}(E)_0 \simeq E$ as $\mfrak{h}$-modules, which gives us the gradation on $\mathbb{W}_{\kappa,\mfrak{g}}(N \otimes_\C E)$ defined by
\begin{align*}
  \mathbb{W}_{\kappa,\mfrak{g}}(N \otimes_\C E)_n = \bigoplus_{k=0}^n \mathbb{M}_{\mcal{K}(\widebar{\mfrak{n}})}(N)_{n-k} \otimes_\C \mathbb{M}_{\kappa-\kappa_c,\mfrak{h}}(E)_k
\end{align*}
for $n \in \N_0$. Further, from Theorem \ref{thm:Wakimoto realization operator} we obtain that the gradation on $\mathbb{W}_{\kappa,\mfrak{g}}(N \otimes_\C E)$ is compatible with the grading of $\widehat{\mfrak{g}}_\kappa$, which implies that $\mathbb{W}_{\kappa,\mfrak{g}}(N \otimes_\C E)$ is a positive energy $\widehat{\mfrak{g}}_\kappa$-module whose top degree component is isomorphic to $N \otimes_\C E$ as a $\mfrak{g}$-module. Hence, we have a homomorphism $N \otimes_\C E \rarr \mathbb{W}_{\kappa,\mfrak{g}}(N \otimes_\C E)$ of $\mfrak{g}$-modules which gives rise to a homomorphism
\begin{align*}
  \mathbb{M}_{\kappa,\mfrak{g}}(N \otimes_\C E) \rarr \mathbb{W}_{\kappa,\mfrak{g}}(N \otimes_\C E)
\end{align*}
of $\widehat{\mfrak{g}}_\kappa$-modules by the universal property of the generalized Verma module $\mathbb{M}_{\kappa,\mfrak{g}}(N \otimes_\C E)$.}

The homomorphism $\pi_\mfrak{g}$ of associative algebras gives rise to a bifunctor
\begin{align*}
  \Phi_\mfrak{g} \colon \mcal{M}(\eus{A}_{\widebar{\mfrak{n}}}) \times \mcal{M}(\mfrak{h}) \rarr \mcal{M}(\mfrak{g})
\end{align*}
defined by
\begin{align*}
  \Phi_\mfrak{g}(N,E) = N \otimes_\C E
\end{align*}
for an $\eus{A}_{\widebar{\mfrak{n}}}$-module $N$ and $\mfrak{h}$-module $E$, where the $\mfrak{g}$-module structure on $N \otimes_\C E$ is given by $\pi_\mfrak{g}$. This enable us to consider $\mathbb{W}_{\kappa,\mfrak{g}}$ as a functor from the category $\mcal{M}(\eus{A}_{\widebar{\mfrak{n}}}) \times \mcal{M}(\mfrak{h})$ to the category $\mcal{E}_+\!(\widehat{\mfrak{g}}_\kappa)$. We will call $\mathbb{W}_{\kappa,\mfrak{g}}$ the \emph{Wakimoto functor}.

Since the image of $\Phi_\mfrak{g}$ need not be a subcategory of $\mcal{M}(\mfrak{g})$, we denote by $\mcal{C}(\mfrak{g})$ the full subcategory of $\mcal{M}(\mfrak{g})$ consisting of $\mfrak{g}$-modules isomorphic to $N\otimes_\C E$ and their finite direct sums, where $N$ is an $\eus{A}_{\widebar{\mfrak{n}}}$-module, $E$ is a semisimple finite-dimensional $\mfrak{h}$-module and the $\mfrak{g}$-module structure on $N \otimes_\C E$ is given through the homomorphism $\pi_\mfrak{g}$. In the next, we will consider the Wakimoto functor $\mathbb{W}_{\kappa,\mfrak{g}}$ as a functor from the category $\mcal{C}(\mfrak{g})$ (there is a unique extension of $\mathbb{W}_{\kappa,\mfrak{g}}$ from $\im \Phi_\mfrak{g}$) to the category $\mcal{E}_+\!(\widehat{\mfrak{g}}_\kappa)$.
\medskip

\proposition{For $\alpha \in \Delta_+$ the twisting functor $T^\mfrak{g}_\alpha$ preserves the category $\mcal{C}(\mfrak{g})$.}

\proof{Let $M$ be a $\mfrak{g}$-module belonging to the category $\mcal{C}(\mfrak{g})$. Then we have $M \simeq \sum_{i=1}^n\! N_i \otimes_\C E_i$, where $N_i$ is an $\eus{A}_{\widebar{\mfrak{n}}}$-module and $E_i$ is a semisimple finite-dimensional $\mfrak{h}$-module for $i=1,2,\dots,n$, which gives us $T_\alpha(M) \simeq \sum_{i=1}^n\! T_\alpha(N_i \otimes_\C E_i)$. Therefore, we need to show that $T_\alpha(N \otimes_\C E)$ is also from the category $\mcal{C}(\mfrak{g})$ for an $\eus{A}_{\widebar{\mfrak{n}}}$-module $N$ and a semisimple finite-dimensional $\mfrak{h}$-module $E$.

Since $f_\alpha$ and $p_\alpha = \pi_\mfrak{g}(f_\alpha)$ for $\alpha \in \Delta_+$ are locally $\ad$-nilpotent regular elements in $U(\mfrak{g})$ and $\eus{A}_{\widebar{\mfrak{n}}}$, respectively, we obtain that
\begin{align*}
  (N \otimes_\C E)_{(f_\alpha)} \simeq N_{(p_\alpha)} \otimes_\C E
\end{align*}
as $\mfrak{g}$-modules, where the $\mfrak{g}$-module structure on $N_{(p_\alpha)} \otimes_\C E$ is defined through the homomorphism $\pi_\mfrak{g}$. Hence, we get $T_\alpha(N \otimes_\C E) \simeq (N_{(p_\alpha)}/N) \otimes_\C E$ for $\alpha \in \Delta_+$, which implies immediately that $T_\alpha(N \otimes_\C E)$ is an object of the category $\mcal{C}(\mfrak{g})$.}

Now, we may apply the Wakimoto functor $\mathbb{W}_{\kappa,\mfrak{g}}$ on Verma modules and Gelfand--Tsetlin modules to obtain a free field realization of the corresponding relaxed Verma modules.

Let $\C_\lambda$ be the $1$-dimensional $\mfrak{h}$-module given by a weight $\lambda \in \mfrak{h}^*$ and let $\alpha \in \Delta_+$ be a positive root. Then by \cite{Krizka-Somberg2017} and \cite{Futorny-Krizka2019b} there are isomorphisms
\begin{align}
  M^\mfrak{g}_\mfrak{b}(\lambda) \simeq \eus{A}_{\widebar{\mfrak{n}}}/\eus{I}_{\rm V} \otimes_\C \C_{\lambda+2\rho} \qquad \text{and} \qquad  W^\mfrak{g}_\mfrak{b}(\lambda,\alpha) \simeq \eus{A}_{\widebar{\mfrak{n}}}/\eus{I}_{{\rm GT},\alpha} \otimes_\C \C_{\lambda+2\rho} \label{eq:M and W realization}
\end{align}
of $\mfrak{g}$-modules, where $\eus{I}_{\rm V}$ and $\eus{I}_{{\rm GT},\alpha}$ are left ideals of $\eus{A}_{\widebar{\mfrak{n}}}$ defined by
\begin{align*}
  \eus{I}_{\rm V} = (x_\gamma;\, \gamma \in \Delta_+) \qquad \text{and} \qquad \eus{I}_{{\rm GT},\alpha} = (\partial_{x_\alpha},x_\gamma;\, \gamma \in \Delta_+ \setminus \{\alpha\}),
\end{align*}
which implies that the Verma module $M^\mfrak{g}_\mfrak{b}(\lambda)$ and $\alpha$-Gelfand--Tsetlin module $W^\mfrak{g}_\mfrak{b}(\lambda,\alpha)$ belong to the category $\mcal{C}(\mfrak{g})$. Hence, we may apply the Wakimoto functor $\mathbb{W}_{\kappa,\mfrak{g}}$ on these $\mfrak{g}$-modules and get the corresponding relaxed Wakimoto modules together with the canonical homomorphisms
\begin{align*}
   \varphi_{\rm V}^\lambda \colon \mathbb{M}_{\kappa,\mfrak{g}}(M^\mfrak{g}_\mfrak{b}(\lambda)) \rarr \mathbb{W}_{\kappa,\mfrak{g}}(M^\mfrak{g}_\mfrak{b}(\lambda)) \qquad \text{and} \qquad \varphi_{{\rm GT},\alpha}^\lambda \colon \mathbb{M}_{\kappa,\mfrak{g}}(W^\mfrak{g}_\mfrak{b}(\lambda,\alpha)) \rarr \mathbb{W}_{\kappa,\mfrak{g}}(W^\mfrak{g}_\mfrak{b}(\lambda,\alpha))
\end{align*}
of $\widehat{\mfrak{g}}_\kappa$-modules given by Theorem \ref{thm:Verma-Wakimoto homomorphism}.
Let us note also that $\mathbb{W}_{\kappa,\mfrak{g}}(M^\mfrak{g}_\mfrak{b}(\lambda))$ and $\mathbb{W}_{\kappa,\mfrak{g}}(W^\mfrak{g}_\mfrak{b}(\lambda,\alpha))$ are positive energy $\widehat{\mfrak{g}}_\kappa$-modules with the top degree component isomorphic to $M^\mfrak{g}_\mfrak{b}(\lambda)$ and $W^\mfrak{g}_\mfrak{b}(\lambda,\alpha)$, respectively.

Further, we introduce vector spaces
\begin{align}
  \eus{F}_{\widebar{\mfrak{n}}} = \C[\partial_{x_\gamma},\gamma \in \Delta_+] \qquad \text{and} \qquad
  \eus{F}_{\widebar{\mfrak{n}},\alpha} = \C[x_\alpha,\partial_{x_\gamma},\gamma \in \Delta_+ \setminus \{\alpha\}]
\end{align}
and endow them with the structure of $\eus{A}_{\widebar{\mfrak{n}}}$-modules by means of the canonical isomorphisms
\begin{align}
  \eus{A}_{\widebar{\mfrak{n}}}/\eus{I}_{\rm V} \simeq \eus{F}_{\widebar{\mfrak{n}}} \qquad \text{and} \qquad
  \eus{A}_{\widebar{\mfrak{n}}}/\eus{I}_{{\rm GT},\alpha} \simeq \eus{F}_{\widebar{\mfrak{n}},\alpha} \label{eq:Fock realization}
\end{align}
of vector spaces. For a bilinear form $\kappa$ on $\mfrak{h}$ and $\lambda \in \mfrak{h}^*$, we define a smooth $\smash{\widehat{\mfrak{h}}}_\kappa$-module
\begin{align*}
  \pi_\lambda^\kappa = \C[y_{\gamma,n},\, n \in \N, \gamma \in \Pi]
\end{align*}
by the formula
\begin{align*}
  h_{\gamma,n} = \begin{cases}
   \kappa(h,h) n \partial_{y_{\gamma,n}} &  \text{if $n \in \N$}, \\
   \lambda(h) & \text{if $n=0$}, \\
    y_{\gamma,-n} &  \text{if $n \in -\N$}
  \end{cases}
\end{align*}
for $\gamma \in \Pi$. In addition, we introduce smooth $\eus{A}_{\mcal{K}(\widebar{\mfrak{n}})}$-modules
\begin{align*}
  \eus{W}_{\widebar{\mfrak{n}}} &= \C[\partial_{x_{\gamma,-n}}, x_{\gamma,n},\, n \in \N, \gamma \in \Delta_+;\, \partial_{x_{\gamma,0}},\, \gamma \in \Delta_+], \\
  \eus{W}_{\widebar{\mfrak{n}},\alpha} & = \C[\partial_{x_{\gamma,-n}}, x_{\gamma,n},\, n \in \N, \gamma \in \Delta_+;\, x_{\alpha,0},\partial_{x_{\gamma,0}},\, \gamma \in \Delta_+ \setminus \{\alpha\}]
\end{align*}
for $\alpha \in \Delta_+$. Then it is easy to verify that
\begin{align}
  \mathbb{M}_{\kappa,\mfrak{h}}(\C_\lambda) \simeq \pi_\lambda^\kappa \label{eq:h module isomorphism}
\end{align}
as $\smash{\widehat{\mfrak{h}}}_\kappa$-modules and that
\begin{align}
  \mathbb{M}_{\mcal{K}(\widebar{\mfrak{n}})}(\eus{F}_{\widebar{\mfrak{n}}}) \simeq \eus{W}_{\widebar{\mfrak{n}}}, \qquad \mathbb{M}_{\mcal{K}(\widebar{\mfrak{n}})}(\eus{F}_{\widebar{\mfrak{n}},\alpha}) \simeq \eus{W}_{\widebar{\mfrak{n}},\alpha} \label{eq:A module isomorphism}
\end{align}
as $\eus{A}_{\mcal{K}(\widebar{\mfrak{n}})}$-modules.
\medskip

\theorem{\label{thm-realiz}
Let $\kappa$ be a $\mfrak{g}$-invariant symmetric bilinear form on $\mfrak{g}$. Then we have
\begin{align*}
  \mathbb{W}_{\kappa,\mfrak{g}}(M^\mfrak{g}_\mfrak{b}(\lambda)) \simeq \eus{W}_{\widebar{\mfrak{n}}} \otimes_\C \pi_{\lambda+2\rho}^{\kappa-\kappa_c} \qquad \text{and} \qquad
  \mathbb{W}_{\kappa,\mfrak{g}}(W^\mfrak{g}_\mfrak{b}(\lambda,\alpha)) \simeq \eus{W}_{\widebar{\mfrak{n}},\alpha} \otimes_\C \pi_{\lambda+2\rho}^{\kappa-\kappa_c}
\end{align*}
for $\lambda \in \mfrak{h}^*$ and $\alpha \in \Delta_+$.}

\proof{By definition of the Wakimoto functor $\mathbb{W}_{\kappa,\mfrak{g}}$ we have that
\begin{align*}
  \mathbb{W}_{\kappa,\mfrak{g}}(M^\mfrak{g}_\mfrak{b}(\lambda)) &\simeq \mathbb{M}_{\mcal{K}(\widebar{\mfrak{n}})}(\eus{F}_{\widebar{\mfrak{n}}}) \otimes_\C \mathbb{M}_{\kappa-\kappa_c,\mfrak{h}}(\C_{\lambda+2\rho}), \\ \mathbb{W}_{\kappa,\mfrak{g}}(W^\mfrak{g}_\mfrak{b}(\lambda,\alpha)) &\simeq \mathbb{M}_{\mcal{K}(\widebar{\mfrak{n}})}(\eus{F}_{\widebar{\mfrak{n}},\alpha}) \otimes_\C \mathbb{M}_{\kappa-\kappa_c,\mfrak{h}}(\C_{\lambda+2\rho})
\end{align*}
for $\lambda \in \mfrak{h}^*$ and $\alpha \in \Delta_+$, where we used \eqref{eq:M and W realization} and \eqref{eq:Fock realization}. Moreover, by using the isomorphisms \eqref{eq:h module isomorphism} and \eqref{eq:A module isomorphism} we get
\begin{align*}
  \mathbb{W}_{\kappa,\mfrak{g}}(M^\mfrak{g}_\mfrak{b}(\lambda)) \simeq \eus{W}_{\widebar{\mfrak{n}}} \otimes_\C \pi_{\lambda+2\rho}^{\kappa-\kappa_c} \qquad \text{and} \qquad
  \mathbb{W}_{\kappa,\mfrak{g}}(W^\mfrak{g}_\mfrak{b}(\lambda,\alpha)) \simeq \eus{W}_{\widebar{\mfrak{n}},\alpha} \otimes_\C \pi_{\lambda+2\rho}^{\kappa-\kappa_c}
\end{align*}
for $\lambda \in \mfrak{h}^*$ and $\alpha \in \Delta_+$.}

The link between the twisting functor $T_\alpha$ for $\alpha \in \Delta_+ \subset \smash{\widehat{\Delta}}^{\rm re}_+$ and the Wakimoto functor $\mathbb{W}_{\kappa,\mfrak{g}}$ is given in the following theorem.
\medskip

Let us introduce a smooth $\eus{A}_{\mcal{K}(\widebar{\mfrak{n}})}$-module
\begin{align*}
  \eus{E}_{\widebar{\mfrak{n}}} = \C[\partial_{x_{\gamma,-n}}, x_{\gamma,n},\, n \in \N, \gamma \in \Delta_+;\, x_{\gamma,0},\, \gamma \in \Delta_+]
\end{align*}
which is isomorphic to $\mathbb{M}_{\mcal{K}(\widebar{\mfrak{n}})}(\C[\widebar{\mfrak{n}}])$. The commutative algebra $\eus{E}_{\widebar{\mfrak{n}}}$ has the natural structure of a $\smash{\widehat{Q}}$-graded algebra, where $\smash{\widehat{Q}}$ is the affine root lattice. The gradation is uniquely determined by
\begin{align*}
  \deg x_{\gamma,n} = \gamma + n\delta \qquad \text{and} \qquad \deg \partial_{x_{\gamma,n}} = -\gamma-n\delta
\end{align*}
for $\gamma\in \Delta_+$ and $n \in \Z$. Hence, we have a direct sum decomposition
\begin{align*}
  \eus{E}_{\widebar{\mfrak{n}}} = \bigoplus_{\gamma \in \widehat{Q}} \eus{E}_{\widebar{\mfrak{n}},\gamma}.
\end{align*}
We say that a differential operator $P \in \smash{\widetilde{\eus{A}}}_{\mcal{K}(\widebar{\mfrak{n}})}$ has degree $\omega \in \smash{\widehat{Q}}$ if $P(\eus{E}_{\widebar{\mfrak{n}},\gamma}) \subset \eus{E}_{\widebar{\mfrak{n}},\gamma+\omega}$ for all $\gamma \in \smash{\widehat{Q}}$.
Further, for $\alpha \in \Delta_+$ we define a differential operator $p_\alpha \in \smash{\widetilde{\eus{A}}}_{\mcal{K}(\widebar{\mfrak{n}})}$ by
\begin{align*}
  p_\alpha = \pi_{\kappa,\mfrak{g}}(f_{\alpha,0}) = \Res_{z=0} \pi_{\kappa,\mfrak{g}}(f_\alpha(z)) = - \Res_{z=0} \sum_{\gamma\in \Delta_+} \normOrd{\!\bigg[{\ad(u(z)) \over e^{\ad(u(z))}-\id}\,f_\alpha\bigg]_\gamma a_\gamma(z)}.
\end{align*}
Then we may write
\begin{align*}
  p_\alpha = p_{\alpha,0} + p_{\alpha,1},
\end{align*}
where $p_{\alpha,0} = \pi_{\mfrak{g}}(f_{\alpha,0})$. The differential operators $p_{\alpha,0}$, $p_{\alpha,1}$ and $p_\alpha$ have degree $-\alpha$ for $\alpha \in \Delta_+$. We denote by $\eus{N}_\alpha$ the Lie subalgebra of $\smash{\widetilde{\eus{A}}}_{\mcal{K}(\widebar{\mfrak{n}})}$ generated by the set $\{\ad(p_{\alpha,0})^n(p_{\alpha,1});\, n \in \N_0\}$ and by $\eus{U}_\alpha$ the subalgebra of $\smash{\widetilde{\eus{A}}}_{\mcal{K}(\widebar{\mfrak{n}})}$ generated by $\eus{N}_\alpha$. Moreover, we have the decompositions
\begin{align*}
  \eus{N}_\alpha = \bigoplus_{n \in \N_0} \eus{N}_{\alpha,-n\alpha} \qquad \text{and} \qquad \eus{U}_\alpha = \bigoplus_{n \in \N_0} \eus{U}_{\alpha,-n\alpha}
\end{align*}
given by the degree of differential operators, where $\eus{N}_{\alpha,0} = 0$ and $\eus{U}_{\alpha,0} = \C$.

Let us consider an $\eus{A}_{\widebar{\mfrak{n}}}$-module $N$. Then the smooth $\eus{A}_{\mcal{K}(\widebar{\mfrak{n}})}$-module $\mathbb{M}_{\mcal{K}(\widebar{\mfrak{n}})}(N)$ has the canonical grading
\begin{align*}
  \mathbb{M}_{\mcal{K}(\widebar{\mfrak{n}})}(N) = \bigoplus_{n = 0}^\infty \mathbb{M}_{\mcal{K}(\widebar{\mfrak{n}})}(N)_n,
\end{align*}
where $\mathbb{M}_{\mcal{K}(\widebar{\mfrak{n}})}(N)_0 \simeq N$ as $\eus{A}_{\widebar{\mfrak{n}}}$-modules. Moreover, by definition of $p_\alpha$ we get immediately that the differential operator $p_{\alpha,0}$, $p_{\alpha,1}$ and $p_\alpha$ preserve the subspaces $\mathbb{M}_{\mcal{K}(\widebar{\mfrak{n}})}(N)_n$ for $n \in \N_0$. Besides, we have
\begin{align*}
  \mathbb{M}_{\mcal{K}(\widebar{\mfrak{n}})}(N) \simeq \eus{E}^-_{\widebar{\mfrak{n}}} \otimes_\C N,
\end{align*}
where
\begin{align*}
  \eus{E}^-_{\widebar{\mfrak{n}}} = \C[\partial_{x_{\gamma,-n}}, x_{\gamma,n},\, n \in \N, \gamma \in \Delta_+],
\end{align*}
which gives us
\begin{align*}
  \mathbb{M}_{\mcal{K}(\widebar{\mfrak{n}})}(N)_n \simeq \eus{E}^-_{\widebar{\mfrak{n}},n} \otimes_\C N
\end{align*}
for $n \in \N_0$. Since the Lie algebra $\eus{N}_\alpha$ preserves $\mathbb{M}_{\mcal{K}(\widebar{\mfrak{n}})}(N)_n$ for $n \in \N_0$, we have
\begin{align*}
   q(\eus{E}^-_{\widebar{\mfrak{n}},n,\gamma} \otimes_\C N) \subset \eus{E}^-_{\widebar{\mfrak{n}},n,\gamma - j\alpha}\! \otimes_\C N
\end{align*}
for $q \in \eus{N}_{\alpha,-j\alpha}$, where $\eus{E}^-_{\widebar{\mfrak{n}},n,\gamma} = \eus{E}^-_{\widebar{\mfrak{n}},n} \cap \eus{E}_{\widebar{\mfrak{n}},\gamma}$. However, the vector space $\eus{E}^-_{\widebar{\mfrak{n}},n}$ is finite dimensional, which implies that there exists $j_{n,\gamma} \in \N_0$ such that $\smash{\eus{E}^-_{\widebar{\mfrak{n}},n,\gamma - j\alpha}} = 0$ for $j > j_{n,\gamma}$. Therefore, for each vector $v \in \mathbb{M}_{\mcal{K}(\widebar{\mfrak{n}})}(N)$ there exists an integer $n_v\in \N_0$ such that $\eus{U}_{\alpha,-n\alpha}v =0$ for all $n > n_v$.

Since $p_\alpha$ and $p_{\alpha,0}$ for $\alpha \in \Delta_+$ are locally $\ad$-nilpotent regular elements of the completed Weyl algebra $\smash{\widetilde{\eus{A}}_{\mcal{K}(\widebar{\mfrak{n}})}}$, we may construct localized modules $\mathbb{M}_{\mcal{K}(\widebar{\mfrak{n}})}(N)_{(p_\alpha)}$ and $\mathbb{M}_{\mcal{K}(\widebar{\mfrak{n}})}(N)_{(p_{\alpha,0})}$. By a completely similar way as in \cite[Lemma 4.6]{Futorny-Krizka2019b}, we can show that the element $p_{\alpha,1}\smash{p_{\alpha,0}^{-1}} \in \smash{\widetilde{\eus{A}}_{\mcal{K}(\widebar{\mfrak{n}})}}$ acts locally nilpotently on $\mathbb{M}_{\mcal{K}(\widebar{\mfrak{n}})}(N)_{(p_{\alpha,0})}$.
\medskip

\proposition{\label{prop:localization isomorphism}
Let $\alpha \in \Delta_+$ and let $N$ be an $\eus{A}_{\widebar{\mfrak{n}}}$-module. Then the linear mapping
\begin{align*}
  \Phi_\alpha \colon \mathbb{M}_{\mcal{K}(\widebar{\mfrak{n}})}(N)_{(p_\alpha)} \rarr \mathbb{M}_{\mcal{K}(\widebar{\mfrak{n}})}(N)_{(p_{\alpha,0})}
\end{align*}
defined by
\begin{align*}
  \Phi_\alpha(p_\alpha^{-n}v) = \varphi_\alpha^n(v)
\end{align*}
for $n \in \N_0$ and $v \in \mathbb{M}_{\mcal{K}(\widebar{\mfrak{n}})}(N)$, where the linear mapping $\varphi_\alpha \colon \mathbb{M}_{\mcal{K}(\widebar{\mfrak{n}})}(N)_{(p_{\alpha,0})} \rarr \mathbb{M}_{\mcal{K}(\widebar{\mfrak{n}})}(N)_{(p_{\alpha,0})}$
is given through
\begin{align*}
  \varphi_\alpha(v) = p_{\alpha,0}^{-1} \sum_{k=0}^\infty (-1)^k(p_{\alpha,1} p_{\alpha,0}^{-1})^kv
\end{align*}
for $v \in \mathbb{M}_{\mcal{K}(\widebar{\mfrak{n}})}(N)$, is an isomorphism of $\smash{(\widetilde{\eus{A}}_{\mcal{K}(\widebar{\mfrak{n}})})_{(p_\alpha)}}$-modules.}

\proof{Since the element $p_{\alpha,1}p_{\alpha,0}^{-1}$ acts locally nilpotently on $\mathbb{M}_{\mcal{K}(\widebar{\mfrak{n}})}(N)_{(p_{\alpha,0})}$, the linear mapping $\varphi_\alpha$ is well defined. By an analogous way as in \cite[Lemma 4.7, Lemma 4.8]{Futorny-Krizka2019b}, we can show that $\Phi_\alpha$ is a homomorphism of $\smash{(\widetilde{\eus{A}}_{\mcal{K}(\widebar{\mfrak{n}})})_{(p_\alpha)}}$-modules. Hence, we only need to prove that the linear mapping $\Phi_\alpha$ is injective and surjective. Since we may write
\begin{align*}
  \Phi_\alpha(qp_\alpha^{-n}v) = q\varphi_\alpha^n(v) = q p_{\alpha,0}^{-n}v
\end{align*}
for $n \in \N_0$, $v \in N$ and $q \in \eus{E}^-_{\widebar{\mfrak{n}}}$, we have the surjectivity of $\Phi_\alpha$. To prove the injectivity of $\Phi_\alpha$ let us assume that $\Phi_\alpha(v)=0$ for some $v \in \mathbb{M}_{\mcal{K}(\widebar{\mfrak{n}})}(N)_{(p_\alpha)}$. Then there exists an integer $n \in \N_0$ such that $p_\alpha^nv \in \mathbb{M}_{\mcal{K}(\widebar{\mfrak{n}})}(N) \subset \mathbb{M}_{\mcal{K}(\widebar{\mfrak{n}})}(N)_{(p_\alpha)}$. Hence, we have
\begin{align*}
    0 = p_\alpha^n\Phi_\alpha(v) = \Phi_\alpha(p_\alpha^nv) = p_\alpha^nv \in \mathbb{M}_{\mcal{K}(\widebar{\mfrak{n}})}(N) \subset \mathbb{M}_{\mcal{K}(\widebar{\mfrak{n}})}(N)_{(p_{\alpha,0})},
\end{align*}
which gives us $v=0$. Therefore, the linear mapping $\Phi_\alpha$ is an isomorphism.}

\theorem{\label{thm:twisting functor intertwining Wakimoto}
Let $\alpha \in \Delta_+ \subset \smash{\widehat{\Delta}}^{\rm re}_+$. Then there exists a natural isomorphism
\begin{align*}
  \eta_\alpha \colon T_\alpha \circ \mathbb{W}_{\kappa,\mfrak{g}} \rarr \mathbb{W}_{\kappa,\mfrak{g}} \circ\, T_\alpha^\mfrak{g}
\end{align*}
of functors, where $T_\alpha^\mfrak{g} \colon \mcal{C}(\mfrak{g}) \rarr \mcal{C}(\mfrak{g})$ is the twisting functor for $\mfrak{g}$ assigned to $\alpha$. In particular, we have
\begin{align*}
  T_\alpha(\mathbb{W}_{\kappa,\mfrak{g}}(M^\mfrak{g}_\mfrak{b}(\lambda))) \simeq \mathbb{W}_{\kappa,\mfrak{g}}(W^\mfrak{g}_\mfrak{b}(\lambda,\alpha))
\end{align*}
for $\lambda \in \mfrak{h}^*$ and $\alpha \in \Delta_+$.}

\proof{Let $N$ be an $\eus{A}_{\widebar{\mfrak{n}}}$-module and $E$ be a semisimple finite-dimensional $\mfrak{h}$-module. Then by using the definition of the Wakimoto functor $\mathbb{W}_{\kappa,\mfrak{g}}$ we obtain an isomorphism
\begin{align*}
  \mathbb{W}_{\kappa,\mfrak{g}}(N \otimes_\C E)_{(f_{\alpha,0})} \simeq \mathbb{M}_{\mcal{K}(\widebar{\mfrak{n}})}(N)_{(p_\alpha)} \otimes_\C \mathbb{M}_{\kappa-\kappa_c,\mfrak{h}}(E)
\end{align*}
of $U(\widehat{\mfrak{g}}_\kappa)$-modules, where the $U(\widehat{\mfrak{g}}_\kappa)$-module structure on $\mathbb{M}_{\mcal{K}(\widebar{\mfrak{n}})}(N)_{(p_\alpha)} \otimes_\C \mathbb{M}_{\kappa-\kappa_c,\mfrak{h}}(E)$ is given via the homomorphism $\pi_{\kappa,\mfrak{g}}$. Further, from Proposition \ref{prop:localization isomorphism} we have that
\begin{align*}
  \mathbb{M}_{\mcal{K}(\widebar{\mfrak{n}})}(N)_{(p_\alpha)} \simeq \mathbb{M}_{\mcal{K}(\widebar{\mfrak{n}})}(N)_{(p_{\alpha,0})} \simeq \mathbb{M}_{\mcal{K}(\widebar{\mfrak{n}})}(N_{(p_{\alpha,0})})
\end{align*}
as $\smash{\widetilde{\eus{A}}_{\mcal{K}(\widebar{\mfrak{n}})}}$-modules, which gives us
\begin{align*}
  \mathbb{W}_{\kappa,\mfrak{g}}(N \otimes_\C E)_{(f_{\alpha,0})} &\simeq \mathbb{M}_{\mcal{K}(\widebar{\mfrak{n}})}(N_{(p_{\alpha,0})}) \otimes_\C \mathbb{M}_{\kappa-\kappa_c,\mfrak{h}}(E) \simeq \mathbb{W}_{\kappa,\mfrak{g}}(N_{(p_{\alpha,0})} \otimes_\C E) \\
  &\simeq \mathbb{W}_{\kappa,\mfrak{g}}((N \otimes_\C E)_{(f_\alpha)})
\end{align*}
as $U(\widehat{\mfrak{g}}_\kappa)$-modules. Hence, by definition of $T_\alpha$ and $T_\alpha^\mfrak{g}$ we get
\begin{align*}
  T_\alpha(\mathbb{W}_{\kappa,\mfrak{g}}(N \otimes_\C E)) \simeq \mathbb{W}_{\kappa,\mfrak{g}}(T_\alpha^\mfrak{g}(N \otimes_\C E))
\end{align*}
for $\alpha \in \Delta_+$. Moreover, it is obvious that this isomorphism is natural. The rest of the statement follows immediately.}

With the help of Theorem \ref{thm:twisting functor intertwining Wakimoto} we can establish a relation between the relaxed Verma module $\mathbb{M}_{\kappa,\mfrak{g}}(W^\mfrak{g}_\mfrak{b}(\lambda,\alpha))$ and the relaxed Wakimoto module $\mathbb{W}_{\kappa,\mfrak{g}}(W^\mfrak{g}_\mfrak{b}(\lambda,\alpha))$.
\medskip

\corollary{\label{thm:alpha Verma-Wakimoto isomorphism}
Let $\lambda \in \mfrak{h}^*$ and let us assume that the Verma module $\mathbb{M}_{\kappa,\mfrak{g}}(M^\mfrak{g}_\mfrak{b}(\lambda))$, or equivalently the Wakimoto module $\mathbb{W}_{\kappa,\mfrak{g}}(M^\mfrak{g}_\mfrak{b}(\lambda))$, is a simple $\widehat{\mfrak{g}}_\kappa$-module. Then the relaxed Verma module $\mathbb{M}_{\kappa,\mfrak{g}}(W^\mfrak{g}_\mfrak{b}(\lambda,\alpha))$
is isomorphic to the relaxed Wakimoto module $\mathbb{W}_{\kappa,\mfrak{g}}(W^\mfrak{g}_\mfrak{b}(\lambda,\alpha))$ for $\alpha \in \Delta_+$.}

\proof{If the Verma module $\mathbb{M}_{\kappa,\mfrak{g}}(M^\mfrak{g}_\mfrak{b}(\lambda))$, or equivalently the Wakimoto module $\mathbb{W}_{\kappa,\mfrak{g}}(M^\mfrak{g}_\mfrak{b}(\lambda))$, is a simple $\widehat{\mfrak{g}}_\kappa$-module, then  $\mathbb{M}_{\kappa,\mfrak{g}}(M^\mfrak{g}_\mfrak{b}(\lambda)) \simeq \mathbb{W}_{\kappa,\mfrak{g}}(M^\mfrak{g}_\mfrak{b}(\lambda))$. By applying the twisting functor $T_\alpha$ for $\alpha \in \Delta_+$ on both sides of the isomorphism, we obtain
\begin{align*}
 \mathbb{M}_{\kappa,\mfrak{g}}(W^\mfrak{g}_\mfrak{b}(\lambda,\alpha)) \simeq  T_\alpha(\mathbb{M}_{\kappa,\mfrak{g}}(M^\mfrak{g}_\mfrak{b}(\lambda))) \simeq T_\alpha(\mathbb{W}_{\kappa,\mfrak{g}}(M^\mfrak{g}_\mfrak{b}(\lambda))) \simeq \mathbb{W}_{\kappa,\mfrak{g}}(W^\mfrak{g}_\mfrak{b}(\lambda,\alpha)),
\end{align*}
where the first isomorphism follow from Theorem \ref{thm:twisting functor intertwining} and the last isomorphism is a consequence of Theorem \ref{thm:twisting functor intertwining Wakimoto}.}

Let us note that Corollary \ref{thm:alpha Verma-Wakimoto isomorphism} provides a free field realization of the relaxed Verma module $\mathbb{M}_{\kappa,\mfrak{g}}(W^\mfrak{g}_\mfrak{b}(\lambda,\alpha))$ whenever $\mathbb{M}_{\kappa,\mfrak{g}}(M^\mfrak{g}_\mfrak{b}(\lambda))$, or equivalently $\mathbb{W}_{\kappa,\mfrak{g}}(M^\mfrak{g}_\mfrak{b}(\lambda))$, is a simple $\widehat{\mfrak{g}}_\kappa$-module. On the other hand, this does not give any information for the critical level $\kappa_c$, since $\mathbb{M}_{\kappa_c,\mfrak{g}}(M^\mfrak{g}_\mfrak{b}(\lambda))$ is never a simple $\widehat{\mfrak{g}}_{\kappa_c}$-module. However, applying \cite[Proposition 9.5.1]{Frenkel2007-book} to the longest element of the Weyl group of $\mfrak{g}$ we have that the Verma module $\mathbb{M}_{\kappa_c,\mfrak{g}}(M^\mfrak{g}_\mfrak{b}(\lambda))$ and the Wakimoto module $\mathbb{W}_{\kappa_c,\mfrak{g}}(M^\mfrak{g}_\mfrak{b}(\lambda))$ are isomorphic if $\langle \lambda +\rho, \alpha^\vee \rangle \notin -\N$ for all $\alpha \in \Delta_+$.
Hence, we have the following extension of Corollary \ref{thm:alpha Verma-Wakimoto isomorphism} in the case of the critical level.
\medskip

\corollary{\label{thm:alpha Verma-Wakimoto isomorphism critical}
Let $\lambda \in \mfrak{h}^*$ satisfy $\langle \lambda +\rho, \alpha^\vee \rangle \notin -\N$ for all $\alpha \in \Delta_+$, i.e.\ $\lambda$ is dominant weight. Then the relaxed Verma module $\mathbb{M}_{\kappa_c, \mfrak{g}}(W^\mfrak{g}_\mfrak{b}(\lambda,\alpha))$ is isomorphic to the relaxed Wakimoto module $\mathbb{W}_{\kappa_c,\mfrak{g}}(W^\mfrak{g}_\mfrak{b}(\lambda,\alpha))$ for $\alpha \in \Delta_+$.}

%\proof{If $\lambda \in \mfrak{h}^*$ satisfies the condition $\langle \lambda +\rho, \alpha^\vee \rangle \notin -\N$ for all $\alpha \in \Delta_+$, then we have that $\mathbb{M}_{\kappa_c,\mfrak{g}}(M^\mfrak{g}_\mfrak{b}(\lambda)) \simeq \mathbb{W}_{\kappa_c,\mfrak{g}}(M^\mfrak{g}_\mfrak{b}(\lambda))$. By applying the twisting functor $T_\alpha$ for $\alpha \in \Delta_+$, we obtain
%\begin{align*}
 %\mathbb{M}_{\kappa_c,\mfrak{g}}(W^\mfrak{g}_\mfrak{b}(\lambda,\alpha)) \simeq  T_\alpha(\mathbb{M}_{\kappa_c,\mfrak{g}}(M^\mfrak{g}_\mfrak{b}(\lambda))) \simeq T_\alpha(\mathbb{W}_{\kappa_c,\mfrak{g}}(M^\mfrak{g}_\mfrak{b}(\lambda))) \simeq \mathbb{W}_{\kappa_c,\mfrak{g}}(W^\mfrak{g}_\mfrak{b}(\lambda,\alpha)),
%\end{align*}
%where the first and last isomorphism follow from Theorem \ref{thm:twisting functor intertwining} and Theorem \ref{thm:twisting functor intertwining Wakimoto}, respectively, which should be shown.}

\vspace{-2mm}

%%%%%%%%%%%%%%%%%%%%%%%%%%%%%%%%%%%%%%%%%%%%%%%%%%%%%%%%%%%%%%%%%%%%%%%%%%%%%%%%%%%%%%%%%%
%%%%%%%%%%%%%%%%%%%%%%%%%%%%%%%%%%%%%%%%%%%%%%%%%%%%%%%%%%%%%%%%%%%%%%%%%%%%%%%%%%%%%%%%%%

\section{Positive energy representations of $\mcal{L}_\kappa(\mfrak{g})$}
\label{sec-pos-energy}

In this section we describe families of positive energy representations of the simple affine vertex algebra $\mcal{L}_\kappa(\mfrak{g})$ of an admissible level $\kappa$ associated to a simple Lie algebra $\mfrak{g}$.

%%%%%%%%%%%%%%%%%%%%%%%%%%%%%%%%%%%%%%%%%%%%%%%%%%%%%%%%%%%%%%%%%%%%%%%%%%%%%%%%%%%%%%%%%%

\subsection{Admissible representations}

Let $\mfrak{g}$ be a complex simple Lie algebra and let $\kappa$ be a $\mfrak{g}$-invariant
symmetric bilinear form on $\mfrak{g}$. Since $\mfrak{g}$ is a simple Lie algebra, we have $\kappa = k\kappa_0$ for $k \in \C$, where $\kappa_0$ is the normalized $\mfrak{g}$-invariant symmetric bilinear form on $\mfrak{g}$ satisfying
\begin{align*}
  \kappa_\mfrak{g} = 2h^\vee \kappa_0.
\end{align*}
Let $\widehat{\mfrak{g}}_\kappa$ be the affine Kac--Moody algebra associated to the Lie algebra $\mfrak{g}$ of level $\kappa$. For $\lambda \in \widehat{\mfrak{h}}^*$, we define its integral root system $\widehat{\Delta}(\lambda)$ by
\begin{align*}
 \widehat{\Delta}(\lambda) =\{\alpha\in \widehat{\Delta}^{\rm re};\, \bra \lambda+\widehat{\rho}, \alpha^\vee \ket \in \Z\},
\end{align*}
where $\widehat{\rho}=\rho+h^\vee \Lambda_0$. Further, let $\widehat{\Delta}(\lambda)_+ = \widehat{\Delta}(\lambda) \cap \widehat{\Delta}^{\rm re}_+$ be the set of positive roots of $\widehat{\Delta}(\lambda)$ and $\smash{\widehat{\Pi}}(\lambda) \subset \smash{\widehat{\Delta}}(\lambda)_+$ be the set of simple roots. Then we say that a weight $\lambda \in \smash{\widehat{\mfrak{h}}^*}$ is \emph{admissible} (\cite{Kac-Wakimoto1989}) provided
\begin{enumerate}[topsep=0pt,itemsep=0pt,parsep=0pt]
 \item[i)] $\lambda$ is \emph{regular dominant}, that is $\bra \lambda + \widehat{\rho}, \alpha^\vee \ket \notin -\N_0$ for all $\alpha\in \widehat{\Delta}^{\rm re}_+$;
 \item[ii)] the $\Q$-span of $\smash{\widehat{\Delta}}(\lambda)$ contains $\smash{\widehat{\Delta}^{\rm re}}$.
\end{enumerate}
In particular, if $\lambda=k\Lambda_0$ is an admissible weight for $k \in \C$, then $k$ is called an \emph{admissible number}. The admissible numbers were described in \cite{Kac-Wakimoto1989,Kac-Wakimoto2008} as follows. The complex number $k \in \C$ is admissible if and only if
\begin{align*}
  k+h^\vee ={p \over q} \ \text{with } p,q\in \N,\ (p,q)=1,\
  p\geq
    \begin{cases}
      h^\vee & \text{if $(r^\vee,q)=1$}, \\
      h & \text{if $(r^\vee,q)=r^\vee$},
    \end{cases}
\end{align*}
where $r^\vee$ is the lacing number of $\mfrak{g}$, i.e.\ the maximal number of edges in the Dynkin diagram of the Lie algebra $\mfrak{g}$. Since admissibility of a number $k \in \C$  depends only on $\mfrak{g}$, we shall say that $k$ is an admissible number for $\mfrak{g}$.

Further, let us assume that $k \in \Q$ is an admissible number for $\mfrak{g}$. We say that a $\mfrak{g}$-module $E$ is \emph{admissible of level $k$} if $\mathbb{L}_{k\kappa_0,\mfrak{g}}(E)$ is an $\mcal{L}_{k\kappa_0}\!(\mfrak{g})$-module, or equivalently if $E$ is an $A(\mcal{L}_{k\kappa_0}\!(\mfrak{g}))$-module. In particular, the simple highest weight $\mfrak{g}$-module $L^\mfrak{g}_\mfrak{b}(0)$ with zero highest weight is an admissible $\mfrak{g}$-module of level $k$. Moreover, since $A(\mcal{L}_{k\kappa_0}\!(\mfrak{g})) \simeq U(\mfrak{g})/I_k$, where $I_k$ is a two-sided ideal of $U(\mfrak{g})$, we obtain that a $\mfrak{g}$-module $E$ is admissible of level $k$ if and only if the ideal $I_k$ is contained in the annihilator $\Ann_{U(\mfrak{g})}\!E$.
\medskip

Admissible highest weight $\mfrak{g}$-modules of level $k$ were classified in \cite{Arakawa2016} as follows. Let ${\rm Pr}_k$ be the set of admissible weights $\lambda \in \smash{\widehat{\mfrak{h}}^*}$ of level $k$ such that there is an element $y \in \smash{\widetilde{W}}$ of the extended affine Weyl group $\smash{\widetilde{W}}$ of $\mfrak{g}$ satisfying $\smash{\widehat{\Delta}(\lambda)} = \smash{y(\widehat{\Delta}(k\Lambda_0))}$. Further, let us define the subset
\begin{align*}
  \widebar{{\rm Pr}}_k = \{\widebar{\lambda};\, \lambda \in {\rm Pr}_k\}
\end{align*}
of $\mfrak{h}^*$, where $\widebar{\lambda} \in \mfrak{h}^*$ denotes the canonical projection of $\lambda \in \smash{\widehat{\mfrak{h}}^*}$ to $\mfrak{h}^*$.
\medskip

\theorem{\cite{Arakawa2016}\label{thm-hweight}
Let $k \in \Q$ be an admissible number for $\mfrak{g}$. Then the simple highest weight $\mfrak{g}$-module $L^\mfrak{g}_\mfrak{b}(\lambda)$ with highest weight $\lambda \in \mfrak{h}^*$ is admissible of level $k$ if and only if $\lambda \in \widebar{{\rm Pr}}_k$.}

For $\lambda \in \mfrak{h}^*$, let $L^\mfrak{g}_\mfrak{b}(\lambda)$ be the simple $\mfrak{g}$-module with highest weight $\lambda$ and
\begin{align}
  J_\lambda = \Ann_{U(\mfrak{g})}\! L^\mfrak{g}_\mfrak{b}(\lambda)
\end{align}
the corresponding primitive ideal of $U(\mfrak{g})$. A theorem of Duflo \cite{Duflo1977} states that for any primitive ideal $I$ of $U(\mfrak{g})$ there exists $\lambda \in \mfrak{h}^*$ such that $I=J_\lambda$. This implies that a simple weight $\mfrak{g}$-module $E$ is admissible of level $k$ if and only if $\Ann_{U(\mfrak{g})}\!E=J_\lambda$ for some $\lambda \in \widebar{{\rm Pr}}_k$. Besides, for $\lambda, \mu \in \widebar{{\rm Pr}}_k$ we have $J_\lambda=J_\mu$ if and only if there exists $w\in W$ such that $ \mu=w\cdot \lambda$ (see Proposition 2.4 in \cite{Arakawa-Futorny-Ramirez2017}). Hence, we may define an equivalence relation on $\widebar{{\rm Pr}}_k$ by
\begin{align*}
  \lambda \sim \mu \Longleftrightarrow \text{there exits $w \in W$ such that $\mu = w \cdot \lambda$}
\end{align*}
and set
\begin{align*}
  [\widebar{{\rm Pr}}_k] = \widebar{{\rm Pr}}_k / \sim
\end{align*}
for an admissible number $k$ of $\mfrak{g}$.
\medskip

Let $\mfrak{p}=\mfrak{l} \oplus \mfrak{u}$ be a standard parabolic subalgebra of $\mfrak{g}$, where $\mfrak{u}$ is the nilradical of $\mfrak{p}$ and $\mfrak{l}$ is the Levi subalgebra of $\mfrak{p}$, and let $k \in \Q$ be an admissible number for $\mfrak{g}$. We denote by $\Omega_k(\mfrak{p})$ the set of those weights $\lambda \in \widebar{{\rm Pr}}_k \cap \Lambda^+(\mfrak{p})$ for which the generalized Verma module $M^\mfrak{g}_\mfrak{p}(\lambda)$ is simple, i.e.\ $M^\mfrak{g}_\mfrak{p}(\lambda) \simeq L^\mfrak{g}_\mfrak{b}(\lambda)$ as $\mfrak{g}$-modules. For $\lambda \in \Omega_k(\mfrak{p})$, we obtain immediately that $\lambda$ is a regular weight, which by \cite{Jantzen1977} gives us that
$M^\mfrak{g}_\mfrak{p}(\lambda)$ is a simple $\mfrak{g}$-module if and only if $\bra \lambda + \rho, \alpha^\vee \ket \notin \N$ for all $\alpha \in \Delta_+^\mfrak{u}$.
\medskip

\theorem{\label{thm-adm}
Let $\mfrak{p}$ be a standard parabolic subalgebra of $\mfrak{g}$ and let $k \in \Q$ be an admissible number for $\mfrak{g}$. Then the $\mfrak{g}$-module $W^\mfrak{g}_\mfrak{p}(\lambda,\alpha)$ is admissible of level $k$ for $\lambda \in \Omega_k(\mfrak{p})$ and $\alpha \in \Delta^\mfrak{u}_+$.}

\proof{Indeed, the Zhu's algebra $A(\mcal{L}_{k\kappa_0}\!(\mfrak{g}))$ is isomorphic to $U(\mfrak{g})/I_k$, where $I_k$ is a two-sided ideal of $U(\mfrak{g})$. Let us consider the simple highest weight $\mfrak{g}$-module $L^\mfrak{g}_\mfrak{b}(\lambda)$ with highest weight $\lambda \in \mfrak{h}^*$. For $\lambda \in \Omega_k(\mfrak{p})$, the $\mfrak{g}$-module $L^\mfrak{g}_\mfrak{b}(\lambda)$ is admissible of level $k$ by Theorem \ref{thm-hweight}. Hence, we get $I_k \subset J_\lambda$, where $J_\lambda = \Ann_{U(\mfrak{g})}\! L^\mfrak{g}_\mfrak{b}(\lambda)$. Further, we need to show that the $\mfrak{g}$-module $W^\mfrak{g}_\mfrak{p}(\lambda,\alpha)$ for $\alpha \in \Delta^\mfrak{u}_+$ is admissible of level $k$, or equivalently that $I_k \subset \Ann_{U(\mfrak{g})}\!W^\mfrak{g}_\mfrak{p}(\lambda,\alpha)$. Since $W^\mfrak{g}_\mfrak{p}(\lambda,\alpha)$ is obtained from $M^\mfrak{g}_\mfrak{p}(\lambda) \simeq L^\mfrak{g}_\mfrak{b}(\lambda)$ by the twisting functor $T_\alpha$, the statement follows from Corollary \ref{cor-annihilator}. Hence, $W^\mfrak{g}_\mfrak{p}(\lambda,\alpha)$ is a module over $U(\mfrak{g})/I_k$, which implies the required statement.}

\vspace{-2mm}

%%%%%%%%%%%%%%%%%%%%%%%%%%%%%%%%%%%%%%%%%%%%%%%%%%%%%%%%%%%%%%%%%%%%%%%%%%%%%%%%%%%%%%%%%%

\subsection{Richardson orbits and associated varieties}

\begin{figure}[ht]
\centering
{\begin{tikzpicture}
[yscale=1.1,xscale=1.7,vector/.style={circle,draw=white,fill=black,ultra thick, inner sep=0.8mm},vector2/.style={circle,draw=white,fill=white,ultra thick, inner sep=1mm}]
\begin{scope}
  \node (A) at (0,0)  {$\mcal{O}_{{\rm zero}}$};
  \node (B) at (0,1)  {$\mcal{O}_{{\rm min}}$};
  \node (C1) at (0,1.7)  {};
  \node (C) at (0,2)  {$\,\dots$};
  \node (C2) at (0,2.3) {};
  \node (D) at (0,3)  {$\mcal{O}_{{\rm subreg}}$};
  \node (E) at (0,4)  {$\mcal{O}_{{\rm reg}}$};
  \node at (2,0) {$0$};
  \node at (2,1) {$2h^\vee - 2$};
  \node at (2,2) {$\dots$};
  \node at (2,3) {$\dim \mfrak{g} - \rank \mfrak{g} - 2$};
  \node at (2,4) {$\dim \mfrak{g} - \rank \mfrak{g}$};
  \draw [thin, -] (B) -- (C1);
  \draw [thin, -] (D) -- (C2);
  \draw [thin, -] (A) -- (B);
  \draw [thin, -] (D) -- (E);
  \node at (0,4.7) {nilpotent orbit};
  \node at (2,4.7) {dimension};
\end{scope}
\end{tikzpicture}}
\caption{Hasse diagram of nilpotent orbits}
\label{fig:hasse diagram general}
\vspace{-2mm}
\end{figure}
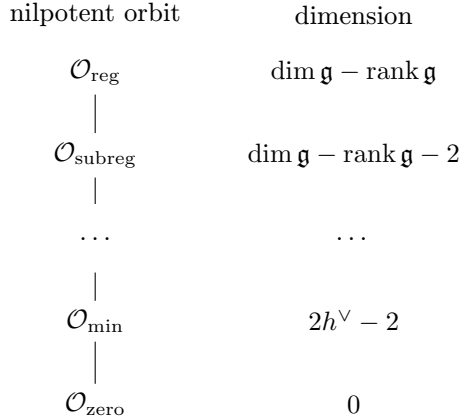

Let $G$ be a complex connected semisimple algebraic group with its Lie algebra $\mfrak{g}$. We denote by $\mcal{N}(\mfrak{g})$ the nilpotent cone of $\mfrak{g}$, i.e.\ the set of nilpotent elements of $\mfrak{g}$. It is an irreducible closed algebraic subvariety of $\mfrak{g}$ and a finite union of $G$-orbits. There is a unique nilpotent orbit of $\mfrak{g}$, denoted by $\mcal{O}_{\rm reg}$ and called the \emph{regular nilpotent orbit} of $\mfrak{g}$, which is a dense open subset of $\mcal{N}(\mfrak{g})$. Next, since $\mfrak{g}$ is simple, there exists a unique nilpotent orbit of $\mfrak{g}$ that is a dense open subset of $\mcal{N}(\mfrak{g}) \setminus \mcal{O}_{\rm reg}$, denoted by $\mcal{O}_{\rm subreg}$ and called the \emph{subregular nilpotent orbit} of $\mfrak{g}$. Besides, there is a unique nonzero nilpotent orbit of $\mfrak{g}$ of minimal dimension, denoted by $\mcal{O}_{\rm min}$ and called the \emph{minimal nilpotent orbit} of $\mfrak{g}$, such that it is contained in the closure of all nonzero nilpotent orbits of $\mfrak{g}$. By $\mcal{O}_{\rm zero}$ we denote the \emph{zero nilpotent orbit} of $\mfrak{g}$.
For the dimension of these distinguished nilpotent orbits of $\mfrak{g}$ see Figure \ref{fig:hasse diagram general}.

Let us consider the PBW filtration on the universal enveloping algebra $U(\mfrak{g})$ of $\mfrak{g}$ and the associated graded algebra $\gr U(\mfrak{g}) \simeq S(\mfrak{g}) \simeq \C[\mfrak{g}^*]$. The \emph{associated variety} $\mcal{V}(I)$ of a left ideal $I$ of $U(\mfrak{g})$ is defined as the zero locus in $\mfrak{g}^*$ of the associated graded ideal $\gr I$ of $S(\mfrak{g})$. Moreover, if $I$ is a two-sided ideal of $\mfrak{g}$, then $I$ and $\gr I$ are invariant under the adjoint action of $G$. Consequently, the associated variety is a union of $G$-orbits of $\mfrak{g}^*$. Obviously, we have
\begin{align}
  \mcal{V}(I) = \Specm (S(\mfrak{g})/\gr I) = \Specm(S(\mfrak{g})/\sqrt{\gr I}),
\end{align}
where $\sqrt{\gr I}$ denotes the radical of $\gr I$. Since the Cartan--Killing form $\kappa_\mfrak{g}$ is a $\mfrak{g}$-invariant symmetric bilinear form on $\mfrak{g}$, it provides a one-to-one correspondence between adjoint orbits of $\mfrak{g}$ and coadjoint orbits of $\mfrak{g}^*$. For an adjoint orbit $\mcal{O}$ of $\mfrak{g}$ we denote by $\mcal{O}^*$ the corresponding coadjoint orbit of $\mfrak{g}^*$. In addition, for a primitive ideal $I$ of $U(\mfrak{g})$ the associated variety $\mcal{V}(I)$ is the closure of $\mcal{O}^*$ for some nilpotent orbit $\mcal{O}$ of $\mfrak{g}$, see \cite{Joseph1985}. If $E$ is a simple $\mfrak{g}$-module and $\mcal{V}(\Ann_{U(\mfrak{g})}\!E)=\widebar{\mcal{O}^*}$ for a nilpotent orbit $\mcal{O}$ of $\mfrak{g}$, we say that $E$ belongs to $\mcal{O}$. For a description of nilpotent orbits of admissible $\mfrak{g}$-modules we refer to \cite{Arakawa2015}.
\medskip

\theorem{\cite{Arakawa2015b}\label{thm:nilpotent orbits}
Let $k \in \Q$ be an admissible number for $\mfrak{g}$ with denominator $q \in \N$. Then there exists a nilpotent orbit $\mcal{O}_q$ of $\mfrak{g}$ such that
\begin{align*}
  \mcal{V}(I_k) = \widebar{{\mcal{O}}_q^*}
\end{align*}
and we have
\begin{align*}
  \widebar{{\mcal{O}}_q} = \begin{cases}
    \{x\in\mfrak{g};\, (\ad x)^{2q}=0\}  & \text{if $(r^\vee,q)=1$}, \\
    \{x\in\mfrak{g};\, \pi_{\theta_s}\!(x)^{2q/r^\vee}=0\}  & \text{if $(r^\vee,q)=r^\vee$},
  \end{cases}
\end{align*}
where $\pi_{\theta_s} \colon \mfrak{g} \rarr \End L^\mfrak{g}_\mfrak{b}(\theta_s)$ is the simple finite-dimensional $\mfrak{g}$-module with highest weight $\theta_s$.}

Let $\mcal{O}$ be a nilpotent orbit of $\mfrak{g}$ and let $k \in \Q$ be an admissible number for $\mfrak{g}$ with denominator $q \in \N$. We define the subset
\begin{align}
  \widebar{{\rm Pr}}_k^\mcal{O} =\{\lambda\in \widebar{{\rm Pr}}_k;\, \mcal{V}(J_\lambda) = \widebar{\mcal{O}^*}\}
\end{align}
of $\mfrak{h}^*$. Then a simple $\mfrak{g}$-module $E$ in the nilpotent orbit $\mcal{O}$ is admissible of level $k$ if and only if $\Ann_{U(\mfrak{g})}\!E=J_\lambda$ for some $\lambda \in \widebar{{\rm Pr}}_k^\mcal{O}$. Further, as $I_k \subset J_\lambda$ for $\lambda \in \smash{\widebar{{\rm Pr}}_k}$, we have $\mcal{V}(J_\lambda) \subset \mcal{V}(I_k) = \smash{\widebar{{\mcal{O}}_q^*}}$ by Theorem \ref{thm:nilpotent orbits}, which gives us
\begin{align}
  \widebar{{\rm Pr}}_k = \bigsqcup_{\mcal{O} \subset \widebar{{\mcal{O}}_q}}  \widebar{{\rm Pr}}_k^\mcal{O}. \label{eq:Pr_k decomposition}
\end{align}
Therefore, we need to describe the subset $\widebar{{\rm Pr}}_k^\mcal{O}$ of $\widebar{{\rm Pr}}_k$ for a nilpotent orbit $\mcal{O}$ of $\mfrak{g}$.
\medskip

For $x \in \mfrak{g}$, we denote by $\mfrak{g}^x$ the centralizer of $x$ in $\mfrak{g}$. Furthermore, for a subset $X$ of $\mfrak{g}$ we define the set
\begin{align*}
  X^{\rm reg} = \{x \in X;\, \dim \mfrak{g}^x = \min\nolimits_{y \in X} \dim \mfrak{g}^y \}
\end{align*}
and call it the \emph{set of regular elements} in $X$.
\medskip

\theorem{\cite{Borho-Brylinski1982}\label{thm:annihilator generalized Verma}
Let $\mfrak{p}$ be a standard parabolic subalgebra of $\mfrak{g}$. Then for $\lambda \in \Lambda^+(\mfrak{p})$ the associated variety $\mcal{V}(\Ann_{U(\mfrak{g})}\!M^\mfrak{g}_\mfrak{p}(\lambda))$ is the closure of $\mcal{O}_\mfrak{p}^*$, where
\begin{align*}
  \mcal{O}_\mfrak{p} = (G.\mfrak{p}^\perp)^{\rm reg}
\end{align*}
and $\mfrak{p}^\perp$ is the orthogonal complement of $\mfrak{p}$ with respect to the Cartan-Killing form. In particular, we have
\begin{align*}
  \widebar{{\mcal{O}}_{\mfrak{p}}} = G.\mfrak{p}^\perp
\end{align*}
and the associated variety $\mcal{V}(\Ann_{U(\mfrak{g})}\!M^\mfrak{g}_\mfrak{p}(\lambda))$ is irreducible. The orbit $  \mcal{O}_\mfrak{p}$ is the {\it Richardson orbit} determined by $\mfrak{p}$.}

Now, we find the standard parabolic subalgebras for which the corresponding Richardson orbits  are the distinguished nilpotent orbits $\mcal{O}_{\rm zero}$, $\mcal{O}_{\rm min}$, $\mcal{O}_{\rm subreg}$ and $\mcal{O}_{\rm reg}$. Let us recall that while the regular orbit $\mcal{O}_{\rm reg}$ and the subregular orbit $\mcal{O}_{\rm subreg}$ are Richardson orbits for any $\mfrak{g}$, the minimal orbit $\mcal{O}_{\rm min}$ is a Richardson orbit only for $\mfrak{sl}_n$, $n \geq 2$. Since all these nilpotent orbits are uniquely determined by their dimensions and since $\dim \mcal{O}_\mfrak{p} = 2\dim \mfrak{u}$ for a parabolic subalgebra $\mfrak{p}=\mfrak{l}\oplus \mfrak{u}$ of $\mfrak{g}$ with the nilradical $\mfrak{u}$, we obtain easily by comparison of the dimensions that
\begin{enumerate}[topsep=3pt,itemsep=0pt]
  \item[i)] $\mcal{O}_{\rm zero} = \mcal{O}_\mfrak{g}$,
  \item[ii)] $\mcal{O}_{\rm min} = \mcal{O}_{\mfrak{p}_\alpha^{\rm max}}$ for $\alpha \in \{\alpha_1,\alpha_{n-1}\}$ and $\mfrak{g} = \mfrak{sl}_n$, $n \geq 2$,
  \item[iii)] $\mcal{O}_{\rm subreg} = \mcal{O}_{\mfrak{p}_\alpha^{\rm min}}$ for $\alpha \in \Pi$,
  \item[iv)] $\mcal{O}_{\rm reg} = \mcal{O}_\mfrak{b}$,
\end{enumerate}
where $\mfrak{p}_\alpha^{\rm max}$ for $\alpha \in \Pi$ is the parabolic subalgebra of $\mfrak{g}$ associated to the subset $\Sigma = \Pi \setminus \{\alpha\}$, and $\mfrak{p}_\alpha^{\rm min}$ for $\alpha \in \Pi$ is the parabolic subalgebra of $\mfrak{g}$ associated to the subset $\Sigma = \{\alpha\}$.
\medskip

\lemma{\label{lem:intersection}
Let $\mfrak{p}_1$, $\mfrak{p}_2$ be standard parabolic subalgebras of $\mfrak{g}$ and let $k \in \Q$ be an admissible number for $\mfrak{g}$. Then we have  $\Omega_k(\mfrak{p}_1) \cap \Omega_k(\mfrak{p}_2) \neq \emptyset$ if and only if $\mfrak{p}_1 = \mfrak{p}_2$.}

\proof{Let $\mfrak{p}$ be a standard parabolic subalgebra of $\mfrak{g}$. If $\lambda \in \Omega_k(\mfrak{p})$ we have $\bra \lambda + \rho, \alpha^\vee \ket \notin \N$ for all $\alpha \in \Delta_+^\mfrak{u}$ and $\bra \lambda+\rho, \alpha^\vee \ket \in \N$ for all $\alpha \in \Delta_+^\mfrak{l}$, which immediately implies that $\Omega_k(\mfrak{p}_1) \cap \Omega_k(\mfrak{p}_2) \neq \emptyset$ if and only if $\mfrak{p}_1 = \mfrak{p}_2$.}

Let $\mcal{O}$ be a Richardson orbit of $\mfrak{g}$ and let $k \in \Q$ be an admissible number for $\mfrak{g}$. Further, let us consider a standard parabolic subalgebra $\mfrak{p}$ of $\mfrak{g}$ satisfying $\mcal{O}_\mfrak{p} = \mcal{O}$. Then for $\lambda \in \Omega_k(\mfrak{p})$ we obtain $\mcal{V}(\Ann_{U(\mfrak{g})}\!L^\mfrak{g}_\mfrak{b}(\lambda)) = \widebar{{\mcal{O}}^*}$ by Theorem \ref{thm:annihilator generalized Verma}, since $M^\mfrak{g}_\mfrak{p}(\lambda) \simeq L^\mfrak{g}_\mfrak{b}(\lambda)$ for any $\lambda \in \Omega_k(\mfrak{p})$, which gives us
\begin{align*}
\Omega_k(\mfrak{p}) \subset \widebar{{\rm Pr}}_k^\mcal{O}.
\end{align*}
Hence, by Lemma \ref{lem:intersection} we get
\begin{align}
  \bigsqcup_{\substack{\Sigma \subset \Pi\\ \mcal{O}_{\mfrak{p}_\Sigma} = \mcal{O}}} \Omega_k(\mfrak{p}_\Sigma) \subset \widebar{{\rm Pr}}_k^\mcal{O},
\end{align}
where $\mfrak{p}_\Sigma$ is the standard parabolic subalgebra of $\mfrak{g}$ associated to a subset $\Sigma$ of $\Pi$.
\medskip

\proposition{\label{prop-orbit-borel}
Let $k \in \Q$ be an admissible number for $\mfrak{g}$. Then we have $\Omega_k(\mfrak{b}) = \widebar{{\rm Pr}}_k^{\smash{\mcal{O}_{\rm reg}}}$.}

\proof{Since the regular nilpotent orbit $\mcal{O}_{\rm reg}$ of $\mfrak{g}$ is a Richardson orbit and $\mcal{O}_{\rm reg} = \mcal{O}_\mfrak{b}$, we obtain $\Omega_k(\mfrak{b}) \subset \widebar{{\rm Pr}}_k^{\smash{\mcal{O}_{\rm reg}}}$. On the other hand, for $\lambda \in \widebar{{\rm Pr}}_k^{\smash{\mcal{O}_{\rm reg}}}$ we have $\mcal{V}(J_\lambda) = \smash{\widebar{\mcal{O}_{\rm reg}^*}}$ by definition, where $J_\lambda = \Ann_{U(\mfrak{g})}\!L^\mfrak{g}_\mfrak{b}(\lambda)$, and $\mcal{V}(I_\lambda) = \smash{\widebar{\mcal{O}_{\rm reg}^*}}$ by Theorem \ref{thm:annihilator generalized Verma}, where $I_\lambda = \Ann_{U(\mfrak{g})}\!M^\mfrak{g}_\mfrak{b}(\lambda)$. Hence, we get $\sqrt{\gr J_\lambda} = \sqrt{\gr I_\lambda}$. Since $\gr I_\lambda$ is a prime ideal of $S(\mfrak{g})$ by \cite[Theorem 5.6]{Borho-Brylinski1982}, we have $\gr J_\lambda \subset \gr I_\lambda$ which gives us $J_\lambda \subset I_\lambda$. Moreover, the primitive ideal $J_\lambda$ is the unique maximal two-sided ideal of $U(\mfrak{g})$ containing $I_\lambda$ by \cite[Proposition 2.4]{Arakawa-Futorny-Ramirez2017}. Therefore, we get $I_\lambda = J_\lambda$ for $\lambda \in \widebar{{\rm Pr}}_k^{\smash{\mcal{O}_{\rm reg}}}$. As $\lambda$ is regular dominant, we have $M^\mfrak{g}_\mfrak{b}(\lambda) \simeq M^\mfrak{g}_\mfrak{b}(\lambda)/I_\lambda M^\mfrak{g}_\mfrak{b}(\lambda)$ and $L^\mfrak{g}_\mfrak{b}(\lambda) \simeq M^\mfrak{g}_\mfrak{b}(\lambda)/J_\lambda L^\mfrak{g}_\mfrak{b}(\lambda)$, which implies $M^\mfrak{g}_\mfrak{b}(\lambda) \simeq L^\mfrak{g}_\mfrak{b}(\lambda)$ and thus $\lambda \in \Omega_k(\mfrak{b})$.}

Besides, from the decomposition \eqref{eq:Pr_k decomposition} is follows that the set $\widebar{{\rm Pr}}_k^{\smash{\mcal{O}_{\rm reg}}}$ is non-empty if and only if $\mcal{O}_{\rm reg} = \mcal{O}_q$, where $q$ is the denominator of $k$, or equivalently if and only if
\begin{align*}
    q \geq \begin{cases}
             h & \text{if $(r^\vee,q)=1$},\\
             {}^L h^\vee r^\vee & \text{if $(r^\vee,q)=r^\vee$},
           \end{cases}
\end{align*}
where $h$ is the Coxeter number of $\mfrak{g}$ and ${}^L h^\vee$ is the dual Coxeter number of the Langlands dual Lie algebra ${}^L\mfrak{g}$ of $\mfrak{g}$ (see \cite{Arakawa-Futorny-Ramirez2017}).
\medskip

\proposition{Let $k \in \Q$ be an admissible number for $\mfrak{g}$. Then we have $\Omega_k(\mfrak{g}) = \widebar{{\rm Pr}}_k^{\smash{\mcal{O}_{\rm zero}}}$.}

\proof{For the zero nilpotent orbit $\mcal{O}_{\rm zero}$ of $\mfrak{g}$, we have $\mcal{O}_{\rm zero} = \mcal{O}_\mfrak{g}$, which gives us $\Omega_k(\mfrak{g}) \subset \widebar{{\rm Pr}}_k^{\smash{\mcal{O}_{\rm zero}}}$. On the other hand, for $\lambda \in \widebar{{\rm Pr}}_k^{\smash{\mcal{O}_{\rm zero}}}$ we have $\mcal{V}(J_\lambda) = \smash{\widebar{\mcal{O}_{\rm zero}^*}}$ by definition, where $J_\lambda = \Ann_{U(\mfrak{g})}\!L^\mfrak{g}_\mfrak{b}(\lambda)$, which implies $\smash{\sqrt{\mfrak{(g)}}} = \sqrt{\gr J_\lambda}$. Hence, for $\alpha \in \Delta_+$ there exists a positive integer $n_\alpha \in \N$ such that $f_\alpha^{n_\alpha} \in \gr J_\lambda$. Therefore, we get immediately that the $\mfrak{g}$-module $L^\mfrak{g}_\mfrak{b}(\lambda)$ is finite dimensional implying $\lambda \in \Omega_k(\mfrak{g})$.}

\vspace{-2mm}

%%%%%%%%%%%%%%%%%%%%%%%%%%%%%%%%%%%%%%%%%%%%%%%%%%%%%%%%%%%%%%%%%%%%%%%%%%%%%%%%%%%%%%%%%%

\subsection{Admissible representations for $\mfrak{sl}_n$}

In this subsection we will consider the simple Lie algebra $\mfrak{g}= \mfrak{sl}_n$ for $n \geq 2$. Let us note that all nilpotent orbits of $\mfrak{g}$ are Richardson orbits.
\medskip

Let $k \in \Q$ be an admissible number for $\mfrak{g}$ with denominator $q \in \N$. Then the nilpotent orbit $\mcal{O}_q$ from Theorem \ref{thm:nilpotent orbits} is given by
\begin{align*}
  \mcal{O}_q = \mcal{O}_{\lambda_q},
\end{align*}
where $\lambda_q$ is the partition of $n$ defined through $\lambda_q = [q^r,s]$, where $n=qr+s$ with $r,s\in \N_0$ and $0 \leq s \leq q-1$. Hence, we immediately get that
\begin{align*}
  \mcal{O}_q = \begin{cases}
    \mcal{O}_{\rm reg} &  \text{if $q \geq n$}, \\
    \mcal{O}_{\rm subreg} & \text{if $q=n-1$}, \\
    \mcal{O}_{\rm zero} & \text{if $q=1$}
  \end{cases}
\end{align*}
for $n \geq 2$. For details about nilpotent orbits in semisimple Lie algebras see \cite{Collingwood-McGovern1993-book}.

\begin{figure}[ht]
\centering
\subcaptionbox{$\mfrak{g}=\mfrak{sl}_2$\label{fig:sl2}}[0.25\textwidth]
{\begin{tikzpicture}
[yscale=1.2, xscale=1]
\begin{scope}
  \node (A) at (0,0)  {$\mcal{O}_{[1^2]}$};
  \node (B) at (0,1)  {$\mcal{O}_{[2]}$};
  \node at (1,0) {$0$};
  \node at (1,1) {$2$};
  \node at (2,0) {\dynkin[x/.style={thin}]{A}{o}};
  \node at (2,1) {\dynkin[x/.style={thin}]{A}{x}};
 \draw [thin, -] (A) -- (B);
\end{scope}
\end{tikzpicture}}
\hfill
\subcaptionbox{$\mfrak{g}=\mfrak{sl}_3$\label{fig:sl3}}[0.25\textwidth]
{\begin{tikzpicture}
[yscale=1.2, xscale=1]
\begin{scope}
  \node (A) at (0,0)  {$\mcal{O}_{[1^3]}$};
  \node (B) at (0,1)  {$\mcal{O}_{[2,1]}$};
  \node (C) at (0,2)  {$\mcal{O}_{[3]}$};
  \node at (1,0) {$0$};
  \node at (1,1) {$4$};
  \node at (1,2) {$6$};
 \node at (2,0) {\dynkin[x/.style={thin}]{A}{oo}};
  \node at (2,0.9) {\dynkin[x/.style={thin}]{A}{ox}};
 \node at (2,1.1) {\dynkin[x/.style={thin}]{A}{xo}};
  \node at (2,2) {\dynkin[x/.style={thin}]{A}{xx}};
  \node at (0,4.3) {$\mcal{O}$};
  \node at (1,4.3) {$\dim$};
  \node at (2,4.3) {$\mfrak{p}_\Sigma$};
  \draw [thin, -] (A) -- (B);
  \draw [thin, -] (B) -- (C);
\end{scope}
\end{tikzpicture}}
\hfill
\subcaptionbox{$\mfrak{g}=\mfrak{sl}_4$\label{fig:sl4}}[0.25\textwidth]
{\begin{tikzpicture}
[yscale=1.2, xscale=1]
\begin{scope}
  \node (A) at (0,0)  {$\mcal{O}_{[1^4]}$};
  \node (B) at (0,1)  {$\mcal{O}_{[2,1^2]}$};
 \node (C) at (0,2)  {$\mcal{O}_{[2^2]}$};
  \node (D) at (0,3)  {$\mcal{O}_{[3,1]}$};
  \node (E) at (0,4)  {$\mcal{O}_{[4]}$};
  \node at (1,0) {$0$};
  \node at (1,1) {$6$};
  \node at (1,2) {$8$};
  \node at (1,3) {$10$};
  \node at (1,4) {$12$};
  \node at (2,0) {\dynkin[x/.style={thin}]{A}{ooo}};
  \node at (2,0.9) {\dynkin[x/.style={thin}]{A}{oox}};
  \node at (2,1.1) {\dynkin[x/.style={thin}]{A}{xoo}};
  \node at (2,2) {\dynkin[x/.style={thin}]{A}{oxo}};
  \node at (2,2.8) {\dynkin[x/.style={thin}]{A}{oxx}};
  \node at (2,3) {\dynkin[x/.style={thin}]{A}{xox}};
  \node at (2,3.2) {\dynkin[x/.style={thin}]{A}{xxo}};
 \node at (2,4) {\dynkin[x/.style={thin}]{A}{xxx}};
  \draw [thin, -] (A) -- (B);
  \draw [thin, -] (B) -- (C);
  \draw [thin, -] (C) -- (D);
  \draw [thin, -] (D) -- (E);
\end{scope}
\end{tikzpicture}}
\hspace{\the\parindent}
\caption{Richardson orbits for $\mfrak{sl}_n$}
\label{fig:hasse diagram sl}
\vspace{-2mm}
\end{figure}
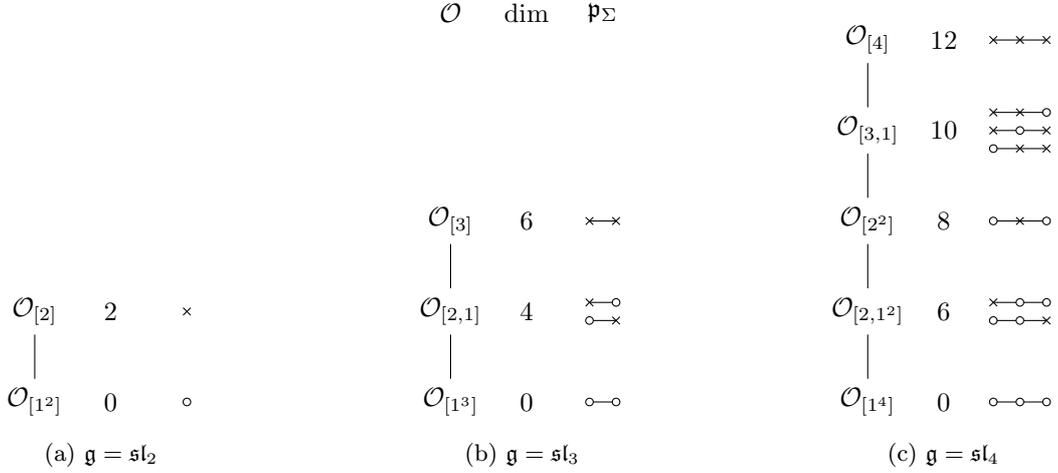

Let us introduce the subset
\begin{align*}
  {\rm Pr}_{k,\Z} = \{\lambda \in {\rm Pr}_k;\, \bra \lambda, \alpha^\vee \ket \in \Z \text{ for $\alpha \in \Pi$}\}
\end{align*}
of admissible weights of level $k$. Then by \cite{Kac-Wakimoto1989} we have
\begin{align*}
  {\rm Pr}_{k,\Z} = \{\lambda \in \widehat{\mfrak{h}}^*;\, \lambda(c)=k,\, \bra \lambda, \alpha^\vee \ket \in \N_0 \text{ for $\alpha \in \Pi$},\, \bra \lambda, \theta^\vee \ket \leq p-h^\vee \},
\end{align*}
where $p = (k+h^\vee)q$. Further, from \cite{Kac-Wakimoto1989} we obtain
\begin{align*}
  {\rm Pr}_k = \bigcup_{\substack{y \in \widetilde{W} \\ y(\widehat{\Delta}(k\Lambda_0)_+) \subset \widehat{\Delta}^{\rm re}_+}} {\rm Pr}_{k,y}, \quad {\rm Pr}_{k,y} = \{y\cdot \lambda;\, \lambda \in {\rm Pr}_{k,\Z}\},
\end{align*}
where the extended affine Weyl group $\smash{\widetilde{W}}$ of $\mfrak{g}$ is defined as
\begin{align*}
  \smash{\widetilde{W}} = W \ltimes P^\vee.
\end{align*}
For $\mu \in P^\vee$, we denote by $t_\mu$ the corresponding element of $\smash{\widetilde{W}}$.
The action of $t_\mu$ on $\smash{\widetilde{\mfrak{h}}^*}$ for $\mu \in P^\vee$ is given by
\begin{align*}
  t_\mu(\gamma) = \gamma + (\gamma,\delta)\mu - \bigg({(\alpha,\alpha) \over 2}\,(\gamma,\delta) + (\gamma,\mu)\!\!\bigg)\delta
\end{align*}
for $\gamma \in \smash{\widetilde{\mfrak{h}}^*}$. Moreover, for $y,y' \in \smash{\widetilde{W}}$ satisfying $y(\widehat{\Delta}(k\Lambda_0)_+) \subset \widehat{\Delta}_+^{\rm re}$, $y'(\widehat{\Delta}(k\Lambda_0)_+) \subset \widehat{\Delta}_+^{\rm re}$ we have
\begin{align}
  {\rm Pr}_{k,y} \cap {\rm Pr}_{k,y'} \neq \emptyset \quad \Longleftrightarrow \quad  {\rm Pr}_{k,y} = {\rm Pr}_{k,y'} \quad \Longleftrightarrow \quad y' = yt_{q\omega_j} w_j \label{eq:Pr equality}
\end{align}
for some $j \in \{1,2,\dots,n-1\}$, where $w_j$ is the unique element of the Weyl group $W$ of $\mfrak{g}$ which preserves the set $\{\alpha_1,\alpha_2,\dots,\alpha_{n-1},-\theta\}$ and $w_j(-\theta) = \alpha_j$. The set of simple roots of $\mfrak{g}$ is $\Pi=\{\alpha_1, \alpha_2, \dots, \alpha_{n-1}\}$. Let us recall that by \cite[Proposition 2.8]{Arakawa-Futorny-Ramirez2017}
we have
\begin{align*}
  [\widebar{{\rm Pr}}_k] = \bigcup_{\substack{\eta \in P_+^\vee \\ (\eta,\theta) \leq q-1}} [\widebar{{\rm Pr}}_{k,t_{-\eta}}],
\end{align*}
where $P_+^\vee$ is the set of dominant coweights of $\mfrak{g}$.
\medskip

\lemma{\label{lem:distinguished element}
If $\widebar{y} \in W$ is not the unit element, then there exists $j \in \{0,1,\dots,n-1\}$ such that $\widebar{y}w_j(\theta) \in \Delta_+$.}

\proof{For a proof, see \cite[Lemma 2.7]{Arakawa-Futorny-Ramirez2017}.}

\theorem{\label{thm-omega-sln}
Let $k \in \Q$ be an admissible number for $\mfrak{g}$ with denominator $q \in \N$. Then we have
\begin{align*}
  \Omega_k(\mfrak{p}_\Sigma) = \bigcup_{\substack{\widebar{y} \in W,\, \eta \in P^\vee_+, \\  (\eta,\theta) \leq q-1,\, \widebar{y}(\theta) \in \Delta_+,\ \widebar{y}(\Delta^\eta_0) = \Delta_\Sigma, \\ \Delta_0^\eta \cap \Delta_+ \subset \widebar{y}^{-1}(\Delta_+) \cap \Delta_+}} \widebar{{\rm Pr}}_{k,\widebar{y}t_{-\eta}}
\end{align*}
for $\Sigma \subset \Pi$, where $\Delta_\Sigma$ is the subroot system of $\Delta$ generated by $\Sigma$ and $\Delta^\eta_0 = \{\alpha \in \Delta;\, (\eta,\alpha)=0\}$ for $\eta \in P^\vee$.}

\proof{Let $y=\widebar{y}t_{-\eta}$ with $\widebar{y} \in W$ and $\eta \in P^\vee$. Then the condition $y(\smash{\widehat{\Delta}}(k\Lambda_0)_+) \subset \smash{\widehat{\Delta}}^{\rm re}_+$ is equivalent to
\begin{align*}
0 \leq (\eta,\alpha) \leq q-1 \text{ if $\widebar{y}(\alpha) \in \Delta_+$} \qquad \text{and}  \qquad 1 \leq (\eta,\alpha) \leq q \text{ if $\widebar{y}(\alpha) \in \Delta_-$}
\end{align*}
for all $\alpha \in \Delta_+$. In particular, for $\eta \in P^\vee$ and $\widebar{y} \in W$ such that $\widebar{y}(\theta) \in \Delta_+$ we obtain easily that $\widebar{y}t_{-\eta}(\smash{\widehat{\Delta}}(k\Lambda_0)_+) \subset \smash{\widehat{\Delta}}^{\rm re}_+$ if and only if $\eta \in P^\vee_+$, $(\eta,\theta) \leq q-1$ and $\Delta_0^\eta \cap \Delta_+ \subset \widebar{y}^{-1}(\Delta_+) \cap \Delta_+$.

For $\lambda \in \mfrak{h}^*$, we define the subset $\Delta(\lambda) = \{\alpha\in\Delta;\, \bra \lambda+\rho, \alpha^\vee \ket \in \Z\}$ of $\Delta$. Then we immediately get $\Delta(\lambda) = \smash{\widehat{\Delta}}(\lambda+k\Lambda_0) \cap \Delta$. Further, let us assume that $\lambda \in \widebar{{\rm Pr}}_{k,\widebar{y}t_{-\eta}}$ with $\widebar{y} \in W$, $\eta \in P^\vee_+$, $(\eta,\theta) \leq q-1$, $\widebar{y}(\theta) \in \Delta_+$, $\Delta_0^\eta \cap \Delta_+ \subset \widebar{y}^{-1}(\Delta_+) \cap \Delta_+$ and $\widebar{y}(\Delta^\eta_0) = \Delta_\Sigma$ for a fixed subset $\Sigma$ of $\Pi$. Then we may write
\begin{align*}
  \Delta(\lambda) &= \smash{\widehat{\Delta}}(\lambda+k\Lambda_0) \cap \Delta = \widebar{y}t_{-\eta}(\smash{\widehat{\Delta}}(k\Lambda_0)) \cap \Delta = \{\widebar{y}(\alpha) + (mq+ (\eta,\alpha))\delta;\, \alpha \in \Delta,\, m\in \Z\} \cap \Delta \\
  &= \{\widebar{y}(\alpha);\, \alpha \in \Delta,\, (\eta,\alpha) = 0\} = \widebar{y}(\Delta^\eta_0) =\Delta_\Sigma.
\end{align*}
Hence, by \cite{Jantzen1977} we have that $M^\mfrak{g}_{\mfrak{p}_\Sigma}(\lambda) \simeq L^\mfrak{g}_\mfrak{b}(\lambda)$ since $\lambda$ is regular dominant, i.e. $\bra \lambda + \rho,\alpha^\vee \ket \notin -\N_0$ for $\alpha \in \Delta_+$, which gives us $\widebar{{\rm Pr}}_{k,\widebar{y}t_{-\eta}} \subset \Omega_k(\mfrak{p}_{\Sigma})$.

On the other hand, if $\lambda \in \Omega_k(\mfrak{p}_\Sigma)$, then it is easy to see that $\Delta(\lambda) = \Delta_\Sigma$. Since we also have $\lambda \in \widebar{{\rm Pr}}_k$, there exists $y \in \smash{\widetilde{W}}$ such that $\smash{\widehat{\Delta}}(\lambda+k\Lambda_0) = y(\smash{\widehat{\Delta}}(k\Lambda_0))$, $y(\smash{\widehat{\Delta}}(k\Lambda_0)_+) \subset \smash{\widehat{\Delta}}_+^{\rm re}$ and $\lambda \in \widebar{{\rm Pr}}_{k,y}$. Let $y=\widebar{y}t_{-\eta}$ with $\widebar{y} \in W$ and $\eta \in P^\vee$. We may assume that $\widebar{y}(\theta) \in \Delta_+$. Indeed, we have
\begin{align*}
  \widebar{{\rm Pr}}_{k,\widebar{y}t_{-\eta}} = \widebar{{\rm Pr}}_{k,\widebar{y}t_{-\eta}t_{q\omega_j}w_j} = \widebar{{\rm Pr}}_{k,\widebar{y}w_jw_j^{-1}t_{q\omega_j-\eta}w_j} = \widebar{{\rm Pr}}_{k, \widebar{y}w_j t_{-\eta_j}}
\end{align*}
by \eqref{eq:Pr equality} for $j \in \{0,1,\dots,n-1\}$, where $\eta_j = w_j^{-1}(\eta-q\omega_j)$. Further, if $\widebar{y}(\theta) \in \Delta_-$, then there exists $j \in \{0,1,\dots,n-1\}$ satisfying $\widebar{y}w_j(\theta) \in \Delta_+$ by Lemma \ref{lem:distinguished element}.

Hence, by using the assumption $\widebar{y}(\theta) \in \Delta_+$, we get that $\lambda \in \widebar{{\rm Pr}}_{k,\widebar{y}t_{-\eta}}$ with $\eta \in P_+^\vee$, $(\eta,\theta) \leq q-1$, $\Delta_0^\eta \cap \Delta_+ \subset \widebar{y}^{-1}(\Delta_+) \cap \Delta_+$ and $\widebar{y}(\Delta_0^\eta) = \Delta(\lambda) = \Delta_\Sigma$.}

\corollary{\label{cor:Omega sets}
Let $k \in \Q$ be an admissible number for $\mfrak{g}$ with denominator $q \in \N$. Then we have
\begin{enumerate}[topsep=3pt,itemsep=0pt]
  \item[i)] $\Omega_k(\mfrak{b}) = \bigcup_{\substack{\widebar{y} \in W,\, \eta \in P^\vee_+, \\  (\eta,\theta) \leq q-1,\, \widebar{y}(\theta) \in \Delta_+,\, \Delta_0^\eta = \emptyset}} \widebar{{\rm Pr}}_{k,\widebar{y}t_{-\eta}}$;
  \item[ii)] $\Omega_k(\mfrak{g}) = \widebar{{\rm Pr}}_{k,e}$.
\end{enumerate}
}

Let us note that as an immediate consequence of Corollary \ref{cor:Omega sets} we have that $\Omega_k(\mfrak{b}) \neq \emptyset$ if and only if $q \geq n$.

%%%%%%%%%%%%%%%%%%%%%%%%%%%%%%%%%%%%%%%%%%%%%%%%%%%%%%%%%%%%%%%%%%%%%%%%%%%%%%%%%%%%%%%%%%
%%%%%%%%%%%%%%%%%%%%%%%%%%%%%%%%%%%%%%%%%%%%%%%%%%%%%%%%%%%%%%%%%%%%%%%%%%%%%%%%%%%%%%%%%%

\section*{Acknowledgments}

V.\,F.\ is supported in part by the CNPq (304467/2017-0) and by the Fapesp (2018/23690-6). The authors are gratefully acknowledge the hospitality and excellent working conditions of the International Center for Mathematics of SUSTech and the University of Sichuan (Chengdu, China) where part of this work was done. Both authors are grateful to Edward Frenkel for interest and fruitful discussions.

%%%%%%%%%%%%%%%%%%%%%%%%%%%%%%%%%%%%%%%%%%%%%%%%%%%%%%%%%%%%%%%%%%%%%%%%%%%%%%%%%%%%%%%%%%
%%%%%%%%%%%%%%%%%%%%%%%%%%%%%%%%%%%%%%%%%%%%%%%%%%%%%%%%%%%%%%%%%%%%%%%%%%%%%%%%%%%%%%%%%%

%\bibliographystyle{amsalpha}
%%\bibliographystyle{amsplain}
%\bibliography{reference}

\providecommand{\bysame}{\leavevmode\hbox to3em{\hrulefill}\thinspace}
\providecommand{\MR}{\relax\ifhmode\unskip\space\fi MR }
\providecommand{\MRhref}[2]{%
  \href{http://www.ams.org/mathscinet-getitem?mr=#1}{#2}
}
\providecommand{\href}[2]{#2}

\end{document}